\documentclass[11pt]{arxiv-article}

\numberwithin{equation}{section}

\DeclareFontFamily{U}{mathx}{}
\DeclareFontShape{U}{mathx}{m}{n}{
  <5> <6> <7> <8> <9> <10>
  <10.95> <12> <14.4> <17.28> <20.74> <24.88>
  mathx10
}{}
\DeclareSymbolFont{mathx}{U}{mathx}{m}{n}
\DeclareMathAccent{\widebar}{0}{mathx}{"73}

\newcommand{\prob}[1]{\mathfrak{#1}\xspace}
\newcommand{\order}[1]{\mathcal{O}\left(#1\right)}
\DeclareMathOperator{\expec}{\mathbb{E}}
\DeclarePairedDelimiterX\infdivx[2]{(}{)}{#1\;\delimsize\|\;#2}
\newcommand{\infdiv}[2]{\operatorname{D}_{\hspace{0.1mm}}\infdivx{#1}{#2}}

\DeclareMathOperator*{\rank}{rank}
\DeclareMathOperator*{\argmin}{argmin}

\definecolor{bostonuniversityred}{rgb}{0.8, 0.0, 0.0}

\begin{document}

\onehalfspace

\thispagestyle{empty}

\begin{flushright}
\phantom{Version: \today}
\\
\end{flushright}
\vskip .2 cm
\subsection*{}
\begin{center}
{\Large {\bf \textcolor{bostonuniversityred}{
The Boltzmann structure of sampling:\\[0.8ex]
Intrinsic $p$-value and its emergent closed-form expression
}
} }
\\[0pt]

\bigskip
\bigskip 
{\large
{\bf Orestis Loukas}\footnote{E-mail: orestis.loukas@uni-marburg.de} 
}
\\[2.5pt]
{\it Institute for Medical Bioinformatics and Biostatistics\\
Philipps-Universität Marburg\\
Hans-Meerwein-Straße 6, 35032 Germany}

\bigskip
\end{center}

\begin{abstract}
We consider observables $X$ whose realizations in a sampled dataset are restricted, for example by measurement resolution, to a finite set of distinguishable categories within their possibly infinite theoretical domain. When the probabilities of observable categories are known, we study the coarse-grained probability mass assigned to families of possible datasets generated through a forward sampling process.
The combinatorial construction naturally induces an intrinsically discrete $p$-value, 
defined directly from the sampling process itself rather than through additional probabilistic structure on observables. 

Concretely, specifying the sample means of $d+k$ arbitrary functions $g_\alpha(X)$ determines a linear family of datasets characterized by $\widebar{g_\alpha(X)}$. 
Estimating the probability mass of such linear families reduces to a weighted sum over all integer lattice points contained within the associated polyhedron.
To circumvent intractable large-$N$ combinatorics, we derive via saddle-point techniques a density that approximates these probability masses in the continuum limit of forward sampling processes in the universality class of multinomial sampling. 

As a demonstration, we consider conditional sampling, in which $d$ structural means $\widebar{g_\alpha(X)}$ are held fixed, while the remaining $k$ means $\widebar{g_i(X)}$ vary over admissible datasets. The information geometry emerging from the saddle-point density, together with the spherical symmetry arising at large $N$ from the intrinsic $p$-value construction, enables efficient computation of the $p$-value in the Laplace approximation via the 
$\chi^2_k$ distribution. 
In particular, the resulting $\chi^2_k$ statistic is given in a semi-analytic closed form by $2N$ times the Kullback-Leibler divergence between the information projections associated with the linear families defined by the
observed $\widebar{g_i(X)}$ and structural $\widebar{g_\alpha(X)}$ means. Crucially, these information projections can be computed efficiently via standard numerical routines that converge to the desired accuracy for sufficiently well-behaved sample means.
\end{abstract}


\section{Introduction}

Observables in a scientific theory typically enter through the expectation values of functions of the observables. Given a relevant dataset, these expectations are naturally estimated by the corresponding sample means, yielding generalized moments that range from a simple prevalence or conditional mean to complex mean-field relations and macroscopic empirical laws.

In this work, we pursue an alternative approach to the study of generalized moments defined by arbitrary functions of observables.
Rather than dealing with the traditional notions of random variables and their underlying distributions~\cite{Kolmogorov1950}, we choose to characterize the  uncertainty associated with these generalized moments through elementary combinatorial arguments.
From this perspective, uncertainty is quantified through the combinatorics of observed entities that can be distinguished only up to the categories induced by the  observables under investigation.

As has been known for more than a century, the combinatorics imposed by a broad class of sampling processes effectively coincide with those of multinomial sampling at large enough sample sizes, so that any competing physical scale is suppressed. The resulting universality implies that these processes can, in the appropriate regime, be regarded as sampling with replacement from a finite urn or equivalently without replacement from an infinite urn. At the same time, the urn models serving as a Gedankenexperiment for these combinatorial arguments naturally accommodate more general sampling mechanisms involving additional physical scales. They thus offer a systematic way to account for additional known scales by perturbatively extending the classical urn paradigm beyond the effective limit of sampling with replacement.

Our starting point is the probability mass associated with families of datasets supporting prescribed values of the target generalized moments. From a statistical-mechanical perspective, this object is naturally interpreted as a collective quantity, analogous to how macroscopic mass arises as the coarse-grained sum of many microscopic constituents (such as atoms).
By contemplating  the  probability masses of dataset families associated with prescribed values of generalized moments, the familiar combinatorics of urn models reveal an unexpected structure: they come equipped with a natural notion of a $p$-value. Despite the vast literature on urn models and multinomial-driven statistics, this intrinsic $p$-value appears to have remained unnoticed or, at best, underappreciated. Crucially, its definition requires no underlying random-variable model, while naturally extending both beyond the multinomial setting to broader classes of sampling and beyond prescribed values of generalized moments to algebraic relations among them. 

Not only does a sampling-based $p$-value for generalized moments intrinsically exist in any urn model, but its computation is also tractable --\, at least in the ideal scenario of prescribed generalized moments at large sample size\,-- due to an  emergent  spherical symmetry in the space of probability masses. This allows us to derive a semi-analytic closed form in Laplace approximation, which is universally --\,meaning for arbitrary generalized moments\,-- expressed  in terms of canonical information-theoretic quantities. The tractability of the intrinsic  $p$-value may be used to motivate large-deviation geometry~\cite{dembo1998large}, but it does not commit to such a geometry as a fundamental principle. 

To make the preceding ideas concrete, we now formulate the forward problem of sampling generalized moments within an empirical framework. This formulation embodies the proposed paradigm shift.

\subsection{An invitation to the forward problem}
\label{sc:ForwardProblem}

Without committing to a specific application, we consider a well-defined dataset to be characterized by a complete data matrix consisting of $N$ observed entities, which are distinguishable up to their recorded categories --\,that is, by virtue of exhibiting inequivalent rows within the data matrix. 
A category can be understood as a nominal or metric  realization $x$ of an observable (or collection of observables) $X$ defined over a theoretical domain $\mathcal X$. 
Whether a chosen observable is also a sensible one remains a separate, application-dependent question that we do not pursue here. 
In fact, $x$ need not even stem from some well-defined measurement operation of an observable, but can represent any measurable manifestation of a more abstract concept. 
In this spirit, we do not postulate random variables and their underlying distribution as primitive objects. 

Upon receiving a concrete dataset, we may instead restrict attention to a finite subset of observable categories $\Omega\subseteq\mathcal X$ that encompasses all observed categories within that dataset. 
Any further dataset generated by an empirical procedure could only yield finitely many realizations of $X$.  Furthermore, finite measurement precision ---\,alongside additional practical constraints\,--- renders sufficiently close elements of $\mathcal X$ empirically indistinguishable. 
Consequently, a restriction to some finite $\Omega$ is always possible without loss of empirical content. As a result, any dataset admits a histogram representation as a count vector $\mathbf n\in\mathbb N^{\vert\Omega\vert}_0$.

Conceptually, we think of an empirical procedure as a physically implementable (e.g.\ by a robot) sampling process that generates datasets according to a map 
\begin{equation}
    \operatorname{Pr} :~ \mathbb N_0^{\vert\Omega\vert}
    \rightarrow [0,1]~.
\end{equation}
Operationally, this map assigns probabilities to datasets according to various rules that are mostly unknown or, for the purposes of one's analysis, irrelevant. In that sense, the probability $\operatorname{Pr}(\mathbf n)$ of sampling a dataset $\mathbf n$ characterizes either effective (arising from lack of knowledge or interest) or intrinsic (inherent in the target system) uncertainty. 

Under the premise of forward sampling, we take the expected relative frequencies of all observable categories to be known, 
\begin{equation}
\label{eq:Intro:expectedCountNumber}
    \sum_{\mathbf n\in\mathbb N^{\vert\Omega\vert}_0}  
    n(x)\, \operatorname{Pr}(\mathbf n)
    = \langle n(x)\rangle \quad\text{for}\quad x\in\Omega~,
\end{equation}
regardless of what else may remain unknown about the target system or the underlying population. After normalizing by the expected total count $N=\sum_{x\in\Omega}\langle n(x)\rangle$, we can express Eq.~\eqref{eq:Intro:expectedCountNumber} through a reference distribution. This concept parallels the more conventional notion of a ``model'' or ``population'' distribution.

\paragraph{Families of datasets and the $p$-value.}
Given some sampling process, the forward problem investigates the probability mass of dataset families $\Lambda\subseteq\mathbb N_0^{\vert\Omega\vert}$
which is estimated by the restricted sum 
\begin{equation}
\label{eq:Intro:FamilyProbabilityMass}
    \operatorname{Pr}(\Lambda) = \sum_{\mathbf n\in\mathbb N_0^{\vert\Omega\vert}}
    \delta(\mathbf n\in \Lambda)\, \operatorname{Pr}(\mathbf n)
\end{equation}
over the non-negative integer orthant.
$\delta$ in discrete summations signifies the indicator function enforcing membership to region $\Lambda$. 
If for example datasets of a fixed (not necessarily large) sample size $N$ are considered, then we simply obtain the Kronecker delta 
\begin{equation}
\label{eq:Intro:KroneckerDelta}
    \delta\!\left(\mathbf n\in\Lambda\right)\equiv\delta\!\left(\sum_{x\in\Omega}n(x),N\right)~.
\end{equation}

Beyond the obvious structural condition of normalization, additional structural conditions may well apply to observable count vectors, depending on the problem at hand. 
Let us denote by $\Lambda_S\subseteq\mathbb N^{\vert\Omega\vert}_0$ the family formed by all datasets compatible with the prespecified structural conditions. 
Requiring the probability to observe any dataset in $\Lambda_S$ 
to be unity, 
the conditional sampling probability can be written as
\begin{equation}
\label{eq:Intro:ConditionalSamplingProbability}
    \operatorname{Pr}(\mathbf n\vert\Lambda_S) = \frac{\delta\!\left(\mathbf n\in\Lambda_S\right)\operatorname{Pr}(\mathbf n)}{\operatorname{Pr}(\Lambda_S)}
    ~.
\end{equation}
Combining Eq.~\eqref{eq:Intro:FamilyProbabilityMass} and Eq.~\eqref{eq:Intro:ConditionalSamplingProbability}, we obtain the coarse-grained probability mass of any nested family $\Lambda\subseteq\Lambda_S$ of datasets in the conditional sampling scheme,
\begin{equation}
    \label{eq:Intro:ConditionalSamplingProbabilityMass}
    \operatorname{Pr}(\Lambda\vert\Lambda_S) = \frac{\operatorname{Pr}(\Lambda)}{\operatorname{Pr}(\Lambda_S)}
    \,\in [0,1]
    ~,
\end{equation}
which is neatly given by the nested-to-structural family mass ratio. 

Given an index set $I$, we can partition the structural family 
by pairwise disjoint ($\Lambda_i\bigcap\Lambda_j=\emptyset$ for $i\neq j\in I$) 
nested families $\Lambda_i\subset\Lambda_S$ for $i\in I$, such that 
\begin{equation}
\label{eq:Intro:StructuralDecomposition}
    \Lambda_S = \bigcup_{i\in I} \Lambda_i ~.
\end{equation} 
For a selected family of datasets $\Lambda_{j}\subset\Lambda_S$ 
we can define 
an inherently discrete  
\begin{equation}
\label{eq:Intro:pValue}
    \text{$p$-value} = \sum_{i\in I_j} \operatorname{Pr}(\Lambda_i\vert\Lambda_S) \quad\text{with}\quad I_j = \left\{i\in I~\vert~ \operatorname{Pr}(\Lambda_i)\leq \operatorname{Pr}(\Lambda_j)\right\}~,
\end{equation}
through the cumulative probability mass of all nested families $\Lambda_i\subset\Lambda_S$ with equal or lower sampling probability than the given one.

\subsection{Paper contribution and outline}

\paragraph{Main result.}
When both the structural family and its partition into nested families are defined by $d$ and $D=d+k$ generalized moments, respectively, the $p$-value of Eq.~\eqref{eq:Intro:pValue} admits a straightforward interpretation. It
gives the cumulative probability of sampling, from the reference distribution (related to Eq.~\eqref{eq:Intro:expectedCountNumber}), the family of datasets supporting the $d+k$ observed generalized moments or any more improbable family, conditioned on supporting the $d$ structural generalized moments. 
At sufficiently large $N\gg d+k$, this quantity can then be approximated by the cumulative distribution function of a $\chi^2_k$ distribution. 

In this leading-order approximation in $N$, the test statistic Eq.~\eqref{eq:LaplaceApproximation:CriticalValue} exhibits the closed form of the \textsc{kl} divergence rescaled by $2N$. 
We numerically determine its two arguments through the information projection ($I$-projection)
of the reference distribution  on the nested family defined by the observed generalized moments and on the ambient family defined by the structural generalized moments, respectively. The $I$-projection on a linear family of distributions can be computed using Newton-Raphson root-finding algorithm.

\paragraph{Outline.} 
In Section~\ref{sc:GeneralizedMoments}, we define the core object whose probability mass we seek to investigate in this paper: the family of all admissible datasets that share the same values of the generalized moments under consideration. We outline both the generalized Diophantine system induced by this combinatorial definition  as well as the linear-algebraic formulation of its continuum limit, corresponding to the idealized limit of infinite sample size.

In Section~\ref{sc:ClassicalSampling}, we review independent and unbiased sampling, which becomes relevant at large $N$ by allowing us to regard the underlying process as effectively multinomial, irrespective of its microscopic details, which appear only as higher-order effects in $1/N$ (see Appendix~\ref{app:largeN}). The resulting $1/N$-expansion, within a compatible continuum limit for $N\gg \vert\Omega\vert$, uncovers the information-theoretic content (see Appendix~\ref{app:InformationGeometry}) of this universal sampling scheme. At the same time,  it already alludes to the tractability of estimating probability masses for $D$ generalized moments under conditional sampling.

In Section~\ref{sc:FourierAnalysis}, we relax the aforementioned limit to $D\ll N,\vert\Omega\vert$. In particular, the number of observable categories may now grow large, provided that the number $D$ of target generalized moments remains much smaller than the sample size $N$.
The probability mass density of the linear family 
is expressed as the inverse Fourier transform of the characteristic function of $\vert\Omega\vert$ independent Poisson processes restricted to the $D$-dimensional Fourier-dual subspace of the generalized moments. After analytic continuation, we compute this quantity to leading order in $N$ using saddle-point methods in Section~\ref{ssc:SaddlePointApproximation}.

In Section~\ref{sc:pValue}, we then rigorously define the probability mass of nested linear families under conditional sampling and show, within the Laplace approximation, that the intrinsic $p$-value is a natural object to consider.  Specifically, our systematic approximation of the discrete quantity Eq.~\eqref{eq:Intro:pValue} yields a semi-analytic closed-form statistic that  depends only on canonical information-theoretic quantities, without ever postulating information-geometric or probabilistic structures borrowed from continuous analysis.

\section{Linear Families of Datasets defined by Generalized Moments}
\label{sc:GeneralizedMoments}

Albeit conceptually clear, the sampling-based definition of $p$-value can prove hard (or practically impossible) to compute for most realistic problems of interest.  
Its computational complexity largely stems from the estimation of probability masses Eq.~\eqref{eq:Intro:FamilyProbabilityMass} for generic families of datasets. 
Depending on the form of the latter, 
the otherwise well-defined forward problem can quickly become computationally prohibitive owing to the combinatorial explosion in the number of admissible histogram configurations $\mathbf n$ at larger sample sizes.

Conceptually, the  combinatorial problem which is required to exactly compute the probability mass of a generic family of datasets 
amounts to characterizing all compatible  
count vectors 
that generate the defining features of that family. 
As we demonstrate, this classical combinatorial problem admits a semi-analytic large-$N$ treatment, when the target features are 
specified as generalized moments.
While derived here for linear conditions on generalized moments, the semi-analytic closed-form solution may be recycled 
to approximate broader classes of problems, including biased sampling and non-convex lattice regions associated with non-linear generalized moment conditions.

\paragraph{Families defined by linear conditions on raw means.}
A scientific theory expresses expectations about observables $X$ through a family of distinct transformations $g_\alpha:\mathcal X\rightarrow\mathbb R$ for $\alpha=0,\ldots,D-1$
that remain well-behaved at least over the observable categories $\Omega$.
Specifying such a collection of functions, 
it is straight-forward to compute raw means in some dataset by  
\begin{equation}
\label{eq:Intro:RawMean}
    m_\alpha 
    = \sum_{x\in\Omega} g_\alpha(x)\,n(x)  \equiv \vec g_\alpha\cdot\mathbf n
    \quad\text{for}\quad \alpha=0,\ldots,D-1~,
\end{equation}
i.e.\ as the scalar product of the column count vector $\mathbf n\in\mathbb N^{\vert\Omega\vert}_0$  with the row vector $\vec g_\alpha\in\mathbb R^{1\times\vert\Omega\vert}$ formed by the $g_\alpha$-image of $\Omega$.
For a given vector of raw means $\boldsymbol m\in\mathbb R^D$, there are generically many (or possibly none)   
compatible datasets $\mathbf n\in\mathbb N^{\vert\Omega\vert}_0$ 
that support them. 
If we are interested in datasets whose raw means $m_\alpha$ satisfy certain polynomial (in)equalities, the corresponding family of count vectors generically defines a (semi-)algebraic set $\Lambda$ intersecting the non-negative integer orthant $\mathbb N_0^{\vert\Omega\vert}$.

In the case that raw means satisfy linear conditions, we can conveniently summarize the prescription of Eq.~\eqref{eq:Intro:RawMean} through the coefficient matrix
\begin{equation}
\label{eq:intro:CoeffientMatrix}
    \mathbf G = 
    \begin{pmatrix}
    \vec g_0 \\
    \vdots\\
    \vec g_{D-1}
    \end{pmatrix} 
    \in \mathbb R^{D\times \vert\Omega\vert}~,
\end{equation} 
whose rows are indexed by $\alpha=0,\ldots,D-1$ and columns by observable $x\in\Omega$. 
Whenever the generalized Diophantine  system of linear equations
\begin{equation}
\label{eq:DiophantineSystem}
    \mathbf G\: \mathbf n = \boldsymbol m \quad\text{for}\quad \mathbf n\in\mathbb N_0^{\vert\Omega\vert} ~. 
\end{equation} 
is consistent, its set of solutions 
\begin{equation}
\label{eq:DatasetFamily}
    \Lambda(\mathbf G;\boldsymbol m)  = \left\{\mathbf n\in\mathbb N^{\vert\Omega\vert}_0~\vert~ \mathbf G\: \mathbf n = \boldsymbol m\right\}
\end{equation}
is isomorphic to the intersection of the algebraic variety representing a (generically unbounded) convex polyhedron with the non-negative integer orthant 
$\mathbb N^{\vert\Omega\vert}_0$.

Given a set of conditions on count vectors encoded by coefficient matrix $\mathbf G$, admissible dataset families $\Lambda_{\boldsymbol m}\equiv\Lambda(\mathbf G,\boldsymbol m)$ can be parametrized by the $D$-dimensional vector of possible raw means 
\begin{equation}
\label{eq:RawMeans:G_image}
    \boldsymbol m\in \mathbf G\:\mathbb N^{\vert\Omega\vert}_0\subset\mathbb R^D
\end{equation}
that are supported in the image $\mathbf G\:\mathbb N^{\vert\Omega\vert}_0$ of the non-negative integer orthant
under the $\mathbf G$-induced linear map.  
In this context, we abbreviate their probability mass Eq.~\eqref{eq:Intro:FamilyProbabilityMass} by 
\begin{equation}
\label{eq:SamplingProbabilityMass:Definition}
    \operatorname{Pr}(\boldsymbol{m})\equiv\operatorname{Pr}(\Lambda_{\boldsymbol m})
    =
    \sum_{\mathbf n\in\mathbb N^{\vert\Omega\vert}} \operatorname{Pr}(\mathbf n)\, \delta(\mathbf n\in\Lambda_{\boldsymbol{m}})
    ~.
\end{equation}
In the conventional limit of large sample sizes, this prescription for the probability mass can be rigorously estimated by the saddle-point approximation 
derived in Section~\ref{sc:FourierAnalysis}. 

In practice given $D=d+k$ functions, we usually wish to investigate how the sampling probability mass $\operatorname{Pr}(\boldsymbol{m})$ changes, when we vary some of their raw means $m_{d+i}$ for $i=0,\ldots,k-1$ while keeping the first $m_\alpha=\hat m_\alpha$ for $\alpha=0,\ldots,d-1$ fixed. Notice that different vectors of raw means from $\mathbf G\:\mathbb N_0^{\vert\Omega\vert}$ in Eq.~\eqref{eq:RawMeans:G_image} generate distinct families of datasets, since a given dataset can support at most one $\boldsymbol{m}$. Hence, varying $m_{d+i}$ induces, as in Eq.~\eqref{eq:Intro:StructuralDecomposition}, a partition of the structural family of datasets determined by the fixed $\hat m_\alpha$ into nested linear families. 

For the probability mass of a nested family of datasets Eq.~\eqref{eq:Intro:ConditionalSamplingProbabilityMass} in such a conditional sampling scheme, we succinctly write 
\begin{equation}
    \operatorname{Pr}(\boldsymbol{m}\vert \hat m_0,\ldots,\hat m_{d-1})
\end{equation}
or even 
\begin{equation}
\label{eq:Intro:ConditionalSampling_ProbabilityMass}
    \operatorname{Pr}(m_{d},\ldots,m_{d+k-1}\vert \hat m_0,\ldots,\hat m_{d-1})~,
\end{equation}
if it is clear that all considered families satisfy anyhow the structural conditions. 
In Section~\ref{sc:pValue}, this conditional probability mass will become the main object of investigation.
Unless the sample size $N$ is implied by a linear combination of the other structural conditions, we shall reserve the first condition for the constant map $g_0(x) = 1$ $\forall\,x\in\Omega$ to force a fixed $N$.   Consequently, all polyhedra considered below are bounded, hence constituting polytopes.

\paragraph{Generalized moments.}
Normalizing histogram counts by the size $N$ of a dataset $\mathbf n$, we turn raw means Eq.~\eqref{eq:Intro:RawMean} to generalized moments 
\begin{equation}
\label{eq:Intro:GeneralizedMoment}
    \mu_\alpha \equiv \widebar{g_\alpha(X)} = \sum_{x\in\Omega} \frac{n(x)}{N} g_\alpha(x)~.  
\end{equation}
In the following, we will often use that $\boldsymbol{m}=N\boldsymbol{\mu}$ when working at fixed sample size $N$.

Besides the constant map $g_0(x)=1$ $\forall\,x\in\Omega$ that extracts the size of the dataset, 
frequently occurring simple setups in descriptive statistics are based on indicator function $g_\alpha(x) = \delta(x\in S)$ for marginalization over regions $S\subseteq\Omega$ or monomial function $g_\alpha(x) = x^n$  to extract the $n$-th moment. Also more complicated setups can be addressed. For example, a function 
\begin{equation}
   g_2(x) = (x - \hat\mu)^2 
\end{equation}
would extract the variance with 
\begin{equation}
\hat\mu 
= \sum_{x\in\Omega} \frac{n(x)}{N}\, x
\end{equation}
being the mean of $X$ in the data, which must be added as a {structural} condition via $g_1(x)=x$. 
Alternatively, choosing
\begin{equation}
g_2(x) = \frac{x \delta(x\in S)}{\hat f} - \frac{x \delta(x\not\in S)}{1 - \hat f}
\end{equation}
for the difference of conditional means presupposes  the structural condition  on the marginal probability of $S\subset\Omega$  
\begin{equation}
\hat f 
= \sum_{x\in\Omega} \frac{n(x)}{N}\, \delta(x\in S) ~,
\end{equation}
which is extracted using indicator function $g_1(x)=\delta(x\in S)$. 

If $X$ denoted the combined voltage $V$ and current observable $I$, the resistance in Ohm's law would be captured as the expectation of function $g(V,I)=V/I$. Its mean and higher moments in any realistic dataset can be then calculated via Eq.~\eqref{eq:Intro:GeneralizedMoment}. As a consequence of the empirical law, the resistance should be independent of the voltage (conversely current), in expectation.
Which observables are suitable for inference, as well as what functional form a fundamental or empirical law might take, constitute highly interesting epistemological questions,
see e.g.\ discussion in Section~4.9 of~\cite{caticha2008lecturesprobabilityentropystatistical}. The semi-analytic formulation pursued here begins immediately after the identification of sensible observables and the functional form of a law in terms of those observables in scientific applications.

Estimating generalized moments via the sample mean of the corresponding function in 
Eq.~\eqref{eq:Intro:GeneralizedMoment} conceptualizes our mathematical expectation 
regarding the metric or nominal properties of observable $X$. 
For example, an expectation $\expec[g_\alpha(X)]$ that corresponds to a moment logically requires that the underlying observable has a metric scale. Likewise, an expectation that corresponds to a prevalence of some group presupposes that the entering observable has a nominal character. Although we collectively refer to $\mu_\alpha$'s as generalized moments, 
our formalism ---\,which computes them as averages over finitely many categories\,---  
effectively treats all observables $X$ as nominal.

\subsection{On the scale of variation of generalized moments}
\label{ssc:VariationScale}

Already elementary number theory suggests a well-defined continuum approximation in the $D$-dimensional space of generalized moments, at least as long as $N\gg D$. 
In the simplest case of one non-trivial generalized moment defined by rational-valued\footnote{Whenever the coefficient matrix contains some irrational element, the continuum approximation becomes more reliable at even smaller $N$.} function $g_1:\Omega\rightarrow\mathbb Q$ (when measurement precision is limited to e.g.\ $L$ decimal places), we exemplify this fact by investigating the $g$-image of $\mathbb N_0^{\vert\Omega\vert}$ subject to normalization condition. 

First, we may rewrite Diophantine system Eq.~\eqref{eq:DiophantineSystem} as 
\begin{equation}
\label{eq:ExampleDiophantineSystem}
    \sum_{x\in\Omega} n(x) = N \quad\text{and}\quad \sum_{x\in\Omega} n(x) \frac{a(x)}{L} = m ~,
\end{equation}
where we parametrized rational matrix elements by 
\begin{equation}
    g_1(x)=\frac{a(x)}{L} \quad\text{with}\quad a(x)\in\mathbb Z \quad\text{and} \quad L \in\mathbb N ~. 
\end{equation}
Consider now a non-trivial change $\Delta \mathbf n\in\mathbb Z^{\vert\Omega\vert}$ of the count configuration $\mathbf n\in\mathbb N_0^{\vert\Omega\vert}$ such that Eq.~\eqref{eq:ExampleDiophantineSystem} gives 
\begin{equation}
\label{eq:ExampleDiophantineSystem2}
    \sum_{x\in\Omega} \Delta n(x) = 0 \quad\text{but}\quad \sum_{x\in\Omega} \Delta n(x) \frac{a(x)}{L} = \Delta m ~.
\end{equation}
By essentially neglecting the non-negativity constraint on counts in the regime $n(x)\gg1$, we can effectively treat the possible changes as independent integer coordinates $\Delta n(x)=0,\pm1,\pm2,\ldots$ for each $x\in\Omega$ only subject to Diophantine system Eq.~\eqref{eq:ExampleDiophantineSystem2}.

Solving former fixed-$N$ condition for some reference $x_0\in\Omega$ and substituting into the latter condition condenses the problem to one Diophantine equation 
\begin{equation}
    \sum_{x\in\Omega} \Delta n(x) \left[a(x)-a(x_0)\right] = L\,\Delta m~.
\end{equation}
According to Bézout's identity, a linear combination of the integer set formed by $a(x)-a(x_0)$ for $x\in\Omega\setminus\{x_0\}$ on the l.h.s.\ ranges over all integer multiples of their greatest common divisor for unconstrained lattice vector $\Delta\mathbf n\in\mathbb Z^{\vert\Omega\vert}$. 
In the aforementioned regime of large counts, the r.h.s.\ must hence be divisible by this greatest common divisor.
Therefore, the effective lattice spacing in $m$ is 
\begin{equation}
    \Delta m_\text{min} = L^{-1}\operatorname{gcd}\left\{\vert a(x)-a(x_0)\vert\right\}_{x\in\Omega} = L^{-1}\operatorname{gcd}\left\{\vert a(x)-a(x')\vert\right\}_{x,x'\in\Omega} ~.
\end{equation}
In total, we thus have for the smallest possible (non-zero) variation of the corresponding generalized moment 
\begin{equation}
    \Delta\mu_\text{min} = \frac{\operatorname{gcd}\left\{\vert a(x)-a(x')\vert\right\}_{x,x'\in\Omega}}{N\,L}~.
\end{equation}

At sufficently large $N$, every single integer multiple of $\Delta m_\text{min}$ becomes achievable within a massive middle range of $m$ (the so-called ``bulk'') supported by count vectors $\mathbf n\in\mathbb N^{\vert\Omega\vert}_0$ well inside the associated polytope. 
Consequently, the possible values of the generalized moment $\mu\in\mathbb Q$ form an increasingly dense grid of points well inside its observable region.
Intuitively, the emergence of such a canonical structure in the space of generalized moments alludes to the successful implementation of saddle-point techniques in the continuum limit of Section~\ref{sc:FourierAnalysis}, even when $\vert\Omega\vert=\order{N}$.

\subsection{The induced linear algebra}
\label{ssc:LinearAlgebra}

Motivated by the previous number-theoretic result that the lattice spacing of generalized moments scales as $\order{N^{-1}}$, it is convenient to express counts through relative frequencies
\begin{equation}
\label{eq:FromCountsToProbabilities}
    p(x) = \frac{n(x)}{N} \in\mathbb Q~,
\end{equation}
which later admit a continuum treatment in the large-$N$ limit. 
Over the finite set $\Omega$ of observable categories, we define next the observable simplex
\begin{equation}
    \label{eq:simplex}
    \mathcal P(\Omega) = \left\{\prob p\in\mathbb R^{\vert\Omega\vert}\quad\text{with}\quad p(x)\geq0 ~\text{ and }~\sum_{x\in \Omega} p(x)=1\right\}~,
\end{equation}
where (real-valued) probability distributions  live. Operating over the same $\Omega$, we abbreviate $\mathcal P\equiv\mathcal P(\Omega)$.
For its $N$-dilation, we succinctly write $N\mathcal P$ to denote the uniformly scaled polytope in $\mathbb R^{\vert\Omega\vert}$ whose integer lattice points correspond to admissible count configurations at total size $N$.

The restriction of the forward problem to an observable simplex induces a linear-algebra formulation.
Due to the inherently discrete and finite nature of $\Omega\subseteq\mathcal X$, one can think of distributions $\prob p:\Omega\rightarrow[0,1]^{\vert\Omega\vert}$ as non-negative $\vert\Omega\vert$-dimensional column vectors.
The well-defined estimator of generalized moment $\mu_\alpha$ 
introduced in Eq.~\eqref{eq:Intro:GeneralizedMoment}
can be then written as the scalar product between the row vector $\vec g_\alpha\in\mathbb R^{1\times\vert\Omega\vert}$ 
and the column vector $\prob p\in\mathbb R^{\vert\Omega\vert}_{\geq0}$ of a distribution on the observable simplex.

In analogy to the generalized Diophantine system of linear equations introduced in Eq.~\eqref{eq:DiophantineSystem}, we arrive at the compact form of linear system
\begin{equation}
\label{eq:LinearSystemDEF}
    \mathbf G \: \prob p = \boldsymbol{\mu}
\end{equation}
 that is controlled by the column vector 
\begin{equation}
    \boldsymbol\mu = 
    \begin{pmatrix}
    \mu_0 \\
    \vdots\\
    \mu_{D-1}
    \end{pmatrix} \in \mathbb R^D
    ~,
\end{equation}
such as the generalized moments are reproduced. In the following for $D<\vert\Omega\vert$, we take any linear system in Eq.~\eqref{eq:LinearSystemDEF} to be consistent, i.e.\ to admit at least one solution over $\mathbb R^{\vert\Omega\vert}$, as well as non-redundant, i.e.\ its coefficient matrix Eq.~\eqref{eq:intro:CoeffientMatrix} to be of full row rank,  $D=\rank\mathbf G$.
To work at fixed (and large) $N$, we shall assume throughout that normalization is implied by the coefficient matrix, i.e.\ 
\begin{equation}
\label{eq:NormalizationSymmetry}
    \exists\,\boldsymbol w\in\mathbb R^D \quad\text{s.t.\ }\quad \boldsymbol w^T\, \mathbf G = \left(1,\ldots,1\right)\in\mathbb R^{1\times\vert\Omega\vert}~.
\end{equation} 

The convex solution space of linear system Eq.~\eqref{eq:LinearSystemDEF} together with the non-negativity constraint gives rise to a linear family of distributions on the observable simplex,
\begin{equation}
    \label{eq:LinearFamilyDEF}
    \mathcal L(\mathbf G;\boldsymbol\mu) 
    = \left\{\prob p\in \mathbb R^{\vert\Omega\vert}_{\geq0} ~\vert~ \mathbf G\:\prob p = \boldsymbol\mu \right\} \subseteq \mathcal P~.
\end{equation}
Geometrically, this is realized as the polytope obtained by intersecting the affine solution space with the non-negative orthant $\mathbb R^{\vert\Omega\vert}_{\geq0}$.
The intersection of its $N$-dilation with the integer lattice is isomorphic under Eq.~\eqref{eq:FromCountsToProbabilities} to the family of datasets Eq.~\eqref{eq:DatasetFamily} which recovers the discrete problem, $\Lambda(\mathbf G;\boldsymbol{m}) = N \mathcal L(\mathbf G;\boldsymbol\mu)\cap \mathbb Z^{\vert\Omega\vert}$.
Besides linearity, we further assume that the conditions on generalized moments do not force the probability of any observable category to vanish, i.e.\ there are no sampling or structural zeros:
\begin{equation}
\label{eq:NoOutliers}
    \forall x\in\Omega\,: \quad \exists\, \prob p \in\mathcal L(\mathbf G,\boldsymbol{\mu}) \quad \textit{s.t.}\quad p(x) > 0~.
\end{equation}
In future work, we will explore how the linear-algebraic language serves as a basis for a systematic treatment of more general sampling problems beyond linear conditions on generalized moments.

\section{Probability masses in the classical theory of sampling}
\label{sc:ClassicalSampling}

Given a sufficiently large sample of size $N$ drawn from a target system or population, its statistical analysis admits an effective description governed by the combinatorial structure among a large number of relevant entities, distinguishable only up to the chosen observational categories.
Therefore, many sampling-related quantities can be treated in simplified form using large-$N$ techniques. This perspective --\, already outlined in the introduction\,-- reflects a standard, often implicit assumption in statistical physics that macroscopic structure admits a combinatorial origin in the counting of microscopic configurations.

By 1878, Ludwig Boltzmann~\cite{boltzmann_weitere_2012,Bach1990} had not only contemplated multinomial sampling but had also performed a large-$N$ expansion of multinomial sampling probabilities. Furthermore, he derived the emerging divergence measure and the exponential form of its non-trivial minimizer. These constitute core structures~\cite{Akaike1985} that would later reappear in the work of Claude Shannon, Solomon Kullback, and Richard Leibler during the first half of the 20th century, before being further formalized by E.\ T.\ Jaynes and Imre Csiszár.  

Most crucially, Boltzmann's considerations did not assume any particular functional form for the energy of the underlying physical system. This demonstrates the generality of his sampling scheme, which we would nowadays describe as the standard sampling of colored balls with replacement from an urn, where energy levels play the role of colors. 
In terms of traditional sampling from an urn (thought of as a statistical ``population''), this classical regime corresponds to the combined limit of large population size $M$ and large sample size $N$. Whenever $M\gg N\gg1$, we may thus recycle Boltzmann's theory to obtain the leading-order description, 
even when independent sampling is reasonably violated, see Appendix~\ref{app:largeN} as well as~\cite{dembo1998large}. 

In the forward problem, we operate under the assumption that the ``population'' composition is known in the sense of Eq.~\eqref{eq:Intro:expectedCountNumber}. Hereafter, it is represented by reference distribution $\prob v\in\mathcal P$, which assigns strictly positive probabilities to observable categories, such that $v(x)=N^{-1}\langle n(x)\rangle>0$ for all $x\in\Omega$. Provided that we can for our practical purposes assume independent and unbiased sampling\footnote{In such a sampling regime where future draws are effectively independent of past draws, while sampled frequencies effectively reflect the population probabilities up to stochastic fluctuations, any deviations from the simplifying features are \(1/N\)-suppressed.}
the count vector $\mathbf n\in\mathbb N^{\vert\Omega\vert}_0$ follows the multinomial distribution
\begin{equation}
\label{eq:MultinomialDistribution}
    \operatorname{mult}(\mathbf n;\prob v) = \frac{N!}{\prod_{x\in\Omega}n(x)!} \prod_{x\in\Omega} v(x)^{n(x)}
\end{equation}
subject to normalization condition 
\begin{equation}
    \sum_{x\in\Omega} n(x) = N~.
\end{equation}

As this paradigm underpins all subsequent departures from the simplifying assumptions, we shall refer to it as the  forward problem of classical sampling. Beyond its central role in the semi-analytic approximation of the coarse-grained probability masses of dataset families, the classical sampling scheme admits a direct correspondence with classical statistical mechanics through Maxwell–Boltzmann statistics. Within this analogy, departures from the classical assumptions correspond to alternative statistical frameworks, such as Bose–Einstein and Fermi–Dirac statisticss~\cite{2008arXiv0808.2102N,Niven2009}.

\subsection{The scaling of sampling probabilities}
\label{ssc:NaiveLargeN}

The limiting behavior of independent unbiased sampling is just a consequence of an ordinary {large-$N$ expansion} of multinomial sampling probability, famously performed in~\cite{sanov_probability_1957}, and in fact, for the first time in~\cite{boltzmann_weitere_2012}. Using Stirling's formula~\eqref{eq:app:Stirling}, we give the expansion in $\log$-space:
\begin{equation}
\label{eq:Multinomial:LargeNexpansion}
    \log \operatorname{mult}(N\prob p;\prob v) = - N \infdiv{\prob p}{\prob v}
    - \frac{\vert\Omega\vert-1}{2}\log\left(2\pi N\right)  - \tfrac{1}{2} \sum_{x\in\Omega}\log p(x)  + \order{1/N}~.
\end{equation}
Over observable simplex Eq.~\eqref{eq:simplex}, the Information or {Kullback-Leibler divergence} (abbreviated \textsc{kl};~\cite{kullback_information_1951} also called relative entropy) of $\prob p$ from $\prob v$,
\begin{equation}
\label{eq:intro:KL}
   \infdiv{\prob p}{\prob v} = \sum_{x\in\Omega} p(x) \log \frac{p(x)}{v(x)} = -H(\prob p) - N^{-1}\ell_{N\prob p}(\prob v) ~,
\end{equation}
is the canonical  paradigm of divergence measure --\,the first to be discovered in this context and thoroughly investigated in physics and information theory. $0\cdot\log0=0$ is understood, as always.
The right-hand side in Eq~\eqref{eq:intro:KL} emphasizes the two contributions to the \textsc{kl} divergence coming from (i) the entropy of $\prob p$, 
\begin{equation}
\label{eq:intro:Entropy}
    H(\prob p) = -\sum_{x\in\Omega} p(x)\log p(x)
\end{equation}
arising from the large-$N$ combinatorics of multinomial coefficient\footnote{The multinomial coefficient $N! \left[\prod_{x\in\Omega} (Np(x))!\right]^{-1}$ counts all possible ways to partition $N$ distinct entities into distinct groups of sizes $Np(x)$ for $x\in\Omega$.}  and (ii) the log-likelihood 
\begin{equation}
\label{eq:intro:logLikelihood}
   \ell_{N\prob p}(\prob v) = \sum_{x\in\Omega} N p(x)\log v(x)
\end{equation}
of a dataset $N\prob p\in\mathbb N^{\vert\Omega\vert}_0$ under the reference distribution $\prob v\in\mathcal P$. 

In the  regime $N\gg\vert\Omega\vert$, where the expansion is formally valid, 
the emergent rather than axiomatic~\cite{9ac0edfe-5642-30a9-b54d-2ebb820d3810} character of \textsc{kl} divergence, which is often overlooked, can be summarized by a scaling law:
\begin{equation}
\label{eq:SamplingProbability:ScalingLaw}
    -\lim_{N\rightarrow\infty} N^{-1}\log \operatorname{Pr}(N\prob p) = \infdiv{\prob p}{\prob v}~.
\end{equation}
In the language of large deviation theory~\cite{LargeDeviations_PhysicsReview,burenev2025introduction}, the r.h.s.\ plays the role of the rate function. It describes how fast the sampling probability ``concentrates'' around its mode, which in absence of non-trivial conditions coincides with the reference distribution $\prob v$. 

Next, we investigate the implications of this scaling law for the coarse-grained probability mass $\operatorname{Pr}(N\boldsymbol{\mu})$ in Eq.~\eqref{eq:SamplingProbabilityMass:Definition} for a linear family $\Lambda(\mathbf G;\boldsymbol{m})$ defined by generalized moments $\boldsymbol{\mu}=\boldsymbol{m}/N$ (recall that we overload $\operatorname{Pr}$ to also denote the induced probability mass of families of datasets).
From Eq.~\eqref{eq:SamplingProbability:ScalingLaw}, we quickly recognize that 
the exponential scaling of sampling probabilities suppresses atypical datasets. 
Given a dataset $N \prob p$ and a more ``exotic'' dataset $N \prob p'$ in the  information-theoretic sense of $\infdiv{\prob p'}{\prob v}> \infdiv{\prob p}{\prob v}$, the contribution of latter to the probability mass Eq.~\eqref{eq:SamplingProbabilityMass:Definition} would be relatively suppressed  by  
\begin{equation*}
    \log\frac{\operatorname{Pr}(N\prob p')}{\operatorname{Pr}(N\prob p)} = 
    \order{ -N \left[\infdiv{\prob p'}{\prob v}-\infdiv{\prob p}{\prob v}\right] }~.
\end{equation*}

In the theoretical limit of infinitely large datasets, 
only the unique (see Appendix~\ref{app:InformationGeometry}) minimum of \textsc{kl} divergence over the associated linear family $\mathcal L_{\boldsymbol\mu}\equiv\mathcal L(\mathbf G;\boldsymbol\mu)$ survives:
\begin{equation}
    \label{eq:SamplingProbabilityMas:ScalingLaw}
        \lim_{N\rightarrow\infty}N^{-1}\log\operatorname{Pr}(N\boldsymbol{\mu}) = -\min_{\prob p\in\mathcal L_{\boldsymbol{\mu}}}\infdiv{\prob p}{\prob v} 
        ~.
\end{equation}
At finite $N$, the minimizing distribution, the $I$-projection of $\prob v$ onto $\mathcal L_{\boldsymbol\mu}$ itself,  
\begin{equation}
     \prob q_{\boldsymbol{\mu}} = \argmin_{\prob p\in\mathcal L_{\boldsymbol\mu}} \infdiv{\prob p}{\prob v}~, 
\end{equation}
does not need to belong to the Diophantine set of solutions, as generally $N q_{\boldsymbol{\mu}}(x)\not\in \mathbb N_0$. At sufficiently large, but still finite $N$, scaling law Eq.~\eqref{eq:SamplingProbabilityMas:ScalingLaw} tells us that  $N\infdiv{\prob q_{\boldsymbol{\mu}}}{\prob v}$ gives a good approximation to the coarse-grained probability mass. In particular,  the contribution from the specifics of the integral polytope become exponentially suppressed.

Keeping certain generalized moments $\hat{\boldsymbol\mu}\in\mathbb R^d$ fixed, 
\begin{equation}
\label{eq:ConditionalSampling:LinearSystem}
    \mathbf n\in\mathbb N_0^{\vert\Omega\vert}:\quad\mathbf G_S\mathbf n \overset{!}{=}N\hat{\boldsymbol\mu}~,
\end{equation}
the conditional sampling probability of datasets within the structural family $\hat{\mathcal L}\equiv\mathcal L_{\hat{\boldsymbol\mu}}$  follows from Eq.~\eqref{eq:Intro:ConditionalSamplingProbability}:
\begin{equation}
\label{eq:ConditionalSamplingProbability:GivenGeneralizedMoments}
    \operatorname{Pr}(\mathbf n\vert\hat{\boldsymbol\mu}) = \frac{\operatorname{Pr}(\mathbf n)\,\delta(\mathbf n\in N\hat{\mathcal L})}{\operatorname{Pr}(N \hat{\boldsymbol\mu})}
    \quad\text{for}\quad
    \mathbf n\in\mathbb N_0^{\vert\Omega\vert}\cap N \hat{\mathcal L}~.
\end{equation}
If only normalization on the simplex applies ($d=1$), we recover the unconditional multinomial sampling. 
Combining Eq.~\eqref{eq:SamplingProbability:ScalingLaw} with Eq.~\eqref{eq:SamplingProbabilityMas:ScalingLaw}, we readily derive the scaling law for conditional sampling,
\begin{equation}
\label{eq:ConditionalSampling:ScalingLaw}
    \lim_{N\rightarrow\infty}\log N^{-1}\operatorname{Pr}(N\prob p\vert\hat{\boldsymbol\mu}) = -\infdiv{\prob p}{\hat{\prob q}}~,
\end{equation}
after applying Pythagorean identity Eq.~\eqref{eq:PythagorasIdentity} for the $I$-projection of $\prob v$ onto ambient $\hat{\mathcal L}$,
\begin{equation}
\label{eq:ConditionalSampling:iProjection}
     \hat{\prob q} 
     = \argmin_{\prob p\in\hat{\mathcal L}} \infdiv{\prob p}{\prob v}
\end{equation}
that governs the probability mass of the structural linear family.

\subsection{The ``naive''  approximation to the probability mass}
\label{ssc:NaiveGaussianization}
This subsection mainly provides intuition for the subsequent rigorous derivation in Sections~\ref{sc:FourierAnalysis} and~\ref{sc:pValue}. It may therefore be skipped by the reader familiar with Jaynes' theorem of entropy concentration~\cite{jaynes_concentration_1989}, see also chapter~11 of~\cite{jaynes_probability_2003}.

Starting from scaling law Eq.~\eqref{eq:ConditionalSampling:ScalingLaw}, it is easy to generalize Jaynes' framework to a non-uniform reference distribution $\prob v$. 
Most of the sampling probability Eq.~\eqref{eq:ConditionalSamplingProbability:GivenGeneralizedMoments}
shall peak in the vicinity of the $I$-projection Eq.~\eqref{eq:ConditionalSampling:iProjection}. 
Therefore, we should investigate probability fluctuations within the structural linear family around the $I$-projection, i.e.\ for $\prob p\in\hat{\mathcal L}\equiv\mathcal L(\mathbf G_S;\hat{\boldsymbol\mu})$.
To this end, the large-$N$ expansion  of conditional sampling probability Eq.~\eqref{eq:ConditionalSamplingProbability:GivenGeneralizedMoments} following from Eq.~\eqref{eq:Multinomial:LargeNexpansion} becomes an expansion around the governing $I$-projection\footnote{In the ideal scenario of Section~\ref{app:IdealProblem}, the $I$-projection provides a stable expansion point, well inside the interior of $\mathcal L$, as outlined in Section~\ref{app:iProjections}.}: 
\begin{align}
\label{eq:SamplingProbability:QuadraticForm}
    \log\text{Pr}(N\prob p\vert\hat{\boldsymbol\mu}) 
    =\,&
    -N \infdiv{\prob p}{\hat{\prob q}} - \log Z
    \\[1ex]
    =\,&
    - \sum_{x\in\Omega} \frac{\left[\sqrt{N}\left(p(x)-\hat q(x)\right)\right]^2}{2\hat q(x)} + \order{N^{-1/2}} - \log Z 
    ~,
    \nonumber
\end{align}
upon remarking that normalization on the simplex requires $\sum_{x\in\Omega}\left[\prob p(x)- \hat q(x)\right]=0$. The normalizing partition function $Z$ will be fixed by consistency of large-$N$ expansion, below.  

Since by definition both $\prob p,\hat{\prob q}$ belong to the structural family $\hat{\mathcal L}$,
probability fluctuations live in the kernel of the structural coefficient matrix from Eq.~\eqref{eq:ConditionalSampling:LinearSystem}, meaning
\begin{equation}
    \left(\prob p - \hat{\prob q}\right)\in\ker\mathbf G_S \quad\Leftrightarrow\quad
    \mathbf G_S \left(\prob p - \hat{\prob q}\right) = \mathbf 0~.
\end{equation}
Hence, the quadratic form in Eq.~\eqref{eq:SamplingProbability:QuadraticForm} exhibits only $\vert\ker\mathbf G_S\vert=\vert\Omega\vert-d$ degrees of freedom, which we may parametrize by 
\begin{equation}
\label{eq:kernelFluctiations:Parametrization}
    \sqrt N \left(p(x) - q(x)\right) = \sum_{a=1}^{\vert\Omega\vert - d} \psi_a k_{a}(x)~,
\end{equation}
where the row-space of $\mathbf K\in\mathbb R^{(\vert\Omega\vert - d)\times \vert\Omega\vert}$ spans the kernel of our coefficient matrix, i.e.\ 
\begin{equation}
    \mathbf G_S \: \mathbf K^T = \mathbf 0~.
\end{equation}

The leading term in large $N$, which is quadratic in the rescaled probability difference 
\begin{equation}
\label{eq:SamplingProbability:ProbabilityFluctuations}
    \sqrt{N}\left(p(x)-\hat q(x)\right)\in[-\sqrt N,\sqrt N]~,
\end{equation} 
invites us to continuously approximate the conditional sampling probability 
via $\rho(\boldsymbol{\psi}\vert\hat{\boldsymbol\mu})d^{\vert\Omega\vert-d}\!\boldsymbol{\psi}$ using a density $\rho$ over $\boldsymbol{\psi}\in\mathbb R^{\vert\Omega\vert-d}$. 
Such a continuum limit --\,often formalized as weak convergence of sums to integrals~\cite{LargeDeviations_PhysicsReview} and widely employed in physics\,--  can be justified a posteriori, when evaluating sums via Laplace's method. Intuitively, the continuum limit is motivated by the increasing density of points 
within the ambient polytope associated to Eq.~\eqref{eq:ConditionalSampling:LinearSystem}. 

To this end, it is instructive to examine the lattice volume contained in the $N$-dilation of the simplex
\begin{equation}
\label{eq:VolumeOfIntegralLattice:StarsAndBars}
       \operatorname{vol}\left(N\mathcal P\cap\mathbb N_0^{\vert\Omega\vert}\right) =
       \sum_{\mathbf n\in\mathbb N^{\vert\Omega\vert}_0} \delta\!\left(\sum_{x\in\Omega}n(x),N\right) 
       = \binom{N + \vert\Omega\vert - 1}{\vert\Omega\vert - 1} \overset{N\gg\vert\Omega\vert}{\approx} \frac{N^{\vert\Omega\vert-1}}{(\vert\Omega\vert-1)!} ~,
\end{equation}
which exactly equals the binomial coefficient, see the stars and bars construction in combinatorics.
Under the ``naive'' continuum limit based on $\Delta p(x) = \Delta n(x)/N$, i.e.\ 
\begin{equation}
\label{eq:NaiveContinuumLimit}
    \Delta\mathbf n = \prod_{x\in\Omega}\Delta n(x)  \quad\overset{N\gg\vert\Omega\vert}{\longrightarrow}\quad 
    d^{\vert\Omega\vert}(N\prob p)=N^{\vert\Omega\vert}\prod_{x\in\Omega}dp(x) ~,
\end{equation}
which readily turns the Kronecker delta to Dirac delta function, we compute 
\begin{equation}
\label{eq:ContinoousApproximation:Intro}
    \int d^{\vert\Omega\vert}(N\prob p)\, \delta\!\left(\sum_{x\in\Omega}Np(x)-N\right) = N^{\vert\Omega\vert-1}  \operatorname{vol}\left(\mathcal P\right)~. 
\end{equation}
This agrees with the large-$N$ expansion of the volume of the Diophantine set in Eq.~\eqref{eq:VolumeOfIntegralLattice:StarsAndBars}, since 
the volume of the standard simplex in $\mathbb R^{\vert\Omega\vert}$ is given by 
\begin{equation}
    \operatorname{vol}\left(\mathcal P\right) = \left[\prod_{x\in\Omega}\int_0^1 dp(x)\right] \delta\!\left(\sum_{x\in\Omega}p(x)-1\right) 
    =\frac{1}{\left(\vert\Omega\vert-1\right)!}
    ~.
\end{equation}
Indeed for unconditional multinomial sampling ($d=1$), the $N$-factor from the ``naive'' differential Eq.~\eqref{eq:NaiveContinuumLimit},
\begin{equation}
    N^{\vert\Omega\vert-1}d^{\vert\Omega\vert}\prob p\,\delta\!\left(\sum_{x\in\Omega}p(x)-1\right)
    \overset{\eqref{eq:kernelFluctiations:Parametrization}}{=} N^{\frac{\vert\Omega\vert-1}{2}} d^{\vert\Omega\vert-1}\boldsymbol{\psi}~,
\end{equation} 
precisely cancels the subleading contribution in the large-$N$ expansion of Eq.~\eqref{eq:Multinomial:LargeNexpansion}, thereby signaling the consistency of the continuum limit.

Plugging the suggested parametrization Eq.~\eqref{eq:kernelFluctiations:Parametrization} for $\prob p\equiv\prob p(\boldsymbol{\psi})\in\hat{\mathcal L}$ into Eq.~\eqref{eq:SamplingProbability:QuadraticForm} thus leads to the density 
\begin{equation}
    \label{eq:Concentration:PsiExpansion}
    \log\rho(\boldsymbol{\psi}\vert\hat{\boldsymbol\mu}) 
    =
    - \tfrac12\sum_{a,b=1}^{\vert\Omega\vert - d} \psi_a\psi_b \sum_{x\in\Omega}\frac{k_a(x)k_b(x)}{\hat q(x)} + \order{N^{-1/2}}  - \log Z
    ~,
\end{equation}
in the aforementioned continuum limit. 
The quadratic form is well-behaved at $\Vert\boldsymbol{\psi}\Vert\rightarrow\infty$ when $N\rightarrow\infty$, since the precision  matrix $\mathbf \Sigma^{-1}\in\mathbb R^{(\vert\Omega\vert-d)\times(\vert\Omega\vert-d)}$ that is induced by elements
\begin{equation}
\label{eq:Concentration:PrecisionMatrix}
    (\mathbf\Sigma^{-1})_{ab} =   \sum_{x\in\Omega}\frac{k_a(x)k_b(x)}{\hat q(x)}
\end{equation}
is positive-definite for the same reasons as the Jacobian of Section~\ref{ssc:iProjector}, see in particular Eq.~\eqref{eq:iProjector:Jacobian}. 
Therefore, the density $\rho$ can be used to approximate the coarse-grained probability mass
\begin{equation}
    \sum_{\mathbf n\in\mathbb N_0^{\vert\Omega\vert}} \operatorname{Pr}(\mathbf n\vert\hat{\boldsymbol\mu}) \,\delta(\mathbf n\in N\mathcal R) 
    ~,
\end{equation}
within nested families $\mathcal R\subseteq\hat{\mathcal L}$  that contain the $I$-projection, i.e.\ $\prob p(\boldsymbol{\psi}=\boldsymbol0)\equiv\hat{\prob q}\in{\mathcal R}$, so that Laplace method is applicable, by replacing the discrete sum with the corresponding continuum integral. 
Using the $N$-leading multivariate Gaussian density from Eq.~\eqref{eq:Concentration:PsiExpansion}, the partition function can then be easily estimated:
\begin{equation}
    \sum_{\mathbf n\in\mathbb N_0^{\vert\Omega\vert}} \operatorname{Pr}(\mathbf n\vert\hat{\boldsymbol\mu}) \overset{!}{=} 1 \quad\rightarrow\quad
    Z \approx 
    \int d^{\vert\Omega\vert-d}\boldsymbol{\psi}\, \rho(\boldsymbol{\psi}\vert\hat{\boldsymbol\mu}) ~.
\end{equation}
$d^{\vert\Omega\vert-d}\boldsymbol{\psi}$ denotes the differential volume in the kernel space of the coefficient matrix. 

Specifically, we may consider a nested linear family $\tilde{\mathcal L}\equiv\mathcal L(\mathbf G;\tilde{\boldsymbol\mu})\subset\hat{\mathcal L}$ defined by generalized moments $\tilde{\boldsymbol\mu}\in\mathbb R^{d+k}$ which imply the generalized moments $\hat{\boldsymbol\mu}\in\mathbb R^d$ of the ambient family $\hat{\mathcal L}\subseteq\mathcal P$ in the linear-algebraic sense. For a thorough discussion of nested sampling, see Section~\ref{ssc:ConditionalSampling}.
Using scaling law Eq.~\eqref{eq:SamplingProbabilityMas:ScalingLaw} for the probability masses of ambient and nested linear families
leads after application of the Pythagorean identity Eq.~\eqref{eq:PythagorasIdentity} to the compact expression 
\begin{equation}
\label{eq:nestedKL}
    -\lim_{N\rightarrow\infty}N^{-1}\frac{\operatorname{Pr}(N\tilde{\boldsymbol\mu})}{\operatorname{Pr}(N\hat{\boldsymbol\mu})} = \infdiv{\tilde{\prob q}}{\hat{\prob q}}
\end{equation}
in terms of the inner $I$-projection 
\begin{equation}
    \tilde{\prob q} \equiv \argmin_{\prob p\in\tilde{\mathcal L}}\infdiv{\prob p}{\prob v} = \argmin_{\prob p\in\tilde{\mathcal L}}\infdiv{\prob p}{\hat{\prob q}}~,
\end{equation}
alongside the ambient $I$-projection Eq.~\eqref{eq:ConditionalSampling:iProjection}.

As argued above, a discrete summation over admissible $N\tilde{\boldsymbol\mu}$ at fixed $\hat{\boldsymbol\mu}$ can be approximated by an integral over the conditional density in Eq.~\eqref{eq:Concentration:PsiExpansion}, provided that the integration region includes the ambient $I$-projection.
Using the leading Gaussian density in $\vert\Omega\vert-d$ dimensions, we can integrate out 
the $\vert\Omega\vert-d-k$ degrees of freedom associated with the  kernel of nested linear families. The resulting marginal density is again Gaussian, now in $k$ degrees of freedom spanning the orthogonal complement of the kernel of nested linear family $\tilde{\mathcal L}$ within the kernel of the ambient linear family $\hat{\mathcal L}$. 
The statistic of the minimized \textsc{kl} divergence from $\hat{\prob q}$ that appears in Eq.~\eqref{eq:nestedKL} then follows from the standard property of Gaussian marginalization (Appendix~\ref{app:calculus}):  the resulting quadratic form is given by the minimum of the original quadratic form Eq.~\eqref{eq:SamplingProbability:QuadraticForm}. The latter minimization problem reproduces the $I$-projection of $\hat{\prob q}$ onto the nested family. 

From Sections~\ref{ssc:LaplaceApproximation} and~\ref{ssc:pValue_SphericalSymmetry}, this ``Gaussianization''  allows the computation (within Laplace approximation) of the $p$-value associated with an experimental effect encoded in $\tilde{\boldsymbol\mu}$, under structural generalized moments $\hat{\boldsymbol\mu}$, via a $\chi^2_k$ distribution.  
Most notably, this $\chi^2$ approximation to the $p$-value holds, as we shall rigorously show, even in the regime $\vert\Omega\vert=\order{N}$. 
This is rather surprising, since the classical Sanov theory could break down in this regime~\cite{dembo1998large,jiao2015minimax} --\,as already reflected in the failure of the ``naive'' continuum limit Eq.~\eqref{eq:NaiveContinuumLimit} to approximate the volume of the lattice region $\mathbb N_0^{\vert\Omega\vert}\cap N\mathcal P$. 

In particular, as the $\Omega$ space grows larger,  an increasing fraction of datasets $N\prob p$ concentrates on the faces of the target polytope corresponding to vanishing frequencies. Consequently many counts satisfy $n(x)\leq 1$, implying that $\Delta p(x)/p(x)=\order{1}$. The necessary hierarchy $\Delta p(x)\ll p(x)$, required for a well-defined continuum approximation of $p(x)=n(x)/N$, is therefore no longer guaranteed.
In spite of the failure of the entropic (Sanov-Jaynes) description to capture the sampling probabilities of datasets in the double-scaling limit, Gaussian fluctuations in the $k$ directions orthogonal to the nested kernel remain well-behaved. In turn, this enables reliable large-$N$ estimates (via Laplace approximation) for sampling probabilities in the space of generalized moments.

\section{Systematic evaluation of probability mass densities} 
\label{sc:FourierAnalysis}

Starting from the characteristic function of the probability mass of linear dataset families in Section~\ref{ssc:CharacteristicFunctions}, we rigorously define its inverse Fourier transform in Section~\ref{ssc:InverseFourier}. 
Within the continuum approximation (well justified for $D\ll N$), the latter admits the interpretation of the probability mass density over the generalized moments defining linear families.
After analytic continuation, we perform a detailed saddle-point analysis in Section~\ref{ssc:SaddlePointApproximation} to extract the leading-order contribution in $N$ to this probability mass density. 

The rigorous treatment of this section establishes the validity of the ``Gaussianization'' outlined in Section~\ref{ssc:NaiveGaussianization}, extending its applicability to the double-scaling limit $\order{\vert\Omega\vert}=N$, when both the sample size and the number of observable categories grow large. This scaling regime is relevant for high-dimensional data structures and especially for metric attributes $X$ with  theoretically infinite domain $\vert\mathcal X\vert=\infty$.

\subsection{The characteristic function of linear families of datasets}
\label{ssc:CharacteristicFunctions}

It is first didactic to note that the Fourier transform of the multinomial sampling probability Eq.~\eqref{eq:MultinomialDistribution} is given by 
\begin{align}
\label{eq:CharacteristicFunction:Multinomial}
    \widetilde{\Phi}_{\mathcal P}({\mathbf t}) = \sum_{\mathbf n \in\mathbb N_0^{\vert\Omega\vert}} \delta\!\left(\sum_{x\in\Omega}n(x), N\right)\operatorname{mult}(\mathbf n; \prob v) \exp\left\{i {\mathbf t}\cdot\mathbf n\right\}
    =
    \left(\sum_{x\in\Omega} v(x) \exp\left\{i  t(x)\right\}\right)^N~,
\end{align}
parametrized by column vector ${\mathbf t}\in\mathbb R^{\vert\Omega\vert}$, the formal dual to probability vector $\prob p\in\mathcal P$.
As usual, the Kronecker delta constrains the summation over the non-negative integer orthant to datasets of fixed $N$. 
Viewed as a function of real-valued vector $\mathbf t$, the Fourier transform $\widetilde{\Phi}_{\mathcal P}$ is also known as the {characteristic function} of the multinomial distribution. 
The inverse Fourier transform works by applying the Multinomial theorem in the other direction and using the identity for complex phases 
\begin{equation}
\label{eq:Fourier:Kronecker_IntegralRepresentation}
    \int_{-\pi}^\pi \frac{d t}{2\pi} \exp\left\{i t\left(n - m \right)\right\} =
    \begin{cases}
        1 & n=m\in\mathbb N_0\\
        0& n \neq m\in\mathbb N_0\\
    \end{cases}~.
\end{equation}
This is nothing but the integral representation of Kronecker delta $\delta_{n,m}$.

Even more generally, we could start from the characteristic function of a collection of independent Poissonian counts. Recall that for a single Poisson point process 
with expected number of events given by $\lambda>0$, the sampling probability is 
\begin{equation}
   \operatorname{Poisson}(n;\lambda) =  \frac{e^{-\lambda}\lambda^n}{n!} \quad\text{for}\quad n\in\mathbb N_0~.
\end{equation}
Its characteristic function of dual parameter $t\in\mathbb R$ reads 
\begin{equation}
\label{eq:CharacteristicFunction:Poisson}
    \sum_{n\in\mathbb N_0} \operatorname{Poisson}(n;\lambda)\, e^{i n t} 
    = \exp\left\{\lambda e^{it} - \lambda\right\}~.
\end{equation}
Requiring the expected counts in Eq.~\eqref{eq:Intro:expectedCountNumber} to reflect the probabilities of reference distribution $\prob v\in\mathcal P$, the Poissonian sampling probability of a count vector $\mathbf n\in\mathbb N_0^{\vert\Omega\vert}$ 
simply equals the unconstrained product 
\begin{equation}
    \operatorname{Pr}(\mathbf n) = \prod_{x\in\Omega} \operatorname{Poisson}(n(x); Nv(x))~.
\end{equation}
Similar to Eq.~\eqref{eq:CharacteristicFunction:Multinomial}, it is straight-forward to Fourier-transform using the dual vector $\mathbf t\in\mathbb R^{\vert\Omega\vert}$ to obtain an unconstrained characteristic function of the non-negative integer orthant: 
\begin{align}
\label{eq:CharacteristicFunction:Poissons}
    \widetilde{\Phi}(\mathbf t) =&\,\, 
    \sum_{\mathbf n\in \mathbb N_0^{\vert\Omega\vert}} \operatorname{Pr}(\mathbf n)\, \exp\left\{i\,\mathbf t\cdot \mathbf n\right\}
    \nonumber
    \\[1ex]
     =&\,\, 
     \prod_{x\in\Omega} \sum_{n(x)\in\mathbb N_0}\operatorname{Poisson}(n(x);Nv(x))\, \exp\left\{i\, t(x) n(x) \right\}
    \nonumber
    \\[1ex]
    =&\,\,
    \exp\left\{N\left(\sum_{x\in\Omega}v(x)e^{it(x)} - 1\right)\right\} 
    ~.
\end{align}

To restrict this characteristic function to the space of admissible raw means Eq.~\eqref{eq:Intro:RawMean}, we proceed with a projection-slice theorem. This relates the Fourier transform of the non-negative integer orthant over $\Omega$ established in Eq.~\eqref{eq:CharacteristicFunction:Poissons} to the Fourier transform of Eq.~\eqref{eq:SamplingProbabilityMass:Definition}. Crucially, this reduction to linear dataset families defined by coefficient matrix $\mathbf G\in\mathbb R^{D\times\vert\Omega\vert}$ is carried out entirely in Fourier space, irrespective of the scaling of $N, \vert\Omega\vert$ or $D$.

For a vector ${\boldsymbol s}\in\mathbb R^D$ dual to $\boldsymbol m\in\mathbf G\:\mathbb N_0^{\vert\Omega\vert}$, it is thus straight-forward to define the most general characteristic function through
\begin{align}
\label{eq:FourierTrafo:ProbabilityMass}
    \widetilde{\Phi}_{\Lambda}({\boldsymbol s}) =&\,\,
    \sum_{\boldsymbol m\in\mathbf G\mathbb N^{\vert\Omega\vert}_0}
    \operatorname{Pr}\left(\boldsymbol m\right)\exp\left\{i\,\boldsymbol s\cdot \boldsymbol m\right\}
    \\[1ex]
    =&\,\,
    \sum_{\boldsymbol m\in\mathbf G\mathbb N^{\vert\Omega\vert}}
    \sum_{\mathbf n\in\mathbb N_0^{\vert\Omega\vert}} 
    \delta(\mathbf n\in \Lambda_{\boldsymbol{m}})\,\operatorname{Pr}(\mathbf n)
    \exp\left\{i \,\boldsymbol s^T \mathbf G\,\mathbf n\right\}
    =
    \widetilde{\Phi}
    ({\boldsymbol s}^T\mathbf G)
    \nonumber
    ~.
\end{align}
In the second line, we substituted $\boldsymbol m=\mathbf G\:\mathbf n$ in the Fourier exponent for all datasets in the respective generalized Diophantine set $\Lambda_{\boldsymbol m}$ of Eq.~\eqref{eq:DatasetFamily}. The double summation --\,over admissible generalized moments 
and over member datasets in each linear family\,-- amounts to an overall summation over the whole semi-lattice $\mathbb N_0^{\vert\Omega\vert}$. 
Therefore, the Fourier transform of the coarse-grained probability mass of our families is obtained by evaluating the original Fourier transform in Eq.~\eqref{eq:CharacteristicFunction:Poissons} at ${\mathbf t} = {\boldsymbol s}^T\mathbf G$.

\subsection{The inverse Fourier transform}
\label{ssc:InverseFourier}

In anticipation of the large-$N$ expansion (see in particular Sections~\ref{ssc:VariationScale} and~\ref{ssc:LinearAlgebra}), we parametrize raw means by generalized moments $\boldsymbol{\mu}=\boldsymbol{m}/N$ for some yet to be fixed scale $N$. 
To recover from Eq.~\eqref{eq:FourierTrafo:ProbabilityMass} the coarse-grained probability mass of the linear family of datasets supporting generalized moments $\boldsymbol{\mu}$,
\begin{equation}
    \operatorname{Pr}(N\boldsymbol\mu) 
    = \sum_{\mathbf n\in\mathbb N_0^{\vert\Omega\vert}}
    \operatorname{Pr}(\mathbf n)\,
    \delta\!\left(\mathbf n\in N\mathcal L_{\boldsymbol\mu}\right)~,
\end{equation}
we need to perform the inverse Fourier transformation in a non-canonical setting, where $\boldsymbol m=N\boldsymbol\mu$ are neither generally integer-valued nor equally spaced. 

To do so, we first generalize the integral representation of Kronecker delta Eq.~\eqref{eq:Fourier:Kronecker_IntegralRepresentation} to force the $\alpha$-th generalized moment to take the value $\mu_\alpha$ via a Dirac delta constraint, also signified by $\delta$. Briefly dropping in this paragraph the moment index for clarity, we note the formal limit in the one-dimensional case 
\begin{equation}
\label{eq:DiracDeltaFunction:FormalLimit}
    \lim_{\varepsilon\rightarrow0^+} \int_{-\infty}^{\infty} ds \exp\left\{i N\left(\tilde\mu-\mu\right) s - 
    \varepsilon\, s^2\right\} = 
    \begin{cases}
        \infty & \tilde\mu = \mu \\
        0 & \text{else} 
    \end{cases}~.
\end{equation}
Among many possible representations, we have selected one which introduces a natural regularizer through the quadratic term in $s\in\mathbb R$ ---\,a feature that facilitates the subsequent analytic continuation.
The overall normalization can be fixed by requiring consistency in the space of possible generalized moments, namely
\begin{equation}
    \sum_{N\boldsymbol\mu\in \mathbf G\mathbb N_0} \delta\left(N\tilde\mu-N\mu\right) = 1~.
\end{equation}
Approximating the sum via an integral, we can determine up to $\varepsilon$-smooth terms
\begin{equation}
    \lim_{\varepsilon\rightarrow0^+} \int_{\mu_\text{min}}^{\mu_\text{max}}d\mu\, C \int_{-\infty}^{\infty} ds\exp\left\{i N\left(\tilde\mu-\mu\right) s -
    \varepsilon\, s^2\right\} \overset{!}{=} 1 \quad\Rightarrow\quad C = \frac{N}{2\pi}
    ~.
\end{equation}

Such a continuous approximation remains valid at least as long as $D\ll \vert\Omega\vert, N$. 
In the $D$-dimensional space, the replacement of sums over admissible values of the generalized moments by integrals constitutes the natural continuum analogue of the multidimensional Euler–Maclaurin approximation for lattice sums.
The leading corrections are expected to arise from the boundary $\partial\mathcal L$ of the corresponding polytope. In fact, the saddle-point analysis of the density derived below reveals a stronger, exponentially small suppression of these corrections, whenever the $I$-projection of $\prob v$ onto $\mathcal L$ remains in the interior of the simplex.

In total, the appropriately normalized dimensionless differential in the continuum approximation of our inverse Fourier transform must be 
\begin{equation}
    \left(\frac{N}{2\pi}\right)^D d^D\!\boldsymbol\mu\, d^D\!\boldsymbol s \equiv \prod_{\alpha=0}^{D-1} d(N\mu_\alpha) \frac{ds_\alpha}{2\pi}~.
\end{equation}
In particular, the integral representation of Dirac delta constraint in the space of generalized moments reads 
\begin{equation}
\label{eq:DiracDeltaFunction:IntegralRepresentation}
    \boldsymbol\delta^{(D)}\!\left(N\tilde{\boldsymbol\mu}-N\boldsymbol\mu\right) \equiv N^D \left[\prod_{\alpha=0}^{D-1} \int_{-\infty}^\infty \frac{d s_\alpha}{2\pi}\right] \exp\left\{i N\left(\tilde\mu_\alpha-\mu_\alpha\right) s_\alpha - 
    \varepsilon_\alpha\, s^2_\alpha\right\}~,
\end{equation}
with the (silent) understanding that the $\varepsilon_\alpha\rightarrow0^+$ limit is to be taken at the very end. 
The dimension of the Dirac delta function is given by the dimension of volume in the space of inverse generalized moments, 
\begin{equation}
    \left[\boldsymbol\delta^{(D)}\right]=\prod_{\alpha=0}^{D-1}\left[\mu_\alpha^{-1}\right]~.
\end{equation}

\paragraph{The classical action.}
Given a target effect encoded by generalized moments $\boldsymbol\mu\in\mathbb R^D$,  
we formally retrieve the probability mass density corresponding to the family of compatible datasets
from characteristic function Eq.~\eqref{eq:FourierTrafo:ProbabilityMass}  by employing  prescription Eq.~\eqref{eq:DiracDeltaFunction:IntegralRepresentation}:
\begin{equation}
\label{eq:InverseFourier:ProbabilityMassDensity}
    \rho(\boldsymbol\mu) 
    = 
    \left(\frac{N}{2\pi}\right)^D 
    \int d^D\!\boldsymbol{s}\,
    \exp\left\{-i N\,\boldsymbol s\cdot \boldsymbol \mu\right\}
    \, \widetilde{\Phi}_{\Lambda}(\boldsymbol s)
    \equiv 
   \left(\frac{N}{2\pi}\right)^D \int d^D\!\boldsymbol{s}\,
   \, \exp\left\{-N W(\boldsymbol s)\right\}
   ~,
\end{equation}
with the resulting action of the classical problem being
\begin{equation}
\label{eq:InverseFourier:Action}
    W(\boldsymbol{s})  = \sum_{\alpha=0}^{D-1} \left[\frac{\varepsilon_\alpha\,s_\alpha^2}{N}  + i s_\alpha \mu_\alpha\right] 
    -\sum_{x\in\Omega}v(x)\exp\left\{i\sum_{\alpha=0}^{D-1}s_\alpha  g_\alpha(x) \right\} + 1
     ~.
\end{equation}
When trivially integrating over the full real line, we drop for brevity  the limits of integration.
From the formal limit in the prescription Eq.~\eqref{eq:DiracDeltaFunction:FormalLimit}  of the integral representation of Dirac delta function, 
it follows that $\order{\Vert\boldsymbol\varepsilon\Vert/N} \ll 1$.

\subsection{Saddle-point approximation}
\label{ssc:SaddlePointApproximation}

The overall $N$ in the exponent of the integrand in Eq.~\eqref{eq:InverseFourier:ProbabilityMassDensity} invites us to apply large-$N$ techniques to evaluate the density associated with the coarse-grained  probability mass of the linear family. 
For that, we proceed with an analytic continuation of the classical action $W$ by complexifying the dual variables to
\begin{equation}
    \boldsymbol z\in\mathbb C^D \quad\text{with}\quad \operatorname{Re}(\boldsymbol{z}) = \boldsymbol{s} 
    ~.
\end{equation}
The saddle-point equations for $W(\boldsymbol{z})\in\mathbb C$, 
\begin{equation}
\label{eq:SaddlePoint:Definition}
    \left.\nabla_{\boldsymbol z} W\right\vert_{\boldsymbol{z}=\boldsymbol{z}^*} \overset{!}{=} \mathbf 0~,
\end{equation}
produce a set of $D$ conditions  
\begin{equation}
\label{eq:InverseFourier:SaddlePoint}
    \mu_\alpha = \sum_{x\in\Omega} v(x) \exp\left\{i\sum_{\beta=0}^{D-1} z_\beta^*\,  g_\beta(x) \right\} g_\alpha(x) + 
    \order{\Vert\boldsymbol\varepsilon\Vert}~.
\end{equation}

Real-valued expectations $\mu_\alpha=\mathbb E[g_\alpha(X)]\in\mathbb R$ $\forall\,\alpha=0,\ldots,D-1$ imply that the saddle point of interest lies in the vicinity of the imaginary axis, so that we may write 
\begin{equation}
    \lim_{\boldsymbol\varepsilon\rightarrow\mathbf 0^+} \boldsymbol z^* = -i \boldsymbol\theta 
   ~.
\end{equation}
Using the $\boldsymbol\theta$-parametrized distribution  (by Eq.~\eqref{eq:NormalizationSymmetry} normalization is always implied)
\begin{equation}
\label{eq:Fourier:SaddlePoint_Iprojection}
     p(x;\boldsymbol\theta) = 
     v(x) \exp\left\{ \sum_{\alpha=0}^{D-1}\theta_\alpha g_\alpha(x) \right\} 
\end{equation}
from the exponential family $\mathcal E(\mathbf G;\prob v)$ in Eq.~\eqref{eq:ExponentialFamily} generated around $\prob v\in\mathcal P$ by coefficient matrix $\mathbf G\in\mathbb R^{D\times\vert\Omega\vert}$, we recognize that Eq.~\eqref{eq:InverseFourier:SaddlePoint} predominantly reproduces the familiar conditions Eq.~\eqref{eq:LinearSystemDEF} on generalized moments: 
\begin{equation}
\label{eq:SaddlePoint:GeneralizedMomentConditions}
    \alpha=0,\ldots,D-1\,:\quad \sum_{x\in\Omega} p(x;\boldsymbol{\theta}) g_\alpha(x) + \order{\Vert\boldsymbol\varepsilon\Vert} = \mu_\alpha~.
\end{equation}
As discussed around Eq.~\ref{eq:app:iProjection:Reparametrization}, a distribution 
of this exponential form that satisfies\footnote{Removal of the regularizer kills the real offset in $\mathbf z^*$ as well as any higher corrections in $\boldsymbol\theta$.} the defining conditions in Eq.~\eqref{eq:SaddlePoint:GeneralizedMomentConditions} of linear family $\mathcal L_{\boldsymbol{\mu}}\equiv\mathcal L(\mathbf G;\boldsymbol\mu)$ coincides with the unique $I$-projection of reference  $\prob v$ onto said linear family.
We can thus write 
\begin{equation}
\label{eq:Fourier:Iprojection_identification}
    \prob p_{\boldsymbol\theta}=\argmin_{\prob p\in\mathcal L_{\boldsymbol\mu}}\infdiv{\prob p}{\prob v}\equiv \prob q~.
\end{equation} 
The dual parameters $\boldsymbol\theta\equiv\boldsymbol\theta(\boldsymbol{\mu})\in\mathbb R^D$ remain finite due to the implications of Eq.~\eqref{eq:NoOutliers}. 

We emphasize that the forward $I$-projection naturally emerges as the distribution parameterizing the saddle-point conditions in Eq.~\eqref{eq:SaddlePoint:Definition}, without ever invoking concepts from information geometry or even introducing probability distributions \emph{ad hoc}. In this sense, the information geometry of Csiszár and Amari~\cite{amari2000methods} is not assumed \emph{a priori} through a statistical-manifold formulation. Rather, in the original spirit of Boltzmann's work, it emerges as an intrinsic consequence of the large-$N$ description of independent sampling, with the same geometric structure extending seamlessly to the space of generalized moments.

\paragraph{Contour deformation and saddle-point evaluation.}
Because our regulator in Eq.~\eqref{eq:InverseFourier:Action} is separable and gives absolute convergence, we may repeatedly apply  Cauchy's integral theorem in one complex variable, deforming the integration contour successively in each component of $\boldsymbol{z}$ while keeping the remaining fixed.  
Holding the remaining $z_\beta$ with $\beta\neq\alpha$ fixed, we can specifically apply Cauchy's theorem to the closed contour $\mathcal C_\alpha$ in the complex plane of the $\alpha$-th dual variable:
\begin{equation}
\label{eq:InverseFourier:CauchyTheorem}
    \int_{\mathcal C_\alpha}d z_\alpha\, \exp\left\{-N W(\boldsymbol{z})\right\} = 0~, 
\end{equation}
where $W(z_0,\ldots,z_\alpha,\ldots,z_{D-1})$ is entire\footnote{The exponential of an entire function is also entire}. 
Along the integration contour,
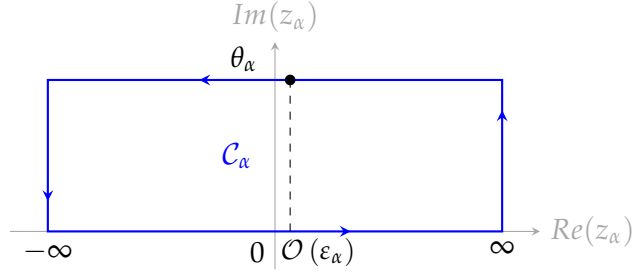
\begin{center}
\begin{tikzpicture}[>=stealth]
    \def\L{3}
    \def\thetaVal{2}
    \def\off{0.5}

    \draw[->, gray!70] (-\L-\off, 0) -- (\L+\off, 0) node[right] {$Re(z_\alpha)$};
    \draw[->, gray!70] (0, -0.5) -- (0, \thetaVal+\off) node[above] {$Im(z_\alpha)$};

    \draw[thick, blue, 
        postaction={decorate, decoration={
            markings,
            mark=at position 0.100 with {\arrow{>}},
            mark=at position 0.375 with {\arrow{>}},
            mark=at position 0.600 with {\arrow{>}},
            mark=at position 0.875 with {\arrow{>}}
        }}
    ] (\L, 0) -- (\L, \thetaVal) -- (-\L, \thetaVal) -- (-\L, 0) -- cycle;

    \node[below] at (\L, 0) {$\infty$};
    \node[below] at (-\L, 0) {$-\infty$};
    \node[left] at (-0.1, 0.28+\thetaVal) {$\theta_\alpha$};
    \fill[black] (0.2,\thetaVal) circle (2pt);   
    \node[left] at (1.2, -0.27) {\small$\order{\varepsilon_\alpha}$};
    \draw[dashed] (0.2,0) -- (0.2,\thetaVal);  
    \node[below left] at (0,0) {$0$};
    
    \node[blue] at (-0.5, \thetaVal/2) {$\mathcal{C}_\alpha$};
\end{tikzpicture}
\end{center}
the integral decomposes into two horizontal (parallel to the real axis) and two vertical (parallel to the imaginary axis) integrals. The contribution from the latter vertical contour segments vanishes, because the integrand decays exponentially as $\operatorname{Re}(z_\alpha)\rightarrow\infty$, owing to the regularizing $s_\alpha^2$ term in the classical action of Eq.~\eqref{eq:InverseFourier:Action}.  

Consequently, after taking the formal $\varepsilon_\alpha\rightarrow0^+$ limit, Eq.~\eqref{eq:InverseFourier:CauchyTheorem} for $\alpha=0$ reduces  to 
\begin{equation}
    \int_{-\infty}^\infty  d s_0\, \exp\left\{-NW(s_0,z_1,\ldots,z_{d-1})\right\} = \int_{-\infty}^\infty  d s_0\, \exp\left\{-NW(s_0-i\theta_0,z_1,\ldots, z_{d-1})\right\} ~.
\end{equation}
Applying Cauchy's theorem in the remaining $D-1$ directions similarly deforms the original multiple integral in Eq.~\eqref{eq:InverseFourier:ProbabilityMassDensity}, shifting the integration contour across the imaginary axis so that it passes through the saddle point of Eq.~\eqref{eq:InverseFourier:SaddlePoint}:
\begin{equation}
\label{eq:InverseFourier:sIntegral_over_density}
    \rho(\boldsymbol\mu) = 
   \left(\frac{N}{2\pi}\right)^{D}
   \int d^D\!\boldsymbol{s}
    \exp\left\{- N W\left(\boldsymbol{s}-i\boldsymbol\theta\right)\right\}~.
\end{equation}

After deforming the integration contour so that it passes through the saddle point, we expand the action in a Taylor series about $-i\boldsymbol\theta$ along the real directions, 
\begin{align}
    W\left(\boldsymbol s - i\boldsymbol{\theta}\right) \approx&\,\, W\left(-i\boldsymbol{\theta}\right) + \tfrac12\sum_{\alpha,\beta=0}^{D-1} 
    \left.\frac{\partial^2 W\left(\boldsymbol s - i\boldsymbol{\theta}\right)}{\partial{s}_\alpha\partial{s}_\beta}\right\vert_{\boldsymbol s=\mathbf0}
    s_\alpha s_\beta ~.
\end{align}
The zeroth-order term evaluates to the minimized \textsc{kl} divergence from reference distribution $\prob v$,
\begin{equation}
\label{eq:SaddlePoint:MinimizedKL}
    W\left(-i\boldsymbol\theta\right) = \boldsymbol\theta \cdot\boldsymbol{\mu} 
    =
    \min_{\prob p\in\mathcal L_{\boldsymbol{\mu}}} \infdiv{\prob p}{\prob v} 
    ~,
\end{equation} 
using the $\boldsymbol{\theta}$-parametrization of the $I$-projection Eq.~\eqref{eq:Fourier:SaddlePoint_Iprojection}.
Owing to Eq.~\eqref{eq:SaddlePoint:GeneralizedMomentConditions} after taking the $\Vert\boldsymbol\varepsilon\Vert\rightarrow0$ limit,
\begin{equation}
\label{eq:SaddlePoint:GeneralizedMoment_Parameter_Duality}
    \sum_{x\in\Omega} p(x;\boldsymbol{\theta}) g_\alpha(x)  = \mu_\alpha \quad\text{for}\quad\prob p_{\boldsymbol\theta}\in\mathcal E(\mathbf G;\prob v)~,
\end{equation}
there is a generically complicated functional dependence $\boldsymbol{\theta}\equiv\boldsymbol{\theta}(\boldsymbol{\mu})$ of dual parameters on the generalized moments.
The Hessian coincides with the Jacobian derived in Eq.~\eqref{eq:iProjector:Jacobian},
\begin{equation}
\label{eq:SaddlePoint:Jacobian}
    J_{\alpha\beta} \equiv \left.\frac{\partial^2 W\left(\boldsymbol s - i\boldsymbol{\theta}\right)}{\partial{s}_\alpha\partial{s}_\beta}\right\vert_{\boldsymbol s=\mathbf0} = \sum_{x\in\Omega} g_\alpha(x) p(x;\boldsymbol\theta) g_\beta(x)~,
\end{equation}
after using the distribution from Eq.~\eqref{eq:Fourier:Iprojection_identification}
As shown in Section~\ref{ssc:iProjector}, this matrix is positive definite.

Finally, evaluating the multiple integral over $\boldsymbol s$ in Eq.~\eqref{eq:InverseFourier:sIntegral_over_density} around the saddle point, where the action is well approximated by its quadratic expansion for large $N$, we retrieve the coarse-grained probability mass associated with the differential region of linear family $\mathcal L(\mathbf G;\boldsymbol{\mu})$: 
\begin{equation}
\label{eq:ProbabilityMass_Differential}
    d^D\!\boldsymbol{\mu}\,\rho(\boldsymbol\mu) = 
    d^D\!\boldsymbol\mu\,
    \sqrt{\det\left(\frac{N  \mathbf J^{-1}(\boldsymbol\mu)}{2\pi}\right)}
    \exp\left\{- N \min_{\prob p\in\mathcal L_{\boldsymbol\mu}}\infdiv{\prob p}{\prob v}
    \right\} 
    ~.
\end{equation}
Therefore, the sampling probability of the linear families (and hence of all compatible datasets) in the vicinity of $\mathcal L(\mathbf G;\boldsymbol{\mu})$ is governed to leading order in $N$ by the minimized \textsc{kl} divergence in Eq.~\eqref{eq:SaddlePoint:MinimizedKL}, with a subleading determinant prefactor accounting for volume corrections. 
Here, as in the subsequent Laplace approximation, we disregard $1/N$-suppressed effects, since these are non-universal and depend on the the precise details of the sampling process at hand. 

Using the generalized-moment-parameter duality in Eq.~\eqref{eq:SaddlePoint:GeneralizedMoment_Parameter_Duality}, we may change variables from the generalized moments to the dual parameters $\boldsymbol{\theta}\in\mathbb R^D$,
\begin{equation}
    \boldsymbol\mu = \mathbf G\: \prob p_{\boldsymbol\theta} \quad\rightarrow\quad 
    d^D\!\boldsymbol\mu = \det\mathbf J \,d^D\! \boldsymbol\theta ~,
\end{equation}
so that the sampling probability associated with the differential element in dual space becomes
\begin{equation}
\label{eq:ProbabilityMass_Differential:DualSpace}
    d^D\!\boldsymbol\theta\,\rho(\boldsymbol\theta) = 
    d^D\!\boldsymbol\theta\,
    \sqrt{\det\left(\frac{N  \mathbf J(\boldsymbol\theta)}{2\pi}\right)}
    \exp\left\{- N \infdiv{\prob p_{\boldsymbol\theta}}{\prob v}
    \right\} 
    ~,
\end{equation}
with $\mathbf G\:\prob p_{\boldsymbol\theta}=\boldsymbol\mu$ for $\prob p_{\boldsymbol\theta}\in\mathcal E(\mathbf G;\prob v)$.
In the regime $D\ll N$, expressions Eq.~\eqref{eq:ProbabilityMass_Differential} and~\eqref{eq:ProbabilityMass_Differential:DualSpace} provide the continuous approximation to the sampling probability mass $\operatorname{Pr}(\boldsymbol{m})$ of the linear family $\Lambda_{\boldsymbol{m}}$ consisting of all datasets $\mathbf n\in\mathbb N_0^{\vert\Omega\vert}$ that exhibit the $D$ raw means ${m}_\alpha$. No restriction is imposed on the cardinality of $\Omega$, which may itself grow with $N$.

\section{The large-\texorpdfstring{$N$}{N} approximation to the \texorpdfstring{$p$}{p}-value}
\label{sc:pValue}

In this section, the derived expressions for the probability mass density in the $\boldsymbol\mu$- and $\boldsymbol{\theta}$-dual spaces are leveraged to approximate the probability mass of nested dataset families under conditional sampling. Working within the continuous approximation inherited from the preceding section, we parametrize the problem in Section~\ref{ssc:ConditionalSampling} using the algebraic framework of linear families of distributions. 
By a slight overload of terminology, we then discuss in Section~\ref{ssc:SamplingProbabilityDensity} the conditional sampling probability density of such linear families of distributions. 

In the continuum limit, these linear families provide an effective parametrization of neighboring families of datasets, whose coarse-grained probability mass is approximated by the Laplace method  employed in Section~\ref{ssc:LaplaceApproximation}. 
Finally, we demonstrate in Section~\ref{ssc:pValue_SphericalSymmetry}  how the sampling-based definition of the $p$-value in Eq.~\eqref{eq:Intro:pValue}, alongside its implied spherical symmetry, leads in Laplace approximation to the semi-analytic, information-theoretic closed form of Eq.~\eqref{eq:LaplaceApproximation:CriticalValue}.

\subsection{The Linear Algebra of Conditional Sampling}
\label{ssc:ConditionalSampling}

We first extend the linear-algebraic language of Section~\ref{ssc:LinearAlgebra} to account for conditional sampling. 
By keeping $d<D$ generalized moments fixed, 
\begin{equation}
\label{eq:StructuralConditions}
    \widebar{g_\alpha(X)}\overset{!}{=}\hat\mu_\alpha \quad\text{for}\quad \alpha=0,\ldots,d-1~,
\end{equation}
we define a {structural} linear family of distributions 
\begin{equation}
\label{eq:StructuralFamily}
    \widehat{\mathcal L}\equiv\mathcal L(\mathbf G_S;\hat{\mu}_0,\ldots,\hat{\mu}_{d-1}) = \left\{\prob p\in\mathcal P ~\vert~ \sum_{x\in\Omega}g_\alpha(x)p(x)=\hat\mu_\alpha \quad\text{for}\quad \alpha=0,\ldots,d-1\right\}~,
\end{equation}
by using coefficient matrix $\mathbf G_S\in\mathbb R^{d\times\vert\Omega\vert}$.
For $k>0$, any linear family constructed by $\mathbf G\in\mathbb R^{(d+k)\times\vert\Omega\vert}$ such that 
\begin{equation}
\label{eq:ExperimentalFamily}
    \mathcal L_{\boldsymbol{\mu}} \equiv \mathcal L(\mathbf G;\boldsymbol{\mu})= \left\{\prob p\in\widehat{\mathcal L} ~\vert~ 
    \sum_{x\in\Omega}g_{d+i}(x)p(x)=\mu_{d+i} \quad\text{for}\quad i=0,\ldots,k-1\right\}
\end{equation}
will be nested within the structural family $\widehat{\mathcal L}$ in the usual linear-algebraic sense, i.e.\ $\mathcal L_{\boldsymbol{\mu}}\subset\widehat{\mathcal L}\subseteq\mathcal P$, by imposing experimental-type conditions on $k$ additional sample means, 
\begin{equation}
\label{eq:ExperimentalConditions}
    \widebar{g_{d+i}(X)}\overset{!}{=}\mu_{d+i} \quad\text{for}\quad i=0,\ldots,k-1~.
\end{equation}

Accordingly, the coefficient matrix of the nested description with $D=d+k$ in Eq.~\eqref{eq:intro:CoeffientMatrix} decomposes into 
\begin{equation}
    \mathbf G = \begin{pmatrix}
        \mathbf G_S \\
        \mathbf G_E 
    \end{pmatrix}
    \quad\text{with}\quad \mathbf G_S\in\mathbb R^{d\times \vert\Omega\vert} \quad\text{and}\quad \mathbf G_E\in\mathbb R^{k\times \vert\Omega\vert}~,
\end{equation}
so that the Jacobian enjoys a block form 
\begin{equation}
\label{eq:Jacobian_BlockForm}
    \mathbf J = 
    \begin{pmatrix}
        \mathbf J_S & \mathbf M\\
        \mathbf M^T & \mathbf J_E\\
    \end{pmatrix}~.
\end{equation}
Due to the assumptions of Section~\ref{app:IdealProblem}, the block Jacobian matrices $\mathbf J_S\in\mathbb R^{d\times d}$ and $\mathbf J_E\in\mathbb R^{k\times k}$ 
are positive definite:
\begin{equation}
    (\mathbf J_S)_{\alpha\beta} = \sum_{x\in\Omega} g_\alpha(x) \,p(x)\, g_\beta(x)
    \quad\text{and}\quad
    (\mathbf J_E)_{ij} = \sum_{x\in\Omega} g_{d+i}(x) \,p(x)\, g_{d+j}(x)
\end{equation}
The mixing matrix is defined by
\begin{equation}
    M_{\alpha,i} = \sum_{x\in\Omega} g_\alpha(x)\, p(x)\, g_{d+i}(x)\quad\text{for}\quad \alpha=0,\ldots,d-1 \quad\text{and}\quad i=0,\ldots,k-1~.
\end{equation}
Arranging that the condition for $\alpha=0$ reflects normalization (i.e.\ $g_0$ is the constant map), we have $J_{\alpha,0}=J_{0,\alpha}=\mu_\alpha$ for any $\alpha=0,\ldots,d+k-1$ in the full $(d+k)\times(d+k)$ space, meaning $(\mathbf J_S)_{0,\alpha}=\hat\mu_\alpha$ for $\alpha=0,\ldots,d-1$ and $M_{0,i}=\mu_{d+i}$ for $i=0,\ldots,k-1$.
Using Schur's complement based on block form of Jacobian in  Eq.~\eqref{eq:Jacobian_BlockForm} 
\begin{equation}
\label{eq:SchurComplement}
    \mathbf \Sigma^{-1} \equiv  \mathbf J_E - \mathbf M^T \mathbf J_S^{-1} \mathbf M~,
\end{equation}
we note the determinant relation
\begin{equation}
\label{eq:SchurComplement:detIdentity}
    \det \mathbf J = \det \mathbf J_S \cdot \det\left( \mathbf \Sigma^{-1} \right)
    ~.
\end{equation}

\subsection{The probability mass density of nested families}
\label{ssc:SamplingProbabilityDensity}

Following the notational convention of Eq.~\eqref{eq:Intro:ConditionalSampling_ProbabilityMass} for the probability mass under conditional sampling, we introduce the conditional probability mass density
\begin{equation}
\label{eq:ConditionalProbabilityMassDensity:Definition}
    \rho(\mu_{d},\ldots,\mu_{d+k-1}\vert\hat\mu_0,\ldots,\hat\mu_{d-1}) 
    = 
    \frac{\rho(\hat\mu_0,\ldots,\hat\mu_{d-1},\mu_d,\ldots,\mu_{d+k-1})}{\rho(\hat\mu_0,\ldots,\hat\mu_{d-1})}~,
\end{equation}
with the common coarse-graining volume element (associated with the ambient family) canceling in the ratio. 
Employing Eq.~\eqref{eq:ProbabilityMass_Differential} from the saddle-point approximation twice, in the $(d+k)$-dimensional space corresponding to the nested family and in the  $d$-dimensional space corresponding to the structural family, we arrive at 
\begin{equation}
\label{eq:ConditionalProbabilityMassDensity:muSpace}
    \rho(\mu_{d},\ldots,\mu_{d+k-1}\vert\hat\mu_0,\ldots,\hat\mu_{d-1}) 
    = 
    \sqrt{\left(\frac{N}{2\pi}\right)^k \frac{\det \widehat{\mathbf J}}{\det\mathbf J(\boldsymbol\mu)}} \exp\left\{-N \min_{\prob p\in\mathcal L_{\boldsymbol\mu}}\infdiv{\prob p}{\hat{\prob q}}\right\}~.
\end{equation}
Here, we simplified the effective action\footnote{In analogy with statistical mechanics, we formally define the effective action given some density $\rho$ through $-\lim_{N\rightarrow\infty}  N^{-1}\log\rho$.} 
arising from the two contributions of Eq.~\eqref{eq:ProbabilityMass_Differential},
\begin{equation}
\label{eq:nestedDensity:Pythagorean}
    \min_{\prob p\in\mathcal L_{\boldsymbol\mu}}\infdiv{\prob p}{\prob v} - \min_{\prob p\in\widehat{\mathcal L}}\infdiv{\prob p}{\prob v} = \min_{\prob p\in\mathcal L_{\boldsymbol\mu}}\left[\infdiv{\prob p}{\prob v} - \infdiv{\hat{\prob q}}{\prob v}\right] =
    \min_{\prob p\in\mathcal L_{\boldsymbol\mu}}\infdiv{\prob p}{\hat{\prob q}}~,
\end{equation}
using the Pythagorean identity Eq.~\eqref{eq:PythagorasIdentity} in the structural linear family $\widehat{\mathcal L}$, onto which the $I$-projection of $\prob v$, denoted by 
\begin{equation}
\label{eq:ambient_iProjection}
    \hat{\prob q} = \argmin_{\prob p\in\widehat{\mathcal L}} \infdiv{\prob p}{\prob v}~,
\end{equation}
remains fixed throughout the present discussion of conditional sampling. 

Consequently, the Jacobian $\widehat{\mathbf J}\in\mathbb R^{d\times d}$ from Eq.~\eqref{eq:SaddlePoint:Jacobian} corresponding to the structural linear family with 
\begin{equation}
\label{eq:ambientJacobian}
    \widehat J_{\alpha\beta} = \sum_{x\in\Omega} g_\alpha(x)\, \hat q(x)\, g_\beta(x) \quad\text{for}\quad\alpha,\beta = 0,\ldots,d-1~,
\end{equation} 
also remains fixed. In contrast, the Jacobian $\mathbf J\equiv \mathbf J(\boldsymbol{\mu})\in\mathbb R^{(d+k)\times(d+k)}$ corresponding to nested linear families $\mathcal L_{\boldsymbol\mu}$ with 
\begin{equation}
    J_{\alpha\beta} = \sum_{x\in\Omega} g_\alpha(x)\, p(x;\boldsymbol{\theta})\, g_\beta(x) \quad\text{for}\quad\alpha,\beta = 0,\ldots,d+k-1~,
\end{equation}
varies with the specified experimental conditions in Eq.~\eqref{eq:ExperimentalConditions} through the varying $I$-projection
\begin{equation}
\label{eq:Laplace:nestediProjection}
    \prob p_{\boldsymbol\theta}=\argmin_{\prob p\in\mathcal L_{\boldsymbol\mu}}\infdiv{\prob p}{\prob v}~.
\end{equation}

To explicitly account for this variation, we switch to the dual parameter space via Eq.~\eqref{eq:SaddlePoint:GeneralizedMoment_Parameter_Duality}.
In the coordinate parametrization of the $D$-dimensional dual space in Eq.~\eqref{eq:Fourier:SaddlePoint_Iprojection} with $D=d+k$, the only independent degrees of freedom turn out to be the parameters 
\begin{equation}
\label{eq:ExperimentalParameters}
    \theta_{d+i}\equiv\lambda_i \quad\text{for}\quad i=0,\ldots,k-1
\end{equation} 
corresponding to the experimental conditions. We collectively denote these experimental parameters by vector $\boldsymbol{\lambda}\in\mathbb R^k$. 
The structural conditions 
\begin{equation}
\label{eq:lambdaParametrization:StructuralConditions}
    \sum_{x\in\Omega}g_\alpha(x) p(x;\theta_0,\ldots,\theta_{d-1},\boldsymbol{\lambda})=\hat\mu_\alpha \quad\text{for}\quad \alpha=0,\ldots,d-1
\end{equation}
in Eq.~\eqref{eq:StructuralConditions} then implicitly determine  the remaining parameters
\begin{equation}
\label{eq:lambdaParametrization:StructuralParameters}
\theta_\alpha\equiv\theta_\alpha(\boldsymbol\lambda)
\quad\text{for}\quad \alpha=0,\ldots,d-1~.
\end{equation}
Likewise, the experimental generalized moments 
\begin{equation}
\label{eq:lambdaParametrization:ExperimentalMoments}
    \mu_{d+i}\equiv\mu_{d+i}(\boldsymbol{\lambda}) \quad\text{for}\quad i=0,\ldots,k-1
\end{equation}
can be naturally regarded through Eq.~\eqref{eq:ExperimentalConditions}, prescribing 
\begin{equation}
\label{eq:ExperimentalConditions:lambdaParametrization}
    \sum_{x\in\Omega}g_{d+i}(x) p(x;\theta_0(\boldsymbol\lambda),\ldots,\theta_{d-1}(\boldsymbol\lambda),\boldsymbol{\lambda})=\mu_{d+i} \quad\text{for}\quad i=0,\ldots,k-1
    ~,
\end{equation}
as functions of our controlling parameters $\boldsymbol\lambda$.

Hence, $\boldsymbol\lambda$ provides a natural coordinate system on the set of nested linear families, 
succinctly denoted in the following by 
\begin{equation}
\label{eq:nestedFamilies:lambdaParametrization}
\mathcal L_{\boldsymbol{\lambda}}\equiv \mathcal L(\mathbf G;\hat\mu_0,\ldots,\hat\mu_{d-1}, \mu_{d}(\boldsymbol\lambda),\ldots, \mu_{d+k-1}(\boldsymbol\lambda))\subset\widehat{\mathcal L}
~.
\end{equation}
Viewed as the unique intersection of the exponential family $\mathcal E(\mathbf G;\prob v)$ with the $\boldsymbol{\lambda}$-parametrized linear family, the nested $I$-projection \eqref{eq:Laplace:nestediProjection} accordingly admits the parametrization
\begin{equation}
\label{eq:lambdaParametrization:ExponentialFamily}
    p(x;\boldsymbol{\lambda}) = v(x) \exp\left\{\sum_{\alpha=0}^{d-1} \theta_\alpha(\boldsymbol{\lambda}) g_\alpha(x) + \sum_{i=0}^{k-1} \lambda_i\,g_{d+i}(x)\right\}~.
\end{equation}
Therefore, the conditional probability mass density in Eq.~\eqref{eq:ConditionalProbabilityMassDensity:muSpace} associated with a nested linear family $\mathcal L_{\boldsymbol{\lambda}}$ can be re-expressed through Eq.~\eqref{eq:ProbabilityMass_Differential:DualSpace} in the dual parameter space 
\begin{equation}
    \label{eq:ConditionalProbabilityMassDensity:dualSpace}
    \rho(\boldsymbol{\lambda}\vert\hat\mu_0,\ldots,\hat\mu_{d-1}) 
    =
    \sqrt{\left(\frac{N}{2\pi}\right)^k\frac{\det\mathbf J(\boldsymbol{\lambda})}{\det\widehat{\mathbf J}}}
    \exp\left\{-N \infdiv{\prob p_{\boldsymbol{\lambda}}}{\hat{\prob q}}\right\} ~.
\end{equation}

Closing this paragraph, we note some useful linear-algebraic relations governing the variation of structural parameters in Eq.~\eqref{eq:lambdaParametrization:StructuralParameters} and experimental generalized moments in Eq.~\eqref{eq:lambdaParametrization:ExperimentalMoments} with respect to the controlling parameters $\lambda_i$ for $i=0,\ldots,k-1$.  Although the functional dependence of $\theta_\alpha$ on $\boldsymbol\lambda$ is generically complicated in the nested description implied by Eq.~\eqref{eq:lambdaParametrization:StructuralConditions},  
it is easy to show that 
\begin{equation}
\label{eq:StructuralParameters:lambdaVariation}
    \frac{\partial\theta_\alpha(\boldsymbol{\lambda})}{\partial\lambda_i} = - \left(\mathbf J_S^{-1}\mathbf M\right)_{\alpha,i}\quad\text{for}\quad \alpha=0,\ldots,d-1 \quad\text{and}\quad i=0,\ldots,k-1~,
\end{equation}
since the $d$ structural conditions Eq.~\eqref{eq:StructuralConditions} remain fixed under conditional sampling, meaning
\begin{equation}
\label{eq:ExperimentalConditions:lambdaVariation}
    \frac{\partial\hat\mu_\alpha}{\partial\lambda_i}\overset{!}{=}0~.
\end{equation}
From the block form Eq.~\eqref{eq:Jacobian_BlockForm} of the Jacobian $\mathbf J(\boldsymbol{\lambda})$ with 
\begin{equation}
\label{eq:nestedJacobian:lambdaParametrization}
    J_{\alpha\beta}(\boldsymbol{\lambda}) = \sum_{x\in\Omega} g_\alpha(x)\, p(x;\boldsymbol\lambda)\, g_\beta(x) \quad\text{for}\quad\alpha,\beta = 0,\ldots,d+k-1~,
\end{equation}
which is associated with the nested family, $\mathbf M$ and $\mathbf J_S$ inherit a nontrivial dependence on $\boldsymbol\lambda$.
Using the derived relation in Eq.~\eqref{eq:StructuralParameters:lambdaVariation}, the variation of experimental generalized moments is then dictated by Schur's complement  Eq.~\eqref{eq:SchurComplement}:
\begin{equation}
\label{eq:SchursComplement:mu_lambda_relation}
    \frac{\partial\mu_{d+j}(\boldsymbol{\lambda})}{\partial\lambda_i} = \left(\boldsymbol{\Sigma}^{-1}(\boldsymbol{\lambda})\right)_{ij} 
    \quad\text{for}\quad i,j=0,\ldots,k-1~,
\end{equation}
which also depends on $\boldsymbol{\lambda}$.

\paragraph{Unconditional multinomial sampling.}
Assume that the only structural condition we keep fixed is normalization, i.e.\ we multinomially sample $N$-sized datasets from $\mathbb N^{\vert\Omega\vert}_0\cap N\mathcal P$.
In this trivial scenario with $d=1$, the only structural generalized moment is $\hat\mu_0=1$, so that  
the Lagrange multiplier dual to the normalization condition on the simplex equals the negative logarithm of the partition sum corresponding to the $I$-projection of $\prob v$ onto $\mathcal L_{\boldsymbol{\lambda}}\subset\mathcal P$,
\begin{equation}
\label{eq:MultinomialSampling:Iprojection:PartitionSum}
    \theta_0(\boldsymbol{\lambda}) = -\log\sum_{x\in\Omega} v(x) \exp\left\{\sum_{i=0}^{k-1}\lambda_i g_{1+i}(x)\right\}~.
\end{equation}

Since the structural linear family covers the complete simplex over $\Omega$,  $\widehat{\mathcal L}=\mathcal P$, the ambient $I$-projection in Eq.~\eqref{eq:ambient_iProjection} trivially coincides with the reference distribution, $\hat{\prob q}=\prob v$. Furthermore, the ambient Jacobian Eq.~\eqref{eq:ambientJacobian} simplifies for $\prob p_{\boldsymbol{\lambda}}\in\mathcal P$ to $\widehat J_{00}=1$, while the determinant of the nested Jacobian in Eq.~\eqref{eq:nestedJacobian:lambdaParametrization} reduces to the determinant of the conventional covariance matrix 
\begin{equation}
    {C}_{ij}(\boldsymbol{\lambda}) = \sum_{x\in\Omega} \left[g_{1+i}(x) - \mu_{1+i}(\boldsymbol{\lambda})\right]p(x;\boldsymbol{\lambda})\left[g_{1+j}(x) - \mu_{1+j}(\boldsymbol{\lambda})\right]
    \quad\text{for}\quad i,j=0,\ldots,k-1~.
\end{equation}
This specializes the conditional probability mass density in Eq.~\eqref{eq:ConditionalProbabilityMassDensity:dualSpace} to the probability mass density over the $N$-dilated polytope corresponding to linear family $\mathcal L(\mathbf G;1,\mu_0(\boldsymbol{\lambda}),\ldots,\mu_{k-1}(\boldsymbol{\lambda}))$ under  unconditional multinomial sampling:
\begin{equation}
    \rho(\boldsymbol{\lambda}\vert\hat\mu_0=1) 
    =
    \sqrt{\det\left(\frac{N \mathbf C(\boldsymbol{\lambda})}{2\pi}\right)} \exp\left\{-N \min_{\prob p\in\mathcal L_{\boldsymbol{\lambda}}}\infdiv{\prob p}{\prob v} \right\}~.
\end{equation}

Regarding the subsequent Laplace approximation, a fully equivalent expression could be directly obtained by applying the saddle point approximation (Section~\ref{ssc:SaddlePointApproximation}) to the inverse Fourier transform of an alternative characteristic function to Eq.~\eqref{eq:CharacteristicFunction:Poissons}. Specifically, it can be derived
by slicing the multinomial characteristic function $\widetilde{\Phi}_{\mathcal P}$ in Eq.~\eqref{eq:CharacteristicFunction:Multinomial} to account for the experimental conditions Eq.~\eqref{eq:ExperimentalConditions}. Hence, the probability mass is formally recovered by the inverse Fourier transform of 
$\widetilde{\Phi}_{\mathcal P}(\boldsymbol s^T\mathbf G_E)$ with $\boldsymbol{s}\in\mathbb R^k$. 
However, special care is required when choosing the (principal) branch of the complex logarithm associated with normalization on the simplex (which automatically leads to Eq.~\eqref{eq:MultinomialSampling:Iprojection:PartitionSum}) upon analytic continuation over \(\boldsymbol{z}\in\mathbb C^k\).

\subsection{Cumulative probability mass under conditional sampling}
\label{ssc:LaplaceApproximation}

Using the conditional probability mass density derived in Eq.~\eqref{eq:ConditionalProbabilityMassDensity:dualSpace}, sums over dataset families can ultimately be replaced by integrals over some target region $\mathcal T\subset\mathbb R^k$ in the dual parameter space. 
This approximation makes particularly sense whenever the target region contains the global minimum of the corresponding effective action.

Applying the Pythagorean identity in Eq.~\eqref{eq:nestedDensity:Pythagorean}, we find that the effective action governing the conditional probability mass density in Eq.~\eqref{eq:ConditionalProbabilityMassDensity:dualSpace} is precisely a \textsc{kl} divergence and is therefore non-negative. 
The origin $\boldsymbol{\lambda}=\boldsymbol{0}$ of the coordinate system introduced in Eq.~\eqref{eq:nestedFamilies:lambdaParametrization} describes  a special member of the exponential family  in Eq.~\eqref{eq:lambdaParametrization:ExponentialFamily}. 
By the saddle-point construction, this distribution satisfies the structural conditions in Eq.~\eqref{eq:StructuralConditions} while its exponential form coincides with that of the $I$-projection of $\prob v$ onto the structural linear family $\widehat{\mathcal L}\equiv\mathcal L(\mathbf G_S;\hat\mu_0,\ldots,\hat\mu_{d-1})$.
Hence, it is identified (see Appendix~\ref{app:iProjections} and in particular discussion around Eq.~\eqref{eq:app:iProjection:Reparametrization}) with the unique $I$-projection of $\prob v$ onto the structural linear family,  
\begin{equation}
\label{eq:MinimizerEffectiveAction:ambien_iProjection}
    \prob p_{\boldsymbol\lambda=\boldsymbol0}  \equiv \hat{\prob q} = \argmin_{\prob p\in\widehat{\mathcal L}}\infdiv{\prob p}{\prob v}
    ~,
\end{equation}
already introduced in Eq.~\eqref{eq:ambient_iProjection}.

Accordingly, the effective action in the $\boldsymbol{\lambda}$-parametrization of Eq.~\eqref{eq:ConditionalProbabilityMassDensity:dualSpace} is minimized at the origin,
\begin{equation}
\label{eq:nested:EffectiveAction_Minimum}
    \left.\nabla_{\boldsymbol{\lambda}}\infdiv{\prob p_{\boldsymbol{\lambda}}}{\hat{\prob q}}\right\vert_{\boldsymbol\lambda=\boldsymbol\lambda^*}\overset{!}{=} 0 
    \quad\Rightarrow\quad \boldsymbol\Sigma^{-1}(\boldsymbol\lambda^*)\, \boldsymbol\lambda^* = \boldsymbol{0}
    \quad\Rightarrow\quad \boldsymbol{\lambda^*} = \boldsymbol{0}~,
\end{equation}
where the \textsc{kl} divergence attains its global minimum. 
The final implication follows from the positive-definiteness of  Schur's complement introduced in Eq.~\eqref{eq:SchurComplement}, guaranteed by the assumptions made in Section~\ref{app:IdealProblem} about the ideal classical problem. To establish the stationarity condition, it is convenient to apply the Pythagorean identity in the structural family in reverse,
\begin{equation}
    \infdiv{\prob p_{\boldsymbol{\lambda}}}{\hat{\prob q}} = \infdiv{\prob p_{\boldsymbol{\lambda}}}{\prob v} - \infdiv{\hat{\prob q}}{\prob v} 
    \overset{\eqref{eq:lambdaParametrization:ExponentialFamily}}{=} \sum_{\alpha=0}^{d-1}\theta_\alpha(\boldsymbol{\lambda}) \hat\mu_\alpha + \sum_{i=0}^{k-1}\lambda_i \mu_{d+i}(\boldsymbol{\lambda}) - \infdiv{\hat{\prob q}}{\prob v}~.
\end{equation}
Differentiating with respect to $\lambda_i$ and employing the identity
\begin{equation}
    \sum_{\alpha=0}^{d-1} \hat\mu_\alpha  \frac{\partial\theta_\alpha(\boldsymbol{\lambda})}{\partial\lambda_i} = - \sum_{\alpha,\beta=0}^{d-1} \left(\mathbf J_s(\boldsymbol{\lambda})\right)_{0\alpha} (\mathbf J_s^{-1}(\boldsymbol{\lambda}))_{\alpha\beta} M_{\beta,i}(\boldsymbol{\lambda}) = - M_{0,i}(\boldsymbol{\lambda}) = -\mu_{d+i}(\boldsymbol{\lambda}) ~,
\end{equation}
which follows from Eq.~\eqref{eq:StructuralParameters:lambdaVariation}, leads to the compact expression
\begin{equation}
    \frac{\partial\infdiv{\prob p_{\boldsymbol{\lambda}}}{\hat{\prob q}}}{\partial\lambda_i} = \sum_{j=0}^{k-1}\frac{\partial\mu_{d+j}(\boldsymbol{\lambda})}{\partial\lambda_i} \lambda_j \overset{\eqref{eq:SchursComplement:mu_lambda_relation}}{=} \sum_{j=0}^{k-1}\left(\boldsymbol{\Sigma}^{-1}(\boldsymbol{\lambda})\right)_{ij} \lambda_j~.
\end{equation}

As expected, minimizing the effective action $\infdiv{\prob p_{\boldsymbol\lambda}}{\hat{\prob q}}$ over the exponential family $\mathcal E(\mathbf G;\prob v)$, parametrized by its intersections with nested linear families $\mathcal L_{\boldsymbol\lambda}\subset\widehat{\mathcal L}$ uniquely recovers the ambient $I$-projection $\hat{\prob q}\in\widehat{\mathcal L}$. 
In the same $\boldsymbol{\lambda}$-coordinate system, the Hessian of the effective action
\begin{equation}
    \frac{\partial^2\infdiv{\prob p_{\boldsymbol\lambda}}{\hat{\prob q}}}{\partial\lambda_i\partial\lambda_j} = \left(\boldsymbol{\Sigma}^{-1}(\boldsymbol{\lambda})\right)_{ij} + \order{\Vert\boldsymbol{\lambda}\Vert}~,
\end{equation}
when evaluated at its minimum Eq.~\eqref{eq:nested:EffectiveAction_Minimum}, reduces to Schur's complement defined in Eq.~\eqref{eq:SchursComplement:mu_lambda_relation}.
Since the gradient of the effective action vanishes at its minimum, a second-order Taylor expansion yields the positive-definite quadratic approximation
\begin{equation}
\label{eq:nested:EffectiveAction:QuadraticExpansion}
    \infdiv{\prob p_{\boldsymbol\lambda}}{\hat{\prob q}}
    \approx  
    \tfrac12 \boldsymbol\lambda^T \boldsymbol{\Sigma}_*^{-1}\boldsymbol\lambda~,
\end{equation}
where the shorthand notation 
\begin{equation}
\label{eq:Covariancematrix:AtTheMinimum}
    \boldsymbol{\Sigma}_* \equiv \boldsymbol{\Sigma}(\boldsymbol\lambda=\boldsymbol0)
\end{equation}
emphasizes the evaluation of the precision matrix at the minimum of the effective action.

\paragraph{The Laplace approximation in dual parameter space.}
At sufficiently large $N$, Laplace's method implies that the integral approximating the cumulative probability mass is dominated by the minimum of the effective action, provided that this minimum lies within the target region.
Evaluating the determinant identity in Eq.~\eqref{eq:SchurComplement:detIdentity} at the minimizer $\boldsymbol{\lambda}=\mathbf0$ from Eq.~\eqref{eq:nested:EffectiveAction_Minimum} gives
\begin{equation}
    \det \mathbf J(\boldsymbol{\lambda}=\mathbf 0) =  \det \widehat{\mathbf J} \cdot \det \boldsymbol{\Sigma}^{-1}_*~.
\end{equation}
In particular,  the structural block $\mathbf J_S(\boldsymbol{\lambda}=\mathbf0)$ coincides with the ambient Jacobian $\widehat{\mathbf J}$ from Eq.~\eqref{eq:ambientJacobian}, by the identification of $I$-projections $\hat{\prob q}=\prob p_{\boldsymbol\lambda=\mathbf0}$ in Eq.~\eqref{eq:MinimizerEffectiveAction:ambien_iProjection}.
Consequently, the determinant ratio in Eq.~\eqref{eq:ProbabilityMass_Differential:DualSpace} simplifies to the correct normalization factor of the Gaussian density obtained from the quadratic approximation of the effective action in Eq.~\eqref{eq:nested:EffectiveAction:QuadraticExpansion}. 

Using the leading-order Gaussian density allows us to write the large-$N$ expansion of the cumulative probability mass: 
\begin{equation}
\label{eq:LapplaceApproximation}
    \int_{\mathcal T} d^k\!\boldsymbol{\lambda}\, \rho(\boldsymbol{\lambda}\vert\hat\mu_0,\ldots,\hat\mu_{d-1}) 
    =  
    \int_{\widetilde{\mathcal T}} d^k\!\tilde{\boldsymbol\lambda}\,  
    \sqrt{\det\left(\frac{\boldsymbol{\Sigma}_*^{-1}}{2\pi}\right)} \exp\left\{-\tfrac12 \tilde{\boldsymbol\lambda}^T \boldsymbol{\Sigma}_*^{-1}\tilde{\boldsymbol\lambda}\right\}
    \, + \,\order{N^{-1/2}} ~.
\end{equation}
Rescaling the dual variables according to  
\begin{equation}
    \tilde{\boldsymbol\lambda}=\sqrt{N}\,\boldsymbol\lambda
\end{equation}
exhibits the suppression of higher-order corrections, which here arise from higher-order terms in the Taylor expansion of the effective action. More generally, the discrete geometry of the target region and the details of the sampling process are also expected to contribute higher-order corrections.

\paragraph{The Laplace approximation in generalized moment space.}
To return from the $\boldsymbol{\lambda}$ coordinate system in Eq.~\eqref{eq:nestedFamilies:lambdaParametrization} to the generalized moment space, we Taylor expand Eq.~\eqref{eq:ExperimentalConditions:lambdaParametrization} for members of the exponential family Eq.~\eqref{eq:lambdaParametrization:ExponentialFamily},  
\begin{equation}
    \mu_{d+i} = \sum_{x\in\Omega}p(x;\boldsymbol{\lambda})\, g_{d+i}(x) ~,
\end{equation}
around the minimum  $\boldsymbol\lambda=\boldsymbol0$ of the effective action in Eq.~\eqref{eq:nested:EffectiveAction_Minimum}.  Combining Eq.~\eqref{eq:StructuralParameters:lambdaVariation}, which resolves the variation of structural parameters $\theta_\alpha(\boldsymbol{\lambda})$, with the definition of the covariance matrix Eq.~\eqref{eq:SchurComplement} evaluated at the minimum as in Eq.~\eqref{eq:Covariancematrix:AtTheMinimum}, we find that 
\begin{equation}
\label{eq:Linearization:Lambda_ExperimentalMu}
    \tilde\lambda_i +\order{N^{-1/2}} = \sqrt{N}\,\sum_{j=0}^{k-1} \left(\boldsymbol\Sigma_*\right)_{ij}\left(\mu_{d+j}-\hat\mu_{d+j}\right) \quad\text{for}\quad i = 0,\ldots,k-1~.
\end{equation}

This expansion expresses the intuitive fact that a $k$-dimensional parameter vector $\boldsymbol{\lambda}=N^{-1/2}\tilde{\boldsymbol\lambda}$ controlling the nested linear family $\mathcal L_{\boldsymbol{\lambda}}$ must be driven by the deviation of the dual experimental generalized moments from their estimates at the ambient $I$-projection Eq.~\eqref{eq:ambient_iProjection},
\begin{equation}
    \hat\mu_{d+i} = \sum_{x\in\Omega} \hat q(x)\, g_{d+i}(x)~.
\end{equation}
Recall that the ambient $\hat{\prob q}$ estimates the category  probabilities $\hat q(x)\equiv p(x;\boldsymbol{\lambda}=\boldsymbol{0})$ in the least-biased sense as the $I$-projection of $\prob v$ onto structural $\widehat{\mathcal L}$. Consequently, it represents the least-biased estimate --\,with respect to reference distribution $\prob v$\,-- of all generalized moments other than the structural moments $\hat\mu_\alpha$ for $\alpha=0,\ldots,d-1$ upon which we conditioned.   

Plugging the linearized expression Eq.~\eqref{eq:Linearization:Lambda_ExperimentalMu} into Eq.~\eqref{eq:LapplaceApproximation}, we obtain the leading-$N$ approximation to the density 
\begin{equation}
\label{eq:LapplaceApproximation:muSpace}
    \sqrt{\det\left(\frac{N\boldsymbol{\Sigma}_*}{2\pi}\right)}\exp\left\{-\frac N2\sum_{i,j=0}^{k-1}\left[\mu_{d+i}-\hat\mu_{d+i}\right]\left(\boldsymbol{\Sigma}_*\right)_{ij}\left[\mu_{d+j}-\hat\mu_{d+j}\right]\right\} 
\end{equation}
which determines the cumulative probability mass,
\begin{equation}
    \int_{\mathcal M}\rho(\mu_{d},\ldots,\mu_{d+k-1}\vert\hat\mu_0,\ldots,\hat\mu_{d-1}) 
    \left(\prod_{i=0}^{k-1}d \mu_{d+i}\right) ~,
\end{equation}
over a target region $\mathcal M\subseteq\mathbf G\mathbb R^{\vert\Omega\vert}$ in the continuum limit.
Hence, the Laplace approximation uncovers an effective Gaussian density over the experimental generalized moments Eq.~\eqref{eq:ExperimentalConditions}, analogous to the centered Gaussian density in the dual parameter space.
In contrast to the centered Gaussian approximation of Eq.~\eqref{eq:LapplaceApproximation} however, the Gaussian approximation for  the target variables $\mu_{d+i}$ is centered  around their least-biased estimate given by the  $I$-projection $\hat{\prob q}$ of $\prob v$ onto the structural linear family $\widehat{\mathcal L}$. 

The validity and subsequent computational tractability of the Laplace approximation (derived for $N\gg d+k$) to the cumulative probability mass of linear families neither relies on a specific form for $g_\alpha(X)$ nor presupposes a central limit theorem for the observables $X$ themselves.

\subsection{The \texorpdfstring{$p$}{p}-value and spherical symmetry}
\label{ssc:pValue_SphericalSymmetry}

Since $\Lambda(\mathbf G;\boldsymbol{m})\cap\Lambda(\mathbf G;\boldsymbol{m}')=\emptyset$ for distinct vectors $\boldsymbol{m}\neq\boldsymbol{m}'$ of raw means, the conditional sampling scheme introduced in Section~\ref{ssc:ConditionalSampling} naturally induces a structural decomposition of the form given by  Eq.~\eqref{eq:Intro:StructuralDecomposition}, namely 
\begin{equation}
    \hat{\Lambda}\equiv\Lambda(\mathbf G_S;\hat m_0,\ldots,\hat m_{d-1}) = \bigcup_{\boldsymbol{m}\in\mathbf G\mathbb N_0^{\vert\Omega\vert}} 
    \delta^d\!\left(\boldsymbol{m},\hat{\boldsymbol m}\right)
    \Lambda(\mathbf G;\hat m_0,\ldots,\hat m_{d-1},m_{d},\ldots,m_{d+k-1})~,
\end{equation}
where the product of Kronecker deltas enforcing the $d$  structural conditions is abbreviated as
\begin{equation}
    \delta^d\!\left(\boldsymbol{m},\hat{\boldsymbol m}\right) = \prod_{\alpha=0}^{d-1} \delta(m_\alpha,\hat m_\alpha)~.
\end{equation}

Suppose that we are interested in the coarse-grained probability mass of a particular linear family characterized by the observed raw effect $\boldsymbol m_\mathrm{obs} = \mathbf G\:\mathbf n_\mathrm{obs}\equiv N\boldsymbol\mu_\mathrm{obs}\in\mathbb R^{d+k}$. Obviously, this task becomes sensible, if the structural generalized moments implied by $\boldsymbol\mu_\mathrm{obs}$ agree with the structural generalized moments defining the conditional sampling scheme in Eq.~\eqref{eq:StructuralConditions}, 
otherwise the probability mass would vanish identically. In this context, where $\prob f\in\widehat{\mathcal L}\subseteq\mathcal P$ denotes the empirical distribution of the provided dataset $\mathbf n_\mathrm{obs}\in\mathbb N_0^{\vert\Omega\vert}$,  defined via the histogram (relative frequencies)
\begin{equation}
\label{eq:Histogram}
    f(x) = \frac{n_\mathrm{obs}(x)}{N}~,
\end{equation}
the linear system $\boldsymbol{\mu}_\mathrm{obs}=\mathbf G\:\prob f$ yields the empirical estimates of the generalized moments. 

According to prescription Eq.~\eqref{eq:Intro:pValue}, the critical region defined by the subset of admissible raw means (related to generalized moments by Eq.~\eqref{eq:Intro:GeneralizedMoment})
\begin{equation}
    \mathcal M = \left\{\left.\boldsymbol m\in\mathbf G\mathbb N_{0}^{\vert\Omega\vert}
    ~\right\vert~ 
    \Lambda_{\boldsymbol{m}}\subset\hat{\Lambda} \quad\text{and}\quad
    \operatorname{Pr}(\boldsymbol{m}) \leq \operatorname{Pr}(\boldsymbol m_\mathrm{obs})\right\}
\end{equation}
naturally leads to the sampling-based 
\begin{align}
    \text{$p$-value} = &\,\,
    \sum_{\boldsymbol m\in\mathcal M}\operatorname{Pr}(m_{d},\ldots,m_{d+k-1}\vert\hat m_0,\ldots,\hat m_{d-1}) 
    \\[1ex]
    = &\,\,
    \left[\operatorname{Pr}(\hat m_0,\ldots,\hat m_{d-1})\right]^{-1} 
    \sum_{N\boldsymbol\mu\in\mathcal M}\sum_{\mathbf n\in\mathbb N_0^{\vert\Omega\vert}}\operatorname{Pr}(\mathbf n)\, \delta(\mathbf n\in N\mathcal L_{\boldsymbol\mu})~,
    \nonumber
\end{align}
with the normalizing coarse-grained probability mass of the structural family (the discrete analogue of Eq.~\eqref{eq:StructuralFamily})  
being
\begin{equation}
    \operatorname{Pr}(\hat m_0,\ldots,\hat m_{d-1})= \sum_{\mathbf n\in\mathbb N_0^{\vert\Omega\vert}}\operatorname{Pr}(\mathbf n)\, \delta(\mathbf n\in N\widehat{\mathcal L})~.
\end{equation}

In the $\boldsymbol\lambda$-parametrization of Eq.~\eqref{eq:nestedFamilies:lambdaParametrization}, the linear family associated with the observed effect is parametrized by some (finite under the assumptions in Section~\ref{app:IdealProblem}) $\boldsymbol{\lambda}_\mathrm{obs}\in\mathbb R^k$,
\begin{equation}
    \label{eq:nested:criticalfamily}
    \mathcal L_\mathrm{obs}\equiv \mathcal L(\mathbf G;\boldsymbol{\mu}_\mathrm{obs})  =
     \mathcal L(\mathbf G;\hat\mu_0,\ldots,\hat\mu_{d-1},\mu_d(\boldsymbol\lambda_\mathrm{obs}),\ldots,\mu_{d+k-1}(\boldsymbol\lambda_\mathrm{obs}))
    ~.
\end{equation} 
$\mathcal L_\mathrm{obs}\neq\emptyset$, since at least the empirical distribution $\prob f$ in Eq.~\eqref{eq:Histogram} belongs to this family, i.e.\ $\prob f\in\mathcal L_\mathrm{obs}$.
In the continuum limit of preceding sections, we immediately deduce from Eq.~\eqref{eq:ConditionalProbabilityMassDensity:dualSpace} that the critical region dual to $\mathcal M$ is predominantly determined by the effective action, namely
\begin{equation}
\label{eq:pValue:CriticalRegion:LambdaParametrization}
    \mathcal T = \left\{\boldsymbol{\lambda}\in\mathbb R^k ~\vert~
    \min_{\prob p\in\mathcal L_{\boldsymbol{\lambda}}}\infdiv{\prob p}{\hat{\prob q}}
    \leq 
    \min_{\prob p\in\mathcal L_\mathrm{obs}}\infdiv{\prob p}{\hat{\prob q}}
    \right\}~.
\end{equation}
Subleading corrections to the specification of the integration region, arising for example from the determinant prefactor in the probability mass density, will contribute at $\order{N^{-1/2}}$ in the Laplace approximation.
Therefore, invoking Eq.~\eqref{eq:LapplaceApproximation}, we may estimate to leading order in $N$ the 
\begin{equation}
    \text{$p$-value} = 1 - 
    \int_{\mathcal T} d^k\!\boldsymbol{\lambda}\,\rho(\boldsymbol{\lambda}\vert\hat\mu_0,\ldots,\hat\mu_{d-1})~.
\end{equation}
In the information-theoretic picture, the complement of our $p$-value accounts for all $I$-projections of $\prob v$ (equivalently of $\hat{\prob q}$) onto nested linear families Eq.~\eqref{eq:nestedFamilies:lambdaParametrization} that are at least as close --\,in the \textsc{kl} sense\,-- to the ambient $I$-projection $\hat{\prob q}$ as the $I$-projection associated with the observed effect.

\paragraph{The emergence of spherical symmetry.}
The integration region Eq.~\eqref{eq:pValue:CriticalRegion:LambdaParametrization} implied by the definition of the $p$-value  evidently contains the minimum of the effective action (namely the governing \textsc{kl} divergence) at the origin $\boldsymbol{\lambda}=\mathbf{0}$.   
Within the Laplace approximation of Eq.~\eqref{eq:LapplaceApproximation}, in which the local geometry around the ambient $I$-projection at the origin is governed by a quadratic form, the constant-action hypersurfaces exhibit an emergent spherical symmetry. 

To this end, we rotate the coordinate system via the linear transformation
\begin{equation}
\label{eq:LaplaceApproximation:OrthogonalTrafo}
    \sqrt{N}\,\boldsymbol{\lambda}=\tilde{\boldsymbol\lambda} = \mathbf O \mathbf D^{-1/2}\boldsymbol{\xi}~,
\end{equation}
using a special orthogonal matrix $\mathbf O\in\mathbb R^{k\times k}$ to diagonalize the precision matrix, i.e.\
\begin{equation}
    \mathbf O^T \boldsymbol\Sigma^{-1}_* \mathbf O = \mathbf D \succ 0~.
\end{equation}
Since the columns of $\mathbf O$  are the orthonormalized eigenvectors of symmetric $\boldsymbol\Sigma^{-1}_*$, it follows that $\mathbf O^T\mathbf O=\mathbf O\mathbf O^T=\mathds1$.
The quadratic approximation of the effective action in Eq.~\eqref{eq:nested:EffectiveAction:QuadraticExpansion} then reduces to
\begin{equation}
    \infdiv{\prob p_{\boldsymbol\lambda}}{\hat{\prob q}} \approx \frac{\Vert\boldsymbol{\xi}\Vert^2}{2N}~.
\end{equation}
Hence, the \textsc{kl} divergence rescaled by $2N$ becomes, in the rotated coordinates $\boldsymbol{\xi}\in\mathbb R^k$, the squared Euclidean radius, which governs the leading contribution to the cumulative integral. 

Accordingly, the integration region $\mathcal T$ transforms into the $k$-dimensional ball 
\begin{equation}
    B^k(\boldsymbol{\xi}_\mathrm{obs}) = \left\{\boldsymbol{\xi}\in\mathbb R^k~\vert~\Vert\boldsymbol{\xi}\Vert^2\leq2N\min_{\prob p\in\mathcal L_\mathrm{obs}}\infdiv{\prob p_{\boldsymbol\lambda}}{\hat{\prob q}}\right\}~,
\end{equation}
with the radius of the boundary sphere fixed  in terms of the sampling probability in the vicinity of the observed linear family Eq.~\eqref{eq:nested:criticalfamily},
\begin{equation}
\label{eq:LaplaceApproximation:CriticalRadius}
     \min_{\prob p\in\mathcal L_\mathrm{obs}} \infdiv{\prob p_{\boldsymbol\lambda}}{\hat{\prob q}} \approx \frac12 \boldsymbol{\lambda}_\mathrm{obs}^T\:\boldsymbol{\Sigma}_*^{-1}\: \boldsymbol{\lambda}_\mathrm{obs} 
     \equiv\frac{1}{2N}\Vert\boldsymbol\xi_\mathrm{obs}\Vert^2~.
\end{equation}
In total, Eq.~\eqref{eq:LapplaceApproximation} amounts to the estimation of
\begin{equation}
\label{eq:LaplaceApproximation:pValue:SphericalSymmetry}
    \text{$p$-value} = 1 - \left(2\pi\right)^{-k/2}\int_{B^k(\Vert\boldsymbol{\xi}_\mathrm{obs}\Vert)} d \boldsymbol{\xi}\, e^{-\tfrac12\Vert\boldsymbol\xi\Vert^2}~.
\end{equation}
Note that the Jacobian of coordinate transformation Eq.~\eqref{eq:LaplaceApproximation:OrthogonalTrafo} produces a determinant factor which cancels exactly with the normalizing determinant prefactor Eq.~\eqref{eq:LapplaceApproximation}, 
\begin{equation}
    \det\left(\mathbf O \mathbf D^{-1/2}\right) = \det\mathbf O \cdot \left(\det \mathbf D\right)^{-1/2} =  \left(\det \mathbf \Sigma_*^{-1}\right)^{-1/2}~.
\end{equation}
This also justifies inverse and square root matrix operations, since the covariance matrix is positive definite, $\boldsymbol\Sigma(\boldsymbol{\lambda})\succ 0$.

In the emerging picture, we can distinguish  spherical shells around the origin, which corresponds to the ambient $I$-projection $\hat{\prob q}$, i.e.\ the least-biased estimate of category probabilities under conditional sampling from $\prob v$. Any linear dataset family in the vicinity of a given spherical shell of radius $\Vert\boldsymbol\xi\Vert$ carries the same coarse-grained probability mass in the leading-order description of the conditional sampling scheme. 
Besides its compatibility with the original premise of significance testing, our sampling-based $p$-value provides a homogeneous description of observable effects in the $k$-dimensional space of coarse-grained probability masses, treating all directions on equal footing. 
This contrasts with the heterogeneous structure in the space of generalized moments, which --\,prior to standardization\,-- are generally expressed in different physical units, and therefore lack a natural notion of joint magnitude. 

For instance, in a joint hypothesis test on the means of height and weight ($k=2$), there is no intrinsic, model-independent physical criterion determining how these observables should be combined into a single notion of ``typicality'' or ``extremeness''.
In the large-sample regime, however, the information-geometric structure resolves this ambiguity by organizing the space of coarse-grained probability masses into geometrically simple regions. In the plane with $k=2$, we obtain a central spherical disc of radius $\Vert\boldsymbol\xi_\mathrm{obs}\Vert$  corresponding to more probable fluctuations around the least-biased (given $\prob v$) estimate of the testable means, and its complementary tail region, which collects all less probable fluctuations.

\paragraph{The emergence of $\chi^2$ statistic.}
Once the spherical symmetry governing the Laplace integral in Eq.~\eqref{eq:LapplaceApproximation}  at large $N$ has been established, the remaining manipulations constitute a standard analytical exercise that connects to the widespread use of $\chi^2$ statistic. 

Introducing spherical coordinates in $\mathbb R^k$ with squared radius $t=\Vert\boldsymbol{\xi}\Vert^2$ transforms the conditional probability mass in the vicinity of nested family $\mathcal L(\mathbf G;\hat\mu_0,\ldots,\hat\mu_{d-1},\boldsymbol{\lambda}=\mathbf O (N\mathbf D)^{-1/2}\boldsymbol{\xi})$ from Eq.~\eqref{eq:LaplaceApproximation:pValue:SphericalSymmetry} into 
\begin{equation}
d^k\boldsymbol{\xi}\, \left(2\pi\right)^{-k/2} e^{-\tfrac12\Vert\boldsymbol\xi\Vert^2} =  d\Omega_{k-1}\,d t\, 2^{-\frac{k}{2}-1} \pi^{-\frac{k}{2}} t^{k/2-1} e^{-t/2}~.
\end{equation}
Since the quadratic form in the integrand equals by definition  the squared radial coordinate $t\geq0$, 
we can trivially perform the angular integration to produce the solid angle 
\begin{equation}
    \omega_{k-1} = \frac{2\pi^{k/2}}{\Gamma(k/2)}~.
\end{equation}
The remaining radial density can then be directly identified with the $\chi^2$ probability density function in the standard convention: 
\begin{equation}
    {d} t\,\frac{2^{-\frac{k}{2}}}{\Gamma(k/2)}   t^{k/2-1} e^{-t/2} \equiv \text{d}t\,\chi^2_k(t)~.
\end{equation}

The derived $\text{d}t\,\chi^2_k(t)$ gives the cumulative probability mass 
assigned to a spherical shell of squared radius $t$ (close to the origin) and thickness $dt$, comprising all nested families  with the same leading-$N$ probability mass density. 
Therefore, we arrive at the familiar expression 
\begin{equation}
    \text{$p$-value} = \int_{t_\mathrm{obs}}^\infty dt \, \chi_k^2(t) 
\end{equation}
with lower integration limit given by the observed squared radius $\Vert\boldsymbol{\xi}_\mathrm{obs}\Vert^2$ from Eq.~\eqref{eq:LaplaceApproximation:CriticalRadius}, 
\begin{equation}
\label{eq:LaplaceApproximation:CriticalValue}
    t_\mathrm{obs} = 
    2N\min_{\prob p\in\mathcal L_\mathrm{obs}} \infdiv{\prob p}{\hat{\prob q}}~.
\end{equation}
The closed-form expression shows that an intrinsic $p$-value can always be estimated under conditional  sampling, using a test statistic given by the rescaled \textsc{kl} divergence between the $I$-projections of the reference distribution $\prob v$ onto the nested family and the ambient family. In an effectively independent and unbiased sampling scheme, these $I$-projections are uniquely determined by the reference distribution $\prob v$ together with linear systems $\mathbf G\:\prob p=\boldsymbol\mu$ and $\mathbf G_S\:\prob p=\boldsymbol\mu_S$. Therefore, the resulting test statistic remains manifestly invariant under reparameterizations of the induced linear-algebraic problem in $\mathbb R^{d+k}$ and $\mathbb R^d$, respectively.  

It is remarkable that no underlying probability measure needs to be postulated, since the construction relies solely on fairly general universality arguments from thermodynamics and statistical mechanics, together with the large-$N$ geometry of coarse-grained probability masses.
Generically, higher-order corrections in $N$, as well as additional physical scales arising when the sampling process deviates from classical multinomial sampling, are expected to break the spherical symmetry. 

The spherical symmetry leading to Eq.~\eqref{eq:LaplaceApproximation:CriticalValue} enables efficient computation of $p$-values by relying on the emergent information-theoretic structure through the $I$-projector presented in Appendix~\ref{ssc:iProjector}. This suggests that future generalizations beyond the ideal classical regime could be formulated as perturbations around the classical solution presented here, thereby treating deviations from symmetry in a controlled manner.
Furthermore, the computational feasibility of the classical solution might be exploited to parametrize non-convex problems, when expectations satisfy higher-order polynomial conditions.


\section*{Acknowledgments}
\vspace{-0.2cm}
We are indebted to Prof.\ Ho Ryun Chung for many valuable discussions, for drawing our attention to critical references in the literature, as well as for proofreading the manuscript.

\printbibliography

\begin{appendix}
\section{On the limiting behavior of sampling processes}
\label{app:largeN}

Here, we merely outline universality classes of sampling schemes in the large-$N$ regime, characterized by the emergence of a common information geometry, independent of microscopic sampling details. 
This universality arises as an intrinsic large-$N$ property of the sampling scheme itself, and need not be derived from central limit theorems or large deviation principles in the classical probabilistic sense.

In particular, we consider known sampling mechanisms that deviate from the independence assumption (Section~\ref{sc:ClassicalSampling}) in a ``reasonable'' manner, governed by the interplay of physical scales in the underlying process. Via brute-force expansions, we show that despite such microscopic differences, these processes are governed by the same \textsc{kl} geometry and hence exhibit the emergent test-statistic of Section~\ref{sc:pValue} in the appropriate limit. 
For that, we make repeated use of Stirling's formula 
\begin{equation}
\label{eq:app:Stirling}
    \log n! = n \log n - n  + \tfrac12\log\left(2\pi n\right) + \order{n^{-1}}\quad\text{for}\quad n\in\mathbb N~,
\end{equation}
to approximate factorials of large counts $n$ by setting $n= N p(x)$, since the analysis remains (under the assumptions of Section~\ref{app:IdealProblem}) away from the boundaries of the simplex, i.e.\ $p(x)>0$ for all observable categories $x\in\Omega$.

\subsection{Sampling from a finite population without replacement}
From an urn of size $M$, we draw without replacement balls of different colors $x\in\Omega$ ($\Omega$ being isomorphic to the  space of observable categories introduced in Section~\ref{sc:ForwardProblem}), whose initial prevalence is given by $Mq(x)\in\mathbb N$.  
After $N$ draws, the probability of observing a dataset described by count vector $N\prob p\in\mathbb N^{\vert\Omega\vert}_0$ follows the standard multivariate hypergeometric distribution 
\begin{equation}
\label{eq:HypergeometricDistribution}
    \operatorname{hypergeom}(N\prob p;M\prob q) = 
    \left[\binom{M}{N} \right]^{-1}\prod_{x\in\Omega}\binom{Mq(x)}{Np(x)} 
    ~.
\end{equation}
As always, distributions lie on the simplex Eq.~\eqref{eq:simplex},  $\prob p,\prob q\in\mathcal P$, conveniently serving as auxiliary parameterizations for computation and asymptotic expansions.

In the double limit of large population and large sample size, i.e.\
\begin{equation}
\label{eq:LargeMNlimit}
    M\gg N \gg 1~,
\end{equation} 
we can identity two relevant expansion parameters
\begin{equation}
\label{eq:EffetiveSmallExpansionParameters}
    \frac{1}{N} \ll1 \quad\text{and}\quad \gamma\equiv\frac{N}{M} \ll1 ~,
\end{equation}
which conveniently organize the large-population limit of hypergeometric $\log$-probability according to
\begin{align}
\label{eq:Hypergeometric:LargeMNexpansion}
    \log\operatorname{hypergeom}(N\prob p;M\prob q) =&\,\, - N \left[\infdiv{\prob p}{\prob q} + \gamma \operatorname{D}_{\chi^2}\infdivx{\prob p}{\prob q} + \order{\gamma^2} \right]
    \\[0.9ex]
    \nonumber
    &\,\,
    - \frac{\vert\Omega\vert-1}{2} \log\left(2\pi N\right) - \tfrac12 \sum_{x\in\Omega}\log p(x) + 
    \frac{\gamma}{2}\left[\sum_{x\in\Omega}\frac{p(x)}{q(x)}-1\right]+ \order{\gamma^2} 
    \nonumber
     \\[0.9ex]
    &\,\,
    + \order{N^{-1}} ~.
    \nonumber
\end{align}

This neatly demonstrates how finite-population effects correct at $\order{\gamma}$ the leading \textsc{kl} divergence (cf.\ scaling law Eq.~\eqref{eq:SamplingProbability:ScalingLaw}) with the addition of Pearson's $\chi^2$ divergence 
\begin{equation}
\label{eq:PearsonChi2Divergence}
    \operatorname{D}_{\chi^2}\infdivx{\prob p}{\prob q} = \sum_{x\in\Omega} \frac{\left[p(x)-q(x)\right]^2}{2 q(x)} = \tfrac12 \left(\sum_{x\in\Omega}\frac{p(x)^2}{q(x)} - 1\right)~.
\end{equation}
For an infinite urn, we recover the multinomial scheme of independent sampling with replacement, 
\begin{equation}
    \lim_{M\rightarrow\infty}\operatorname{hypergeom} (N\prob p;M\prob q)=\operatorname{mult}(N\prob p;\prob q)~,
\end{equation}
where intrinsic probabilities are described by reference distribution $\prob v\equiv\prob q\in\mathcal P$.

\subsection{The general multicolor P\'olya urn}
\label{app::largeN:Polya}

In~\cite{grendar2010polya}, the large-$N$ behavior of multicolor P\'olya urns is systematically investigated. Starting from the unperturbed urn with prevalence vector $M\prob q$, as before, the authors consider a straightforward generalization: at each step of the sequential sampling process a ball is drawn from the urn, its color is recorded, and afterwards the ball is returned to the urn together with the addition or removal of $c\in\mathbb Z$ balls of the same color. 
A sample of $N$ balls then follows the P\'olya-Eggenberger distribution
\begin{equation}
    \operatorname{Pr}(N\prob p; M\prob q;c) = \frac{N!}{\prod\limits_{x\in\Omega}n(x)!} \frac{\prod\limits_{x\in\Omega} q(x)\left[q(x)+c\right]\cdots\left[q(x)+(n(x)-1)c\right]}{M\left(M+c\right)\cdots\left(M+(N-1)c\right)}~.
\end{equation}
In this formulation, $c=-1$ corresponds to the conventional hypergeometric scheme considered above while $c=1$  gives the Dirichlet-Multinomial compound distribution
\begin{equation}
\label{eq:app:DirMult}
    \operatorname{DirMult}(N\prob p;M\prob q) = 
    \frac{\Gamma(M)\,\Gamma(N+1)}{\Gamma(M+N)} \prod_{x\in\Omega} \frac{\Gamma(Mq(x)+Np(x))}{\Gamma(Mq(x))\,\Gamma(Np(x)+1)}~,
\end{equation}
in terms of factorials $\Gamma(n)=(n-1)!$ for $n\in\mathbb N_0$.

In analogy to the classical scaling law Eq.~\eqref{eq:SamplingProbability:ScalingLaw}, the P\'olya-Eggenberger sampling probability is governed by the  P\'olya information divergence, meaning that 
\begin{equation}
    -\lim_{N\rightarrow\infty}\log N^{-1}\operatorname{Pr}(N\prob p; M\prob q;c) = \operatorname{D}_{c\gamma}\infdivx{\prob p}{\prob q} ~,
\end{equation} 
using the expansion parameters identified in Eq.~\eqref{eq:EffetiveSmallExpansionParameters}.
This emerging divergence measure, which can be compactly expressed in terms of \textsc{kl} divergences,   
\begin{equation}
    \operatorname{D}_{c\gamma}\infdivx{\prob p}{\prob q}= \infdiv{\prob p}{\prob v} + \frac{1}{c\gamma}\infdiv{\prob q}{\prob v}
    \quad\text{with}\quad \prob v \equiv \frac{\prob q + c\gamma \prob p}{1+c\gamma}~,
\end{equation}
is also an $f$-divergence~\cite{csiszar1963}, i.e.
\begin{equation}
    \operatorname{D}_{c\gamma}\infdivx{\prob p}{\prob q}  = \sum_{x\in\Omega} q(x) f\!\left(\frac{p(x)}{q(x)}\right)
\end{equation}
for the convex function 
\begin{equation}
    f(t) = t\log t - \left(t+\frac{1}{c\gamma}\right)\log\left(t+\frac{1}{c\gamma}\right) + \frac{1+c\gamma}{c\gamma}\log\frac{1+c\gamma}{c\gamma}~.
\end{equation}

Assuming that $c=\order{1}$, so that the urn remains non-depleted in the double scaling regime described by Eq.~\eqref{eq:LargeMNlimit}, the small-$\gamma$ expansion of the governing P\'olya information divergence follows,
\begin{equation}
    \operatorname{D}_{c\gamma}\infdivx{\prob p}{\prob q} 
    =
    \infdiv{\prob p}{\prob q} - c\gamma \operatorname{D}_{\chi^2}\infdivx{\prob p}{\prob q} + \order{\gamma^2}
    ~,
\end{equation}
which agrees with Eq.~\eqref{eq:Hypergeometric:LargeMNexpansion} for $c=-1$, as expected. This highlights the universality of the \textsc{kl} divergence, beyond more conventional sampling schemes.

Therefore, reasonable modifications of the urn content --\,characterized by order-one $c$ in the present P\'olya-Eggenberger  scheme\,-- can be treated as controlled perturbations of the independent sampling regime.
Concerning the intrinsic $p$-value of Eq.~\eqref{eq:Intro:pValue}, the test statistic associated with the generalized sampling scheme can then be computed semi-analytically by incorporating perturbative corrections to the classical test statistic derived in Section~\ref{sc:pValue}.

\subsection{Wallenius' multivariate noncentral hypergeometric sampling}
\label{app::largeN:Hypergeometric}

From an urn of size $M$, we draw  balls of different colors $x\in\Omega$. 
Assume that the balls we draw --\,without replacement\,-- carry a weight $\omega(x)$ depending on their color. 
Over sequential draws, the probability of observing balls of one particular color is determined by the ratio of the total weight of that color remaining in the urn to the total weight of all remaining balls. Such biased sampling from an urn among competing colors is described by a multivariate extension to Wallenius' noncentral hypergeometric distribution~\cite{wallenius_biased_1963,fog2008sampling}:
\begin{equation}
\label{eq:app:WalleniusHyperGeometric}
    \operatorname{Pr} (N\prob p;M\prob q,\boldsymbol{\omega}) = 
    \left[\prod_{x\in\Omega}\binom{Mq(x)}{Np(x)} \right]
    \int_0^1  \text{d}\tau \prod_{x\in\Omega}\left(1-\exp\left\{\frac{\omega(x)}{\boldsymbol\omega \cdot \left(M\prob q-N\prob p\right)}\log\tau\right\}\right)^{Np(x)}
    ~,
\end{equation}
where we summarized weights by positive vector $\boldsymbol\omega\in\mathbb R^{\vert\Omega\vert}_{+}$.
If the latter were a constant vector, then the mixing integral compactly resolves,
so that we obtain the standard multivariate hypergeometric distribution of Eq.~\eqref{eq:HypergeometricDistribution}, 
which still describes a dependent sampling of balls (since each draw changes the composition of the finite urn). 

We focus on the latter mixing integral in Eq.~\eqref{eq:app:WalleniusHyperGeometric}, since the former product of binomial coefficients over $\Omega$ equals the unnormalized hypergeometric probability analyzed above. For a large population $M\gg1$, the mixing integral in Eq.~\eqref{eq:app:WalleniusHyperGeometric}  can be evaluated order-by-order to
\begin{align}
\label{eq:Wallenius:LargeM_expansion}
    \log\int_0^1d\tau\underbrace{ \prod_{x\in\Omega}\left(1-\tau^{\frac{\omega(x)}{\boldsymbol\omega \cdot \left(M\prob q-N\prob p\right)}}\right)^{Np(x)} }_{\equiv I}
    = 
    -N\log M + \log N! +&\,\, N\sum_{x\in\Omega} p(x)\log\frac{\omega(x)}{\boldsymbol\omega\cdot\prob q} 
    \\ 
    +&\,\, \frac{N(N-1)}{2M}\frac{\boldsymbol\omega\cdot\prob p}{\boldsymbol\omega\cdot\prob q} +  \order{M^{-2}}~,
    \nonumber
\end{align}
after 
expanding 
the product over $\Omega$ in the integrand as 
\begin{equation}
    I = M^{-N}\left(\prod_{x\in\Omega}\left[\frac{\omega(x)}{\prob q\cdot\boldsymbol\omega}\right]^{Np(x)}\right) \left(-\log\tau\right)^N\left[1 + \left(\frac{N^2}{M} + \frac{N}{2M}\log\tau\right) \frac{\boldsymbol\omega\cdot\prob p}{\boldsymbol\omega\cdot\prob q} + \order{M^{-2}}\right]
\end{equation}
and  applying the relation 
\begin{equation}
    \int_0^1\text{d}\tau \left(-\log\tau\right)^N \left(\log\tau\right)^k = (-1)^k \left(N+k\right)!
    \quad\text{for}\quad N,k\in\mathbb N_0 ~.
\end{equation} 
Formally, a combined large-population and large-sample limit exclusively entailing the scales in Eq.~\eqref{eq:LargeMNlimit} 
is understood 
for generic weights whose ratios do not present significant (w.r.t.\ two relevant scales) hierarchies, $\omega(x)/\omega(x')=\order{1}$ for $x,x'\in\Omega$.

Using the effective scales in Eq.~\eqref{eq:EffetiveSmallExpansionParameters}, the expansion of the logarithm of Wallenius' sampling probability thus organizes as in Eq.~\eqref{eq:Hypergeometric:LargeMNexpansion}. Specifically, we obtain
\begin{align}
\label{eq:app:Wallenius:LargeNexpansion}
    \log \operatorname{hypergeom} (N\prob p;M\prob q,\boldsymbol{\omega}) = 
    - N \left[\infdiv{\prob p}{\prob v} + \gamma \left(\operatorname{D}_{\chi^2}\infdivx{\prob p}{\prob v}-\tfrac12\frac{\boldsymbol\omega\cdot\prob p}{\boldsymbol\omega\cdot\prob q}\right) + \order{\gamma^2} \right]
    +\order{1}~,
\end{align}
where the reference distribution $\prob v\in\mathcal P$ now reflects the biased sampling scheme via
\begin{equation}
    v(x) \equiv  \frac{q(x) \omega(x)}{\prob q \cdot \boldsymbol\omega}~,
 \end{equation}
properly normalized by the sum
 \begin{equation}
     \prob q \cdot \boldsymbol\omega = \sum_{x\in\Omega} q(x)\omega(x)~,
 \end{equation}
which counts the initial total weight in the urn. 
In this scenario of competing sampling, the finite-population effects at order $\gamma$ and beyond receive $\boldsymbol{\omega}$-dependent corrections expressing the preferential sampling of colors. In particular at $\order{\gamma}$, the ratio $\tfrac12\frac{\boldsymbol\omega\cdot\prob p}{\boldsymbol\omega\cdot\prob q}$ enters as a subtraction from Pearson's $\chi^2$ divergence Eq.~\eqref{eq:PearsonChi2Divergence}, in order to obtain the net $\gamma$-correction to the $N$-leading \textsc{kl} divergence of the classical ($M\rightarrow\infty$) limit.

\subsection{Coupled Poisson Processes}
\label{app:largeN:KarlisMixture}
We are interested to record the emitted signals over a sufficiently long observation interval $T=\order{N}$. Each of the investigated processes that are labeled by $x\in\Omega$ has a mean count $Nv(x)$, so that we can collectively describe the Poisson point processes by mean vector $N\prob v\in\mathbb R_{+}^{\vert\Omega\vert}$. 
A mixture model can result from the contamination of $\vert\Omega\vert$ independent Poissonians of interest through a silent Poisson point process with latent mean $\bar m\in\mathbb R_{\geq0}$. 

In such a setting, the probability to observe the count vector $N\prob p\in\mathbb N_0^{\vert\Omega\vert }$ of interest, at fixed sample size $N$, is given by constraining~\cite{Karlis_Mixture_2005} according to
\begin{equation}
    \operatorname{Poisson}(N\prob p;N\prob v;\bar m) = \frac{1}{Z_{\bar m}(N\prob v)}\delta\left(\sum_{x\in\Omega} Np(x),N\right) 
    \prod_{x\in\Omega} \frac{(Nv(x))^{Np(x)}}{(Np(x))!} \,J_{\bar m}(N\prob p;N\prob v)
\end{equation}
with the mixing summation 
\begin{equation}
    J_{\bar m}(N\prob p;N\prob v) \equiv 
    \sum_{n=0}^{n_\text{max}} n! \left(\frac{\bar m}{\prod_{x\in\Omega}Nv(x)}\right)^n \prod_{x\in\Omega} \binom{Np(x)}{n}~,
\end{equation}
terminating at 
\begin{equation}
    n_\text{max} = \min_{x\in\Omega}Np(x) 
    \quad\text{since}\quad \binom{Np(x)}{n} = 0 ~\,\forall\, n>Np(x)~.
\end{equation}
The Kronecker delta emphasizes that the total count number observed on all detectors must equal $N$ by ensuring  $\prob p\in\mathcal P$. $Z_{\bar m}(N\prob v)$ denotes the partition function such that the total sampling probability normalizes to one.

If $\bar m=0$, the mixing sum becomes unity, $J_{\bar m=0}(N\prob p;N\prob v)=1$, so that the partition function exactly evaluates (under properly normalized reference distribution $\prob v\in\mathcal P$) to 
\begin{equation}
\label{eq:PoissonMixture:MultinomialLimit}
    Z_0\equiv Z_{\bar m=0}(N\prob v) = 
    \frac{N^N}{N!}~.
\end{equation} 
This is no surprise, as the sampled counts $N\prob p$ from all independent Poisson point processes at fixed sample size $N$ are described by the multinomial distribution, $\operatorname{Poisson}(N\prob p;N\prob v;\bar m=0)=\operatorname{mult}(N\prob p;\prob v)$. The equivalence of independent Poissonians under fixed total count number to  multinomial sampling has been used as a starting point for the systematic evaluation of coarse-grained probability masses, developed in Section~\ref{sc:FourierAnalysis}.

Switching on the contamination $\bar m>0$ generically makes even the numerical evaluation of $Z_{\bar m}$ challenging. Fortunately, we can still draw an important conclusion, as long as we focus on large sampled counts $Np(x)\gg1$. Consistency of the large-$N$ regime means that the partition function can be smoothly expanded around Eq.~\eqref{eq:PoissonMixture:MultinomialLimit}, i.e.\ we may write
\begin{equation}
    Z_{\bar m}(N\prob v) = 
    Z_0 + Z_1(N\prob v;\bar m)
    \quad\text{with}\quad \lim_{\bar m\rightarrow0} Z_1(N\prob v;\bar m) = 0 ~. 
\end{equation}
Always understanding normalization for $\prob p$, we thus have for the sampling probability of the Poisson mixture 
\begin{align}
\label{eq:app:Poisson:LargeNexpansion}
    \log \operatorname{Poisson}(N\prob p;N\prob v; \bar m) = 
    &\,\,- N \infdiv{\prob p}{\prob v} 
     - \frac{\vert\Omega\vert-1}{2} \log\left(2\pi N\right)  -\tfrac12\sum_{x\in\Omega}\log p(x) 
     \\[1ex]
     &\,\,
     + \log J_{\bar m}(N\prob p;N\prob v)
    - \log\left[1 + \frac{Z_1(N\prob v;\bar m)}{Z_0}\right]
    ~.
    \nonumber
\end{align} 
In the first row, we readily recognize the large-$N$ expansion of multinomial sampling probability Eq.~\eqref{eq:Multinomial:LargeNexpansion}.

Consequently, it remains to show that  the mixing summation is $\order{1}$ in the large-$N$ expansion. Since the product of binomial coefficients can be crudely bounded by   
\begin{align}
    \prod_{x\in\Omega}\binom{Np(x)}{n} = &\,\,(n!)^{-\vert\Omega\vert} \prod_{x\in\Omega} \underbrace{Np(x) \left[Np(x)-1\right]\cdots \left[Np(x)-n+1\right]}_{n\text{ terms }\leq Np(x)}
    \nonumber
    \\[1ex]
    <&\,\, (n!)^{-\vert\Omega\vert} \left[\prod_{x\in\Omega}N p(x)\right]^n~,
\end{align}
we indeed recognize for $\vert\Omega\vert\geq2$ that the mixing term 
\begin{equation}
   \log J_{\bar m} 
   < \log\sum_{n=0}^{n_\text{max}} 
    \frac{1}{(n!)^{\vert\Omega\vert-1}}  \left[\bar m\prod_{x\in\Omega}\frac{p(x)}{v(x)}\right]^n
   < \log\sum_{n=0}^{\infty} \frac{1}{n!}\left[\bar m\prod_{x\in\Omega} \frac{p(x)}{v(x)}\right]^n < \bar m\prod_{x\in\Omega} \frac{p(x)}{v(x)}
\end{equation}
cannot modify the $N$-leading expansion in the first row of Eq.~\eqref{eq:app:Poisson:LargeNexpansion}, as long as  $\bar m=\order{1}$. 
Relating to the scales of urn models in Eq.~\eqref{eq:LargeMNlimit}, we may parametrize $\bar m\equiv N^2/M$. Therefore, the \textsc{kl} divergence dictates again in the large-$M$-$N$ regime the sampling probabilities, where the observation time $T$ and hence sample size $N$ outperform  mean contamination $\bar m\ll N$.
\\

In the outlined examples, effectively independent sampling arises not because we have taken almost identical color weights  $\omega(x)\approx\text{const}$ in the urn or assumed almost vanishing contamination $\bar m\approx0$ in detection, but as a result of the large sample size $N$. This embodies the core idea of the classical regime, in which $N$ as a controlling parameter is sufficiently large in order to suppress other physical scales (such as $\omega(x)/\omega(x')$ or $\bar m$) that could otherwise noticeably alter the relevant structure of the sampling procedure. In that way, the scaling law in Eq.~\eqref{eq:SamplingProbability:ScalingLaw} emerges as universal in the classical regime of combinatorial sampling.

\section{Emerging Information geometry}
\label{app:InformationGeometry}

The scaling laws and approximations to the inherently discrete forward problem which follow from a large-$N$ expansion of sampling probabilities 
in Section~\ref{sc:ClassicalSampling} rely on the most canonical type of information geometry. 
Given its persistent emergence in the saddle-point formulation of Section~\ref{sc:FourierAnalysis} and Laplace approximation of Section~\ref{sc:pValue}, we review core properties in information geometry regarding information divergence and information projections onto linear families of distributions. In particular, we present the most straight-forward and stable way to numerically determine the probabilities of observable categories in the forward problem of linear conditions on generalized moments. 

Essentially, what follows is a recapitulation of Csiszár's treatment~\cite{csiszar_i-divergence_1975} of the forward $I$-projection, henceforth simply referred to as $I$-projection. Apart from the uniqueness and existence of the $I$-projection (\cite{csiszar_i-divergence_1975}, Theorem 2.1), we make frequent use of the {Pythagorean identity} for $I$-projections (\cite{cencov_statistical_1982}, Theorem 22.1; \cite{csiszar_i-divergence_1975}, Theorem 2.2), {nestedness} of $I$-projections 
(\cite{csiszar_i-divergence_1975}, Theorem 2.3) 
and the {exponential form} of the $I$-projection (\cite{csiszar_i-divergence_1975}, Theorem 3.1). These properties together uniquely characterize the $I$-projection and its underlying divergence, among all other divergence measures between probability distributions~\cite{amari_alpha-divergence_2009}.
To improve computational efficiency, we provide a generalized nestedness theorem.
Focusing on a finite space $\Omega$ of observable categories lessens the analytical burden and considerably simplifies the formalism.

\subsection{The ideal problem}
\label{app:IdealProblem}

For the purposes of the forward problem under the classical regime,  
we take any linear system from Eq.~\eqref{eq:LinearSystemDEF} to be consistent, 
\begin{equation}
    \label{eq:app:LinearSystem}
    \boldsymbol{\mu}\in\mathbb R^D\::\quad \exists\,\prob p\in\mathbb R^{\vert\Omega\vert} \quad\text{s.t.}\quad 
    \mathbf G\:\prob p=\boldsymbol{\mu} \quad\text{with}\quad D=\rank\mathbf G~,
\end{equation}
i.e.\ to admit at least one solution, as well as non-redundant, i.e.\ its coefficient matrix has full row rank.
In an ideal scenario, the family of distributions $\mathcal L\equiv\mathcal L(\mathbf G;\boldsymbol{\mu})$ defined by this linear system together with non-negativity  
does not force any observable category to have zero probability, meaning that
\begin{equation}
\label{eq:app:NoVanishingProbability}
     \forall\,x\in\Omega\,:\quad \exists\, \prob p\in\mathcal L \quad\text{s.t.}\quad p(x) > 0~.
\end{equation}
Compatibly, any reference distribution of interest $\prob v\in\mathcal P$  will similarly assign non-vanishing probability to all observable categories, $v(x)>0$ $\forall\,x\in\Omega$.

These assumptions have a profound consequence on the structure of linear families. For each observable category $x\in\Omega$, we fix a distribution $\prob p_x\in\mathcal L$ with $p_x(x)>0$. Next, we choose strictly positive coefficients $\lambda(x)>0$ for all $x\in\Omega$ with $\sum_{x\in\Omega}\lambda(x)=1$ and define the linear combination 
\begin{equation}
        \sum_{x\in\Omega} \lambda(x)\,\prob p_x \equiv \prob p_{\vec\lambda}  \in\mathcal L  ~.
\end{equation}
This also belongs to the linear family as a convex combination of distributions in $\mathcal L$, which is by itself convex.
For each $x\in\Omega$, the term $\lambda(x)p_x(x)$ contributes a strict positive amount to the sum over non-negative terms. Hence, $\prob p_{\vec\lambda}$ as a point within the convex hull of $\left\{\prob p_x\right\}_{x\in\Omega}$ assigns positive probability to all observable categories.  
Therefore, the interior of the linear family $\mathcal L$ is non-empty, i.e.\
\begin{equation}
\label{eq:app:NonEmptyInterior}
    \exists\,\prob p\in\mathcal L\subseteq\mathcal P \quad\text{with}\quad p(x) >0 ~~\forall\,x\in\Omega~.
\end{equation}

\subsection{Information projections}
\label{app:iProjections}

Let $\mathcal L$ be a linear family as in Eq.~\eqref{eq:LinearFamilyDEF}, alongside a reference distribution $\prob v\in\mathcal P$ not necessarily belonging to $\mathcal L$. 
We call the minimizer of the \textsc{kl} divergence from $\prob v$ over the linear family,
\begin{equation}
    \label{eq:iProjectionDEF}
        \prob q \equiv \prob q(\mathbf G;\boldsymbol{\mu};\prob v) = \argmin_{\prob p\in\mathcal L} \infdiv{\prob p}{\prob v}
\end{equation}
    the (forward or proper) $I$-projection of $\prob v$ onto $\mathcal L$. Indexed by observable categories $x\in\Omega$, it becomes a $\vert\Omega\vert$-dimensional column vector with components $q(x)\geq0$. 

\paragraph{Form and basic properties.}
As long as $\mathcal L$ does not force any vanishing probability in the sense of Eq.~\eqref{eq:app:NoVanishingProbability}, the $I$-projection of $\prob v$ onto $\mathcal L$ lies in its interior.
To see this, we take a distribution $\prob b\in\mathcal L$ from the boundary which exhibits $b(x_0)=0$ for some category $x_0\in\Omega$. 
Exploiting the convexity of the linear family $\mathcal L$, we introduce a set of distributions
    \begin{equation}
        \mathcal L\ni\prob q_\lambda = \lambda \prob p + (1-\lambda) \prob b \quad\text{for}\quad\lambda\in[0,1]
    \end{equation}
    that interpolates between the boundary distribution $\prob b$ and some distribution $\prob p\in\mathcal L$ in the interior where $p(x_0)>0$ for all $x\in\Omega$. 
    The derivative of \textsc{kl} divergence from the reference $\prob v$ is 
    \begin{equation}
        \frac{d\infdiv{\prob q_\lambda}{\prob v}}{d\lambda} = \sum_{x\in\Omega} \left[p(x) - b(x) \right] \log\frac{q_\lambda(x)}{v(x)} ~.
    \end{equation}
    At the boundary $\lambda=0$, this must explode, since 
    \begin{equation}
        \lim_{\lambda\rightarrow0^+} \frac{d\infdiv{\prob q_\lambda}{\prob v}}{d\lambda}  = \sum_{x\in\Omega} \left[p(x) - b(x) \right] \log\frac{b(x)}{v(x)}  = -\infty~,
    \end{equation}
    given that there is at least one manifestation with $b(x_0)=0$ but $p(x_0),v(x_0)>0$. Because the \textsc{kl} divergence vastly decreases moving away from the boundary, its minimum must lie in the interior of $\mathcal L\subseteq\mathcal P$, meaning $q(x)>0$ $\forall\,x\in\Omega$.

Most importantly, the $I$-projection of $\prob v$ onto $\mathcal L$ is unique in the ideal scenario.
Specifically, the Hessian of the \textsc{kl} divergence from $\prob v$ in the non-empty (Eq.~\eqref{eq:app:NonEmptyInterior}) interior of $\mathcal L$, where the $I$-projection must lie, stays positive definite,
\begin{equation}
    \frac{\partial^2\infdiv{\prob p}{\prob v}}{\partial p(x)^2}  
    =
    \frac{1}{p(x)} > 0 ~.
\end{equation}
Together with the convexity of the linear family in Eq.~\eqref{eq:LinearFamilyDEF}, the strict convexity of the \textsc{kl} divergence in its first argument guarantees uniqueness alongside existence of the I-projection.

By analogy to familiar Euclidean geometry, we proceed with the Pythagorean identity for \textsc{kl} divergences,
\begin{equation}
    \label{eq:PythagorasIdentity}
    \infdiv{\prob p}{\prob v} =  \infdiv{\prob p}{\prob q} + \infdiv{\prob q}{\prob v} \quad \forall\,\prob p\in\mathcal L~,
\end{equation}
which is true, if and only if $\prob q$ is the $I$-projection of $\prob v$ onto $\mathcal L$.
When Eq.~\eqref{eq:PythagorasIdentity} applies for all member distributions of the linear family for some $\prob q\in\mathcal L$, it indeed follows by non-negativity  of  the \textsc{kl} divergence that $$\infdiv{\prob p}{\prob q} = \infdiv{\prob p}{\prob v}-\infdiv{\prob q}{\prob v} \leq  \infdiv{\prob p}{\prob v}~~\forall\,\prob p\in\mathcal L~,$$ which readily reproduces the \textsc{kl} minimizer in Eq.~\eqref{eq:iProjectionDEF}. Hence, $\prob q$ must be the unique $I$-projection.
To constructively show the converse, that $\prob q=\argmin_{\prob p\in\mathcal L} \infdiv{\prob p}{\prob v}$ implies the Pythagorean identity, we consider a family of distributions $\prob p_t = t\,\prob p + \left(1-t\right)\prob q$ in $\mathcal L$ generated by $t\in[0,1]$ (convexity of the linear family means that $\prob p_t\in\mathcal L$ whenever $\prob p,\prob q\in\mathcal L$). 
Because the $I$-projection lies in the interior of $\mathcal L$, the directional derivative around $\prob q$ vanishes in every feasible direction, so that we must have  
\begin{equation}
    0\overset{!}{=}\left.\frac{\partial \infdiv{\prob p_t}{\prob v}}{\partial t}\right\vert_{t=0^+} = \sum_{x\in\Omega} \left[p(x) - q(x) \right] \log\frac{q(x)}{v(x)} \quad \Rightarrow\quad 
    \text{Eq.}~\eqref{eq:PythagorasIdentity}~.
\end{equation}

After excluding categories with $v(x)>0$ which would have been down-weighted to $q(x)=0$ by the conditions defining $\mathcal L$, as instructed in Eq.~\eqref{eq:app:NoVanishingProbability}, the $I$-projection of $\prob v$ onto $\mathcal L$ takes the exponential form 
\begin{equation}
\label{eq:iProjection:ExponentialForm}
    q(x) = v(x) \exp\left\{ \sum_{\alpha=0}^{D-1} \theta_\alpha \,g_\alpha(x) \right\}~.
\end{equation}
To quickly derive this form, we introduce the Lagrangian
\begin{equation*}
    L = \infdiv{\prob p}{\prob v}  - \boldsymbol \theta^T \left[ \mathbf G \: \prob p - \boldsymbol \mu \right]
\end{equation*}
using Lagrange multipliers $\boldsymbol\theta\in\mathbb R^D$ to implement the conditions that define the linear family~\eqref{eq:LinearFamilyDEF}. Exploiting the symmetries of the associated linear system, we can arbitrarily shift Lagrange multipliers. In particular, we can redefine   $\boldsymbol\theta\rightarrow\boldsymbol w + \boldsymbol\theta$ where $\boldsymbol w$ is given in Eq.~\eqref{eq:NormalizationSymmetry} to obtain 
\begin{equation}
\label{eq:iProjection:Lagrangian}
    L = \infdiv{\prob p}{\prob v}  - \boldsymbol \theta^T \left[ \mathbf G \: \prob p - \boldsymbol \mu \right] - \left[\sum_{x\in\Omega} p(x) - 1\right]~,
\end{equation}
which absorbs normalization condition on the simplex.
The variational extremum of $L$ in the probabilities $p(x)$ results in the compact exponential form
of the forward $I$-projection in Eq.~\eqref{eq:iProjection:ExponentialForm}.
The components $g_\alpha(x)$ of the coefficient matrix $\mathbf G$ define the linear family and the $D$ Lagrange multipliers  $\theta_\alpha\equiv\theta_\alpha(\mathbf G;\boldsymbol\mu;\prob v)$ enforce the conditions on expectations $\mathbf G\:\prob q=\boldsymbol\mu$ which are formally recovered by extremizing $L$ in $\boldsymbol\theta$. 

\paragraph{The dual formulation.}
In the dual picture, the parameter vector $\boldsymbol\theta\in\mathbb R^D$ induces an exponential family of distributions on the observable simplex
\begin{equation}
    \label{eq:ExponentialFamily}
    \mathcal E(\mathbf G;\prob v) = \left\{\prob p\in\mathcal P ~\text{ s.t. for } ~\boldsymbol\theta\in\mathbb R^D:~~ p(x) = v(x) \exp \left( \sum_{\alpha=0}^{D-1} \theta_\alpha\, g_\alpha(x) \right) ~~\forall\,x\in\Omega\right\}~.
\end{equation}
In physics literature, one more commonly encounters the equivalent form
\begin{equation}
\label{eq:ExponentialForm:WithPartitionSum}
    p(x)=Z^{-1}v(x)\exp\left\{\sum_{i=1}^{D-1}\vartheta_ig_i(x)\right\} \quad\text{with}\quad Z=\sum_{x\in\Omega}v(x)\exp\left\{\sum_{i=1}^{D-1}\vartheta_ig_i(x)\right\}
\end{equation} 
in terms of parameter vector $\boldsymbol\vartheta\in\mathbb R^{D-1}$. The partition function explicitly implements normalization on the simplex. In this dual formulation, the $I$-projection of $\prob v$ onto $\mathcal L$ represents the unique intersection of the exponential family $\mathcal E(\mathbf G;\prob v)$ with the linear family $\mathcal L$, see Corollary 3.1 in~\cite{csiszar_information_2004}.

Due to the reparametrization symmetry of the linear system\footnote{Under symmetry transformation $\mathcal L(\mathbf G,\boldsymbol\mu)=\mathcal L(\mathbf T\mathbf G,\mathbf T\boldsymbol\mu)$ for full rank matrix $\mathbf T\in\mathbb R^{D\times D}$, the dual parameters transform contravariantly as $\boldsymbol\theta\rightarrow\boldsymbol\theta'\equiv\left(\mathbf T^{-1}\right){\!}^T\boldsymbol\theta$.}, the parametrization in Eq.~\eqref{eq:iProjection:ExponentialForm} in terms of $\boldsymbol\theta$ is not unique. However, the $I$-projection~\eqref{eq:iProjectionDEF} of the linear system~\eqref{eq:LinearFamilyDEF} in probability space, i.e.\ the left-hand side of~\eqref{eq:iProjection:ExponentialForm}, is unique.  
Reversing the logic, we conclude that any distribution $\tilde{\prob q}\in\mathcal L$ of exponential form 
\begin{equation}
\label{eq:app:iProjection:Reparametrization}
    \tilde q(x) = v(x) \exp\left\{\sum_{\alpha=0}^{D-1} \tilde\theta_\alpha\, g_\alpha(x) \right\}~,
\end{equation}
which fulfills all conditions in $\mathcal L$ 
for appropriately chosen $\tilde{\boldsymbol\theta}\in\mathbb R^D$, coincides with the $I$-projection  of $\prob v$ onto $\mathcal L$.
Direct substitution of the exponential form  into the right-hand side of Eq.~\eqref{eq:PythagorasIdentity}, 
\begin{equation}
    \infdiv{\prob p}{\tilde{\prob q}} + \infdiv{\tilde{\prob q}}{\prob v} = \infdiv{\prob p}{\prob v} + \sum_{\alpha=0}^{D-1} \tilde\theta_\alpha \underbrace{\vec g_\alpha\cdot\left(\tilde{\prob q} -  \prob p\right)}_{=0\text{ since }\prob p,\,\tilde{\prob q}\in\mathcal L} ~,
\end{equation}
namely reveals that $\tilde{\prob q}$ satisfies the Pythagorean identity for any $\prob p\in\mathcal L$ and must therefore be the unique $I$-projection of $\prob v$ onto $\mathcal L$.

To relate to the likelihood-based formulation of sampling that dominates the statistics literature, we 
introduce the empirical distribution $\prob f\in\mathcal P$, whose probabilities are estimated by the relative frequencies $f(x)$ of observed categories in the given data matrix. 
Substituting a ``model'' distribution $\prob p_{\boldsymbol\theta}$ from the exponential family Eq.~\eqref{eq:ExponentialFamily} in expression Eq.~\eqref{eq:intro:logLikelihood} we find that 
\begin{equation}
\label{eq:LogLikelihoodOnExponentialFamily}
   \ell_{N\prob f}(\prob p_{\boldsymbol\theta}) - \ell_{N\prob f}(\prob v) = \bar{\boldsymbol\mu}\cdot\boldsymbol\theta  
   \overset{\eqref{eq:ExponentialForm:WithPartitionSum}}{=} \sum_{i=1}^{D-1}\bar{\mu}_i\theta_i - \log\sum_{x\in\Omega} v(x) \exp\left\{\sum_{i=1}^{D-1}\theta_ig_i(x)\right\}~,
\end{equation}
in terms of the empirical estimates $\bar\mu_i=\vec g_i\cdot\prob f$ for non-trivial generalized moments in the received data. 
Hence, the log-likelihood becomes a function of the parameters $\theta_i$ for $i=1,\ldots, D-1$. 
From here, we recognize that the $I$-projection of $\prob v\in\mathcal P$ onto the linear family $\widebar{\mathcal L}\equiv \mathcal L({\mathbf{G}}; \bar{\boldsymbol{\mu}})$ corresponds to the {maximum likelihood estimate} in the corresponding exponential family $\mathcal E(\mathbf G;\prob v)$.  

To see this (cf.~\cite{csiszar_information_2004} Theorem~3.3), we extremize the Lagrangian 
\begin{equation}
    L' = \ell_{N\prob f} (\prob p_{\boldsymbol\theta}) - \lambda \left( \sum_{x\in\Omega} p(x;\boldsymbol{\theta}) - 1 \right) \quad\text{for}\quad \prob p_{\boldsymbol{\theta}}\in\mathcal E~,
\end{equation}
so that the log-likelihood function Eq.~\eqref{eq:LogLikelihoodOnExponentialFamily} is maximized over the exponential family parametrized by $\boldsymbol{\theta}$ in Eq.~\eqref{eq:ExponentialFamily}.
This produces the \textsc{mle} system (in the canonical parametrization where $g_0(x)=1$ $\forall\, x\in\Omega$)
\begin{align}
    \frac{\partial L'}{\partial \lambda} = 0& \quad\Rightarrow\quad \theta_0 = -\log \sum_{x\in\Omega} v(x) \exp\left\{\sum_{i=1}^{D-1}\theta_ig_i(x)\right\}
    \\[1ex]
    \alpha=0,\ldots,D-1:~~ \frac{\partial L'}{\partial\theta_\alpha} = 0&  \quad\Rightarrow\quad
     \bar\mu_\alpha = \sum_{x\in\Omega}g_\alpha(x)f(x)  \overset{!}{=}  \lambda  \sum_{x\in\Omega}g_\alpha(x) v(x) \exp\left\{\sum_{\beta=0}^{D-1}\theta_\beta g_\beta(x)\right\}~.
     \nonumber
\end{align}
Since $\bar\mu_0=\vec g_0\cdot \prob f=1$, it follows from the extremum condition for $\alpha=0$  that   
$\lambda=1$ --\,upon substitution of the partition sum $\theta_0$ from the first line. 
Consequently, the latter $D-1$ equations defining the extremum of $L'$ readily reproduce the extremum equations for the Lagrangian $L$ in Eq.~\eqref{eq:iProjection:Lagrangian} corresponding to the non-trivial (beyond normalization) conditions that define the $I$-projection of $\prob v$ onto $\widebar{\mathcal L}$.
This interpretation of $I$-projections in terms of the \textsc{mle} does not by itself justify why one should restrict attention to the given exponential family Eq.~\eqref{eq:ExponentialFamily}, in the first place.

\subsection{Information projector}
\label{ssc:iProjector}

The preceding paragraphs reviewed the $I$-projection Eq.~\eqref{eq:iProjectionDEF} of $\prob v$ onto $\mathcal L$ with its characteristic exponential form Eq.~\eqref{eq:iProjection:ExponentialForm}. In this section, we provide a straightforward algorithm to efficiently compute the $I$-projection $\prob q$ given some (not necessarily normalized) reference $\prob{v}$ and generalized moments $\boldsymbol\mu$ defining the linear family $\mathcal L$.

Generically, the Lagrange multipliers $\theta_\alpha\in\mathbb R$ viewed as functions of the defining generalized moments $\mu_\alpha$ do not posses a  closed-form solution. 
This lack of an explicit functional form for parameters $\boldsymbol\theta\equiv\boldsymbol\theta(\boldsymbol\mu)$ given sample means $\boldsymbol\mu$ makes the inverse problem notoriously cumbersome to operationally deal with. 
Despite the alleged burden of the inverse problem, the probabilities $q(x)$ assigned by the $I$-projection to observable categories can be numerically determined by cleverly applying the Newton-Raphson root-finding method. Choosing to operate in probability space, we take advantage of the invariance of $q(x)$ under reparametrization of the Lagrange multipliers, cf.\ discussion around Eq.~\eqref{eq:app:iProjection:Reparametrization}. 

As input, the $I$-projector takes the reference distribution $\prob v\in\mathcal P$ or a vector $\mathbf v\in\mathbb R_+^{\vert\Omega\vert}$ of reference weights, alongside the coefficient matrix $\mathbf G\in\mathbb R^{D\times\vert\Omega\vert}$ and a vector of sample means $\boldsymbol{\mu}\in\mathbb R^D$.
Starting from $\prob p^{(0)}=\prob v$, the update rule for $\prob p^{(n+1)}$ after $n$ iterations reads 
\begin{equation}
\label{eq:NewtonStep}
p^{(n+1)}(x) = p^{(n)}(x)\: \exp\left\{- \sum_{\alpha,\beta=0}^{D-1} g_\alpha(x)\: (\mathbf J^{(n)})^{-1}_{\alpha\beta} \left( \sum_{x'\in\Omega} g_\beta(x') \:p^{(n)}(x') - \mu_\beta \right)\right\}~,
\end{equation} 
in terms of previous estimates $p^{(n)}(x)$ and the symmetric positive definite and thus invertible Jacobian matrix $\mathbf J^{(n)}\in\mathbb R^{D\times D}$ with components
\begin{equation}
\label{eq:iProjector:Jacobian}
J^{(n)}_{\alpha\beta} = \sum_{x\in\Omega} g_\alpha(x)\:p^{(n)}(x)\: g_\beta(x)~.
\end{equation}
As long as $\mathcal L\neq\emptyset$ and Eq.~\eqref{eq:app:NoVanishingProbability} holds, 
the numerical estimate locally converges under standard regularity conditions to the $I$-projection of $\prob v$ onto $\mathcal L$, 
\begin{equation}
    \prob p^{(n)}\rightarrow\prob q=\argmin_{\prob p\in\mathcal L}\infdiv{\prob p}{\prob v}~.
\end{equation}

Using that the $I$-projection is a member of the exponential family $\mathcal E(\mathbf G;\prob v)$, we introduce functions 
\begin{equation}
    \label{eq:iProjector:ConstraintFunction}
    F_\alpha^{(n)} =
    \sum_{x\in\Omega} g_\alpha(x) v(x) \exp\left\{\sum_{\beta=0}^{D-1} \theta_\beta^{(n)} g_\beta(x) \right\} - \mu_\alpha
    \quad\text{for}\quad\alpha=0,\ldots,D-1~,
\end{equation}
whose roots characterize the intersection of $\mathcal E(\mathbf G;\prob v)$ with the linear family $\mathcal L(\mathbf G,\boldsymbol\mu)$, see comment below Eq.~\eqref{eq:ExponentialFamily}. 
The $D$-dimensional vector $\mathbf{F}^{(n)} = \mathbf G\:\prob p^{(n)}_{\boldsymbol\theta} - \boldsymbol\mu$ imposes the conditions on sample means defining the linear family as constraints on the dual parameter vector $\boldsymbol\theta\in\mathbb R^{D}$.
To apply Newton-Raphson algorithm to find the roots of $\mathbf F\equiv \mathbf F(\boldsymbol\theta)$, we take the derivative with respect to the dual parameters, 
\begin{equation}
\label{eq:app:iProjector:GradientF}
        \frac{\partial F_\alpha^{(n)}}{\partial\theta_\beta} \equiv J^{(n)}_{\alpha\beta}~,
\end{equation}
which immediately gives the stated Jacobian Eq.~\eqref{eq:iProjector:Jacobian}. This is a manifestly symmetric square matrix. 

The stated update rule after $n$ iterations follows from the standard Newton-Raphson update rule for the dual parameter vector 
\begin{equation}
    \boldsymbol\theta^{(n+1)} = \boldsymbol\theta^{(n)} -  (\mathbf J^{(n)})^{-1} \mathbf F^{(n)}~,
\end{equation}
using the vector-valued function Eq.~\eqref{eq:iProjector:ConstraintFunction} and the inverse of its Jacobian Eq.~\eqref{eq:app:iProjector:GradientF}. After acting on both sides by $\mathbf G^T$, exponentiating componentwise, and multiplying by the reference weights $v(x)$, we obtain
\begin{equation*}
        v(x)\exp \left\{ \sum_{\alpha=0}^{D-1} \theta^{(n+1)}_\alpha g_\alpha(x)\right\} = v(x)\exp\left\{\sum_{\alpha=0}^{D-1} \theta^{(n)}_\alpha g_\alpha(x)\right\} \exp\left\{-  \sum_{\alpha,\beta=0}^{D-1} g_\alpha(x)(\mathbf J^{(n)})^{-1}_{\alpha\beta} F^{(n)}_\beta\right\}~,
\end{equation*}
Identifying the probability estimates on both sides, we recover the update rule Eq.~\eqref{eq:NewtonStep} in probability space. 

That the Jacobian is positive definite --\,needed for a well-behaved numerical optimization\,-- can be easily seen by inspecting its quadratic form,  
\begin{equation}
    \label{eq:iProjector:Jacobian:QuadraticForm}
        \boldsymbol u^T \mathbf J^{(n)} \boldsymbol u =\sum_{\alpha,\beta=0}^{D-1} u_\alpha J_{\alpha\beta}^{(n)} u_\beta = \sum_{x\in\Omega} p^{(n)}(x) \left[ \sum_{\alpha=0}^{D-1} u_\alpha\, g_\alpha(x) \right]^2 ~,
    \end{equation}
which can be written as a weighted sum of squares.
Recall that $p^{(n)}(x)>0$ for all $x\in\Omega$ as long as the iterates remain in the interior of Eq.~\eqref{eq:app:NonEmptyInterior}. Furthermore, $\mathbf G$ has by the assumption in Eq.~\eqref{eq:app:LinearSystem} full row rank, meaning that 
    \begin{equation}
        \forall \,\boldsymbol u \in\mathbb R^D\setminus\{\boldsymbol0\} : \quad\exists\, x\in\Omega \quad\text{s.t.} \quad \sum_{\alpha=0}^{D-1} u_\alpha\, g_\alpha(x) \neq 0 ~.
    \end{equation}
Consequently, there exists at least one category $x$ for which the corresponding summand in Eq.~\eqref{eq:iProjector:Jacobian:QuadraticForm} is strictly positive.  Hence, the quadratic form of the Jacobian is indeed strictly positive, 
\begin{equation}
    \forall \,\boldsymbol u \in\mathbb R^D\setminus\{\boldsymbol0\} : \quad  \boldsymbol u^T \mathbf J^{(n)} \boldsymbol u >0 ~.
\end{equation}
Notice that the Jacobian in Eq.~\eqref{eq:iProjector:Jacobian} converges to a positive definite $\mathbf J$ with the non-trivial generalized moments given by $J_{0i}=J_{i0}=\mu_i$ for $i=1,\ldots,D-1$.

The $I$-projector --\,relying solely on the exponential form of its target Eq.~\eqref{eq:iProjection:ExponentialForm}\,-- streamlines the conditional minimization of \textsc{kl} divergence in its first argument across a wide range of settings. It yields a remarkably stable and fast algorithm under different sets of generalized moment conditions. The rapid convergence  $\prob p^{(n)}\rightarrow\prob q$ contrasts with optimization in the dual parameter space, where the efficiency of  $\boldsymbol\theta^{(n)}\rightarrow{\boldsymbol\theta}$ depends on a ``good'' parametrization, which in practice can be challenging to uncover. 

In the exponential family parametrized as in Eq.~\eqref{eq:ExponentialFamily}, normalization on the simplex is always implied by $\mathbf G$, but is not enforced explicitly at each step of Newton-Raphson optimization. On the one hand, this means that, at a given step of~\eqref{eq:NewtonStep}, the running estimates $\prob p^{(n)}$ need not be normalized within the desired numerical accuracy. On the other hand, it is this form that gives remarkable flexibility to even accommodate categories that are forced to receive zero probability by the conditions defining the linear system Eq.~\eqref{eq:LinearFamilyDEF}. If present, such potential outliers would violate the assumed exponential form of the $I$-projection, which would consequently be forced onto the boundary of the observable simplex.  By numerically suppressing their probabilities below machine precision, our $I$-projector naturally identifies such affected categories.

\subsection{Further properties}

In Section~\ref{sc:pValue}, we considered conditional sampling schemes for $D>d$. In the linear-algebraic language of Section~\ref{ssc:LinearAlgebra},
we summarize the set of conditions on structural generalized moments by $\mathbf G_S\in\mathbb R^{d\times\vert\Omega\vert}$ and the combined set of conditions on both structural and experimental generalized moments by $\mathbf G\in\mathbb R^{D\times\vert\Omega\vert}$. 
Nestedness means that coefficient matrices are related through a transformation matrix $\mathbf T\in\mathbb R^{d\times D}$ with $\rank\mathbf T=d$ according to 
\begin{equation}
\label{eq:NestedTrafo}
\mathbf G_S = \mathbf T\:\mathbf G
\end{equation}
while the combined generalized moments $\boldsymbol{\mu}\in\mathbb R^D$ imply the structural $\boldsymbol{\mu}_S\in\mathbb R^d$, meaning
\begin{equation}
    \boldsymbol\mu_S = \mathbf T\:\boldsymbol\mu~.
\end{equation}
According to Eq.~\eqref{eq:LinearFamilyDEF}, this double implication induces an inner linear family $\mathcal L_E\equiv\mathcal L(\mathbf G; \boldsymbol\mu)$ which is nested within the ambient linear family $\mathcal L_S\equiv\mathcal L(\mathbf G_S;\boldsymbol\mu_S)$.
Every linear family $\mathcal L$ we consider in Eq.~\eqref{eq:LinearFamilyDEF} is nested within the observable simplex Eq.~\eqref{eq:simplex}, $\mathcal L\subseteq\mathcal P$.  

\paragraph{Information projections in nested descriptions.}
An important property of $I$-projections is their consistency under nested linear families (Theorem~2.3 in~\cite{csiszar_i-divergence_1975}, summarized by~Lemma 4.2 in~\cite{csiszar_information_2004}), paralleling the familiar projections in Euclidean geometry. Specifically, we consider two  nested linear families  $\mathcal L_E\subseteq\mathcal L_S$. The $I$-projection of reference distribution $\prob v\in\mathcal P$ onto the inner linear family $\mathcal L_E$ can also be achieved via an intermediate $I$-projection onto the ambient linear family $\mathcal L_S$: 
\begin{equation}
\label{eq:app:nested_iProjections}
    \prob q_E =\argmin_{\prob p\in\mathcal L_E} \infdiv{\prob p}{\prob v} = \argmin_{\prob p\in\mathcal L_E} \infdiv{\prob p}{\prob q_S} \quad\text{with}\quad
    \prob q_S = \argmin_{\prob p\in\mathcal L_S}\infdiv{\prob p}{\prob v}
\end{equation}

Let us denote by $\prob q_E'$ the $I$-projection of $\prob q_S\in\mathcal L_S$ onto the inner linear family $\mathcal L_E$, meaning that 
\begin{equation}
\label{eq:nested:q2prime}
    \infdiv{\prob q_E'}{ \prob q_S} \leq \infdiv{\prob p}{\prob q_S} \quad\forall\,\prob p\in\mathcal L_E~.
\end{equation}
From the Pythagorean identity for the $I$-projection of reference $\prob v$ onto the ambient linear family, we know that
\begin{equation}
    \infdiv{\prob p}{\prob q_S} = \infdiv{\prob p}{\prob v} - \infdiv{\prob q_S}{\prob v} \quad\text{for any $\prob p\in\mathcal L_S$}~.
\end{equation}
Applying this identity to both sides of Eq.~\eqref{eq:nested:q2prime}, since $\prob p,\prob q_E'\in\mathcal L_E\subseteq\mathcal L_S$, reduces the defining inequality for $\prob q_E'$ to 
\begin{equation}
    \infdiv{\prob q_E'}{ \prob v} \leq \infdiv{\prob p}{\prob v} \quad\forall\,\prob p\in\mathcal L_E~.
\end{equation}
This is precisely the defining property of the $I$-projection of $\prob v$ onto $\mathcal L_E$. Since the $I$-projection onto a linear family is unique, we conclude that $\prob q_E'=\prob q_E$.

Especially when $\infdiv{\prob q_E}{\prob v}\gg1$, an important practical consequence of the Pythagorean identity ensures that the ambient $I$-projection $\prob q_S$ is always closer to the inner $I$-projection $\prob q_E$ than the reference distribution $\prob v$ itself, namely $\infdiv{\prob q_E}{\prob q_S} \leq \infdiv{\prob q_E}{\prob v}$. 
To this end, we state the following proposition:
If the $I$-projector in Section~\ref{ssc:iProjector} converges slowly because the $I$-projection lies far from the reference $\prob v$, or if an ambient $I$-projection has already been computed, then the original optimization can be decomposed into a sequence of nested optimization problems; thereby improving both convergence and computational efficiency.
  
After organizing the given set of conditions which define a target linear family $\mathcal L$ into $k$ nested subsets inducing nested linear families on the simplex, $\mathcal L\equiv\mathcal L_k\subset\dots\subset\mathcal L_2 \subset \mathcal L_1$,  we successively impose the corresponding conditions along the nested structure. 
Starting from $\prob q_0\equiv\prob v$, we define a sequence of intermediate $I$-projections via 
\begin{equation}
\label{eq:nested:ChainofIprojections}
        \prob q_n = \argmin_{\prob p\in\mathcal L_n}\infdiv{\prob p}{\prob q_{n-1}}\quad\text{for}\quad n = 1,\ldots k~.
\end{equation} 
By Eq.~\eqref{eq:app:nested_iProjections}, the final estimate $\prob q_{k}=\prob q$ coincides with the $I$-projection of the original reference $\prob v$ onto the innermost linear family $\mathcal L$. Although Eq.~\eqref{eq:nested:ChainofIprojections} is expected to be numerically more stable than directly projecting $\prob v$ onto $\mathcal L$, failure of the $I$-projector at any intermediate step breaks the chain, in which case one should consider a different decomposition of target conditions or an alternative method of numerical optimization.  

\paragraph{Translation invariance and generalized nestedness.}
Another key structural property of the \textsc{kl} divergence is the transitivity of $I$-projections, as formalized in Lemma~4.2 of~\cite{csiszar_information_2004}. In particular point~(ii) of that Lemma reveals a form of ``translation'' invariance in Csiszár’s construction, whereby changes in the vector of generalized moments $\boldsymbol\mu \in \mathbb R^D$ correspond to shifts in the dual exponential parametrization, while the coefficient matrix $\mathbf G$ remains fixed.
This insight motivates a generalized notion of nestedness: once the $I$-projection of $\prob v$ associated with a given vector $\boldsymbol\mu$ has been computed, it can serve as a stable intermediate reference for efficiently updating the $I$-projection under a slightly perturbed vector $\tilde{\boldsymbol\mu}$.

For $D\leq \tilde D$, consider coefficient matrices $\mathbf G\in\mathbb R^{D\times\vert\Omega\vert}$ and $\widetilde{\mathbf G}\in\mathbb R^{\tilde D\times\vert\Omega\vert}$ related as in Eq.~\eqref{eq:NestedTrafo} via $\mathbf G=\mathbf T\,\widetilde{\mathbf G}$. For concreteness, we work in a representation where the first $D$ conditioning functions coincide in both descriptions, i.e.\ the first rows of $\widetilde{\mathbf G}$ agree with $\vec g_\alpha\in\mathbb R^{1\times\vert\Omega\vert}$ for $\alpha=0,\ldots,D-1$. This time, we construct two (not necessarily overlapping) linear families $\mathcal L\equiv\mathcal L(\mathbf G;\boldsymbol\mu)$ and $\widetilde{\mathcal L}\equiv\mathcal L(\widetilde{\mathbf G};\tilde{\boldsymbol\mu})$ on the simplex, defined by generically distinct vectors of generalized moments $\boldsymbol\mu\in\mathbb R^D$ and $\tilde{\boldsymbol\mu}\in\mathbb R^{\tilde D}$. Then, 
    \begin{equation}
    \tilde{\prob q} = 
    \argmin_{\prob p\in\widetilde{\mathcal L}}\infdiv{\prob p}{\prob v}
    = 
    \argmin_{\prob p\in\widetilde{\mathcal L}}\infdiv{\prob p}{\prob q} 
    \quad\text{where}\quad \prob q = \argmin_{\prob p\in\mathcal L} \infdiv{\prob p}{\prob v}
   ~.
\end{equation}
If $\tilde D=D$, the $I$-projection of $\prob v$ slides over the intersection of the same exponential family Eq.~\eqref{eq:ExponentialFamily} with different linear families, where $\mathcal L(\mathbf G,\tilde{\boldsymbol\mu})$ is said to be the {translate} of $\mathcal L(\mathbf G;\boldsymbol\mu)$.

To verify the more general claim, we need to show that the  $I$-projections of $\prob q$ and $\prob v$ onto the same $\widetilde{\mathcal L}$ coincide. First, we write the exponential form of the  $I$-projection of $\prob q$ onto $\widetilde{\mathcal L}$,
\begin{equation*}
    \tilde q(x) = q(x) \exp\left\{ \sum_{\alpha=0}^{\tilde D-1}\lambda_\alpha(\tilde{\boldsymbol\mu}) g_\alpha(x) \right\} ~,
\end{equation*}
making the dependence of Lagrange multipliers on generalized moments explicit to readily relate to the corresponding minimization problem that fixes $\boldsymbol\lambda\in\mathbb R^{\tilde D}$ at optimality.
Substituting the exponential form for $\prob q$ being the $I$-projection of $\prob v$ onto $\mathcal L$,
\begin{equation*}
    \tilde q(x) = v(x) \exp\left\{ \sum_{\alpha=0}^{D-1}\left[\theta_\alpha(\boldsymbol\mu) + \lambda_\alpha(\tilde{\boldsymbol\mu})\right] g_\alpha(x) + \sum_{i=D}^{\tilde D-1} \lambda_i(\tilde{\boldsymbol\mu})  g_i(x)  \right\} \equiv v(x) \exp\left\{ \sum_{\alpha=0}^{\tilde D-1} \tilde\theta_\alpha g_\alpha(x) \right\}~,
\end{equation*}
we recognize that the $\tilde{\prob q}$ assumes the exponential form of the $I$-projection of $\prob v$ onto $\tilde{\mathcal L}$ upon redefining $\tilde\theta_\alpha\equiv \theta_\alpha + \lambda_\alpha$ for $\alpha=0,\ldots, D-1$ and $\tilde\theta_i\equiv \lambda_i$ for the remaining $i=D,\ldots, \tilde D-1$. 
Since $\tilde{\prob q}$ satisfies by construction the conditions on generalized moments defining $\widetilde{\mathcal L}$, and has the form of a member of the exponential family generated by $\tilde{\mathbf G}$ around $\prob v$, it must also be the unique $I$-projection of $\prob v$ onto $\widetilde{\mathcal L}$.

Similar to the previous nested-projection result in Eq.~\eqref{eq:app:nested_iProjections}, the practical relevance of this generalized property becomes apparent whenever the direct $I$-projection of $\prob v$ onto $\mathcal L(\tilde{\mathbf G};\tilde{\boldsymbol\mu})$ is numerically unstable due to $\infdiv{\tilde{\prob q}}{\prob v}\gg1$. If the $I$-projection $\prob q$ is numerically accessible, one may instead obtain $\tilde{\prob q}$ by exploiting a hierarchy 
$\infdiv{\tilde{\prob q}}{\prob q}\leq\infdiv{\tilde{\prob q}}{\prob v}$ that can be arranged for small deviations $\tilde{\mu}_\alpha-\mu_\alpha$ . 
The result is also useful for constructing families of $I$-projections associated with a fixed $\mathbf G$ by varying admissible values of $N\boldsymbol\mu\in\mathbf G \mathbb N_0^{\vert\Omega\vert}$, typically in a neighborhood of the observed estimates.

\section{Some elementary calculus}\label{app:calculus}

For any real symmetric positive-definite matrix $\mathbf A\in\mathbb R^{D\times D}$, the associated quadratic form 
\begin{equation}
\label{eq:app:QuadraticForm}
Q_{\mathbf A}(\mathbf x) = \mathbf x^T \mathbf A \mathbf x = \sum_{i,j=1}^D x_iA_{ij}x_j
\end{equation}
remains strictly  positive for any $\mathbf x\in\mathbb R^D\setminus\lbrace0\rbrace$.
The core formula needed in effective descriptions is the multidimensional Gaussian integral with source term $\mathbf j\in\mathbb R^{D}$
\begin{equation}
\label{eq:GaussianIntegrals:Sourced}
\int d^D \mathbf x \exp \left\{ -\tfrac12Q_{\mathbf A}(\mathbf x) - \mathbf x^T \mathbf j  \right\} 
=
\sqrt{\frac{(2\pi)^D}{\det \mathbf A}} \exp\left\{ \tfrac12 \mathbf j^T \mathbf A^{-1} \mathbf j \right\} ~.
\end{equation}
In the context of multivariate statistics, $\mathbf A$ is called the precision matrix and its inverse the covariance matrix $\mathbf\Sigma\equiv\mathbf A^{-1}$. In this parametrization,  the associated normal distribution with probability density function
\begin{equation}
\label{app:MultivariateNormalDensity}
\sqrt{\det\left(\frac{\mathbf A}{2\pi}\right)}
\exp\left\{ - \tfrac12 \left(\mathbf x - \mathbf A^{-1}\mathbf j  \right)^T  \mathbf A \left(\mathbf x - \mathbf A^{-1}\mathbf j  \right) \right\}
\end{equation}
is centered around mean $\mathbf A^{-1}\mathbf j$.
The expectation of the quadratic form under the associated sourceless normal distribution,
\begin{equation}
\label{eq:appGaussian:QuadraticExpectation}
\sqrt{\det\left(\frac{\mathbf A}{2\pi}\right)}\int \text{d}^D\mathbf x\,\exp\left\{-\tfrac12 Q_{\mathbf A}(\mathbf x)\right\}Q_{\mathbf A}(\mathbf x)
= D
~,
\end{equation}
evaluates to $D$, the dimension of the space.

Any marginal of the normal multivariate distribution is again a normal multivariate distribution in the complementary variables with associated quadratic form given by minimizing the original quadratic form in the variables being integrated out. 
To easily see this,  we start for clarity from the probability density function for $\mathbf j=\mathbf 0$, after appropriately shifting the coordinate system. 
To integrate over only some of the directions in $\mathbb R^D$,  
we need to partition the coordinates of the  ambient space along  the subspaces
\begin{equation}
\mathbf x = \begin{pmatrix}
\mathbf x_1\\
\mathbf x_2
\end{pmatrix}
~,
\end{equation}
where $\mathbf x_a\in\mathbb R^{d_a}$  for $a=1,2$. 
This breaks the precision matrix into  so-called conformable blocks
\begin{equation}
\mathbf A =
\begin{pmatrix}
\mathbf A_1 &  \mathbf B \\
\mathbf B^T & \mathbf A_2
\end{pmatrix}
~,
\end{equation}
where $\mathbf B$ is a $d_1\times d_2$ real matrix. $\mathbf A_a\in\mathbb R^{d_a\times d_a}$ are taken to be symmetric positive-definite matrices for $a=1,2$ and hence invertible. 
Accordingly, the quadratic form \eqref{eq:app:QuadraticForm} decomposes into 
\begin{align}
\label{eq:app:QuadraticForm_Decomposition}
Q_{\mathbf A}(\mathbf x) = \mathbf x^T_1 \mathbf A_1 \mathbf x_1 + \mathbf x^T_2 \mathbf A_2 \mathbf x_2 + 2 \mathbf x^T_1 \mathbf B \mathbf x_2 \equiv Q_{\mathbf A}(\mathbf x_1;\mathbf x_2)~.
\end{align}

Using Eq.~\eqref{eq:GaussianIntegrals:Sourced}, we can marginalize over the Gaussian integral in say $\mathbf x_1$ to obtain 
\begin{equation}
\label{eq:app:Marginal_over_x1}
\sqrt{\det\left(\frac{\mathbf A}{2\pi}\right)}\int \text{d}^{d_1}\mathbf x_1\,\exp\left\{-\frac{Q_{\mathbf A}(\mathbf x)}{2}\right\} =
\sqrt{\left(2\pi\right)^{-d_2}\frac{\det\mathbf A}{\det\mathbf A_1}}
\exp\left\{ - \tfrac12 \mathbf x_2^T\mathbf A/\mathbf A_1\mathbf x_2  \right\} ~,
\end{equation}
where 
\begin{equation}
\label{eq:Gaussian:SchursComplement}
    \mathbf A/\mathbf A_1 \equiv \mathbf A_2 - \mathbf B^T\mathbf A_1^{-1}\mathbf B
\end{equation}
is called the {Schur complement} of block $\mathbf A_1$.
That this resulting precision matrix is positive definite can be easily seen by completing the square for $\mathbf x_1$ in \eqref{eq:app:QuadraticForm_Decomposition}:
\begin{equation}
Q_{\mathbf A}(\mathbf x_1;\mathbf x_2) = \left(\mathbf x_1 + \mathbf A_1^{-1}\mathbf B\mathbf x_2\right)^T\mathbf A_1 \left(\mathbf x_1 + \mathbf A_1^{-1}\mathbf B\mathbf x_2\right)
+ 
\mathbf x_2^T \left(\mathbf A_2 - \mathbf B^T\mathbf A_1^{-1}\mathbf B \right)\mathbf x_2
~.
\end{equation}
Since the original quadratic form is strictly positive for all nonzero vectors in $\mathbb R^{D}$, it remains in particular positive when evaluated at 
\begin{equation}
\mathbf x_1^* = -\mathbf A_1^{-1}\mathbf B\mathbf x_2 = \argmin_{\mathbf y_1\in\mathbb R^{d_1}} Q_{\mathbf A}(\mathbf y_1;\mathbf x_2)~,
\end{equation}
for arbitrary fixed $\mathbf x_2\in\mathbb R^{d_2}\setminus\{\mathbf 0\}$. 
We emphasize that $\mathbf x_1^*$ is the minimizer of the original quadratic form over the $\mathbb R^{d_1}$ subspace, which in Section~\ref{ssc:NaiveGaussianization} amounts to the conditional minimization of \textsc{kl} divergence. 

In total, we recognize that
\begin{equation}
\forall~ \mathbf x_2 \in \mathbb R^{d_2} \setminus\lbrace\mathbf 0\rbrace\,:\quad
0 < Q_{\mathbf A}(\mathbf x_1^*;\mathbf x_2) = \mathbf x_2^T \left(\mathbf A_2 - \mathbf B^T\mathbf A_1^{-1}\mathbf B \right)\mathbf x_2~,
\end{equation}
meaning that $\mathbf A/\mathbf A_1$ has to be positive definite. 
By the calculus rules of determinants, the relation  
\begin{equation}
\det\mathbf A = \det\mathbf A_1 \cdot \det \mathbf A/\mathbf A_1 
\end{equation}
holds for Schur's complement.
This means that the normalization of the marginal integral in \eqref{eq:app:Marginal_over_x1} neatly reduces to 
\begin{equation}
\label{eq:Schur:det_relation}
\left(2\pi\right)^{-d_2}\frac{\det\mathbf A}{\det\mathbf A_1}
= \det\left(\frac{\mathbf A}{2\pi}\right)
\end{equation}
using that $D=d_1+d_2$ for complementary spaces in $\mathbb R^D$. Consequently, 
\begin{equation}
\sqrt{\det\left(\frac{\mathbf A}{2\pi}\right)}\int \text{d}^{d_1}\mathbf x_1\,\exp\left\{-\frac{Q_{\mathbf A}(\mathbf x)}{2}\right\} 
=
\sqrt{\det\left(\frac{\mathbf A/\mathbf A_1}{2\pi}\right)}\exp\left\{-\frac{Q_{\mathbf A/\mathbf A_1}(\mathbf x_2)}{2}\right\} 
~,
\end{equation}
showing the closure for the class of multivariate normal distribution under marginalization.

\end{appendix}

\end{document}